\documentclass[10pt,letterpaper]{amsart}
\usepackage{amsfonts}
\usepackage{amsmath}
\usepackage{amsrefs}
\usepackage{amssymb}
\usepackage{amsthm}
\usepackage{array}
\usepackage{bm}
\usepackage{dcpic,pictexwd}
\usepackage{fancyhdr}
\usepackage{indentfirst}
\usepackage{ifthen}
\usepackage{accents}
\usepackage{marginnote}
\usepackage{fixmath}

\usepackage{graphicx}
\newcommand{\hOmega}{\accentset{\LARGEhat}{\vphantom{\scalebox{1.2}{V}}\Omega}}
\newcommand\hV{{\accentset{\LARGEhat}{\vphantom{\scalebox{1.2}{V}}V}}}
\newcommand\hW{{\accentset{\LARGEhat}{\vphantom{\scalebox{1.2}{W}}W}}}
\newcommand\hZ{{\accentset{\LARGEhat}{\vphantom{\scalebox{1.2}{Z}}Z}}}
\newcommand\hCalV{{\accentset{\LARGEhat}{\vphantom{\scalebox{1.2}{V}}{\mathcal V}}}}

\newcommand\hCalZ{{\accentset{\LARGEhat}{\vphantom{\scalebox{1.2}{Z}}{\mathcal Z}}}}
\newcommand\cV{{\accentset{\LARGEcheck}{\vphantom{\scalebox{1.2}{V}}V}}}
\newcommand\cW{{\accentset{\LARGEcheck}{\vphantom{\scalebox{1.2}{W}}W}}}
\newcommand\cZ{{\accentset{\LARGEcheck}{\vphantom{\scalebox{1.2}{Z}}Z}}}
\newcommand\cCalV{{\accentset{\LARGEcheck}{\vphantom{\scalebox{1.2}{V}}{\mathcal V}}}}

\newcommand\cCalZ{{\accentset{\LARGEcheck}{\vphantom{\scalebox{1.2}{Z}}{\mathcal Z}}}}
\newcommand\cOmega{\accentset{\LARGEcheck}{\vphantom{\scalebox{1.2}{V}}\Omega}}
\newcommand\hL{{\accentset{\LARGEhat}{\vphantom{\scalebox{1.2}{L}}L}}}
\newcommand\hS{{\accentset{\LARGEhat}{\vphantom{\scalebox{1.2}{S}}S}}}
\newcommand\cL{{\accentset{\LARGEcheck}{\vphantom{\scalebox{1.2}{L}}L}}}
\newcommand\cS{{\accentset{\LARGEcheck}{\vphantom{\scalebox{1.2}{S}}S}}}
\newcommand\hJ{{\accentset{\LARGEhat}{\vphantom{\scalebox{1.2}{J}}J}}}
\newcommand\cJ{{\accentset{\LARGEcheck}{\vphantom{\scalebox{1.2}{J}}J}}}

\newcommand{\htau}{\accentset{\Largehat}{\vphantom{\scalebox{1.2}{s}}\tau}}
\newcommand\ctau{\accentset{\Largecheck}{\vphantom{\scalebox{1.2}{s}}\tau}}
\newcommand\hsigma{\accentset{\Largehat}{\vphantom{\scalebox{1.2}{s}}\sigma}}
\newcommand\csigma{\accentset{\Largecheck}{\vphantom{\scalebox{1.2}{s}}\sigma}}
\newcommand\homega{\accentset{\Largehat}{\vphantom{\scalebox{1.2}{s}}\omega}}
\newcommand\comega{\accentset{\Largecheck}{\vphantom{\scalebox{1.2}{s}}\omega}}

\numberwithin{equation}{section}
\newtheorem{Theorem}{Theorem}[section]

\newtheorem*{ThmA}{Theorem A}
\newtheorem*{ThmB}{Theorem B}
\newtheorem*{ThmC}{Theorem C}
\newtheorem*{ThmD}{Theorem D}
\newtheorem*{ThmE}{Theorem E}

\newtheorem{Lemma}[Theorem]{Lemma}

\newtheorem{Corollary}[Theorem]{Corollary}
\newtheorem{Definition}[Theorem]{Definition}
\theoremstyle{definition}

\newtheorem{Remark}[Theorem]{Remark}

\newfont{\deffont}{cmbxti10}
\newfont{\german}{eufm10}
\newfont{\mymath}{cmr12}

\newcommand\lieg{\mathfrak g}

\newcommand\lies{{\german s}}
\newcommand\lieh{{\mathfrak h}}
\newcommand\lieo{{\mathfrak o}}

\renewcommand\qed{\hfill\hbox{\vrule width 4pt height 6pt depth 1.5 pt}}

\newcommand\hook{\mathbin{\raise2.5pt\hbox{\hbox{{\vbox{\hrule height.4pt width6pt depth0pt}}}\vrule height3pt width.4pt depth0pt}\,}}
\newcommand\extd{d\mspace{2mu}}
\newcommand\cTM{T^*\kern-2ptM}

\newcommand\Ad{\text{\rm Ad}}

\newcommand\semibasic{\text{\bf  sb}}

\newcommand\vess{\text{\german v\german e\german s\german s}}

\newcommand\thetaX{\theta_{\kern -1 pt X}}

\newcommand\ann{\text{\rm ann}}
\newcommand\CalPf{\mathcal{P}\text{\it \kern -.3pt f}}
\newcommand\Real{\text{\bf  R}}

\newcommand\diag{\text{\rm diag}}
\renewcommand\:{\colon}

\newcommand\TM{T\kern -2pt M}
\newcommand\real{\text{\bf R}}
\newcommand\mycap{\hbox{\ $\rlap{\kern -.3pt $\cap$}\raise.8pt\hbox{$\scriptstyle+$}$\ } }

\newcommand\Ao{{\kern-2.3pt}\stackrel{\scriptscriptstyle o}{A}{}{\kern-2.3pt}}
\newcommand\Bo{{\kern-2.3pt}\stackrel{\scriptscriptstyle o}{B}{}{\kern-2.3pt}}
\newcommand\Ko{{\kern-2.3pt}\stackrel{\scriptscriptstyle o}{K}{}{\kern-2.3pt}}
\newcommand\Co{{\kern-2.3pt}\stackrel{\scriptscriptstyle o}{C}{}{\kern-2.3pt}}
\newcommand\Qo{{\kern-2.3pt}\stackrel{\scriptscriptstyle o}{Q}{}{\kern-2.3pt}}
\newcommand\Mo{{\kern-2.3pt}\stackrel{\scriptscriptstyle o}{M}{}{\kern-2.3pt}}

\newcommand\Xo{{\stackrel{\scriptscriptstyle o}{X}}{\kern-1.3pt}}
\newcommand\Yo{{\stackrel{\scriptscriptstyle o}{Y}}{\kern-1.3pt}}

\DeclareMathOperator{\rank}{rank}

\DeclareMathOperator{\spn}{span}
\newboolean{proofmode}
\newcommand{\StTag}[1]{ \label{st:#1}
\ifthenelse{\boolean{proofmode}}{\ \marginpar{\quad\scriptsize st:#1} }{}      }
\newcommand{\EqTag}[1]{
\ifthenelse{\boolean{proofmode}}
{ {\label{eq:#1}}
  \stepcounter{equation}
  \tag{\theequation \rlap{\kern 23 pt{\scriptsize eq:#1}}}
}
{\label{eq:#1}}
 }
\newcommand{\EqRef}[1]{\eqref{eq:#1}}
\newcommand{\StRef}[1]{\ref{st:#1}}
\newcolumntype{C}{>\scriptstyle>{$}c <{$} }
\newcolumntype{L}{>\scriptstyle >{$} l <{$} }
\newcommand\CalA{\mathcal{A}}
\newcommand\CalB{\mathcal{B}}
\newcommand\CalC{\mathcal{C}}
\newcommand\CalE{\mathcal{E}}
\newcommand\CalI{\mathcal{I}}
\newcommand\CalL{\mathcal{L}}
\newcommand\CalJ{\mathcal{J}}
\newcommand\CalK{\mathcal{K}}

\newcommand\CalF{\mathcal{F}}
\newcommand\CalO{\mathcal{O}}

\newcommand\CalS{\mathcal{S}}

\newcommand\sd{\mathbin{ \raise0.0pt\hbox{ \vrule height5pt width.4pt depth0pt}\!\times}}
\newcommand\barC{
	\hbox{\kern 2.3 true pt
	\vbox{\hrule width 6.5  true pt height .3 true pt \kern .9 true pt
	\hbox{\kern -0.8 true pt $C$}}}}
\newcommand\barM{
	\hbox{\kern 2.3 true pt
	\vbox{\hrule width 8.5  true pt height .3 true pt \kern .7 true pt
	\hbox{\kern -2.3 true pt $M$}}}}

\newcommand\barK{
	\hbox{\kern 2.5 true pt
	\vbox{\hrule width 6  true pt height .3 true pt \kern 1.4 true pt
	\hbox{\kern -2 true pt $K$}}}}

\newcommand\barW{
	\hbox{\kern 2.5 true pt
	\vbox{\hrule width 6  true pt height .3 true pt \kern 1.4 true pt
	\hbox{\kern -2 true pt $W$}}}}

\newcommand\barU{
	\hbox{\kern .8 true pt
	\vbox{\hrule width 6.5  true pt height .3 true pt \kern .9 true pt
	\hbox{\kern -.8 true pt $U$}}}}
\newcommand\barCalI{
	\hbox{\kern 4.3 true pt
	\vbox{\hrule width 6.5  true pt height .3 true pt \kern .9 true pt
	\hbox{ \kern -4.3 true pt  $\CalI$}}}}
\newcommand\barXi{
	\hbox{\kern 1 true pt
	\vbox{\hrule width 6.5  true pt height .3 true pt \kern .9 true pt
	\hbox{\kern 1 true pt $\Xi$}}}}

\newcommand\Largehat{\smash{\raise -7.5 pt \hbox{\rm\Large\^{}}}}
\newcommand\LARGEhat{\smash{\raise -9.5 pt \hbox{\rm\LARGE\^{}}}}
\newcommand\hugehat{\smash{\raise -7.5 pt \hbox{\rm\huge\^{}}}}
\newcommand\Hugehat{\smash{\raise -7.5 pt \hbox{\rm\Huge\^{}}}}
\newcommand\Largecheck{\smash{\raise -7.5 pt \hbox{\rm\Large\v{}}}}
\newcommand\LARGEcheck{\smash{\raise -9 pt \hbox{\rm\LARGE\v{}}}}
\newcommand\hugecheck{\smash{\raise -7.5 pt \hbox{\rm\huge\v{}}}}
\newcommand\Hugecheck{\smash{\raise -7.5 pt \hbox{\rm\Huge\v{}}}}

\newcommand\hU{{\accentset{\LARGEhat}{U}}}

\newcommand\cU{{\accentset{\LARGEcheck}{U}}}

\newcommand\cSigma{\accentset{\LARGEcheck}{\Sigma}}
\newcommand\hSigma{\accentset{\LARGEhat}{\Sigma}}

\newcommand\bfa{\boldsymbol{a}}
\newcommand\bfA{\mathbold{A}}
\newcommand\bfb{\boldsymbol{b}}
\newcommand\bfB{\boldsymbol{B}}
\newcommand\bfc{\boldsymbol{c}}
\newcommand\bfC{\boldsymbol{C}}
\newcommand\bfE{\boldsymbol{E}}
\newcommand\bfF{\boldsymbol{F}}
\newcommand\bfG{\boldsymbol{G}}
\newcommand\bfH{\boldsymbol{H}}
\newcommand\bfI{\boldsymbol{I}}
\newcommand\bfJ{\boldsymbol{J}}

\newcommand\bfM{\boldsymbol{M}}
\newcommand\bfN{\boldsymbol{N}}
\newcommand\bfq{\mathbf{q}}
\newcommand\bfp{\mathbf{p}}
\newcommand\bfr{\mathbf{r}}
\newcommand\bfR{\boldsymbol{R}}
\newcommand\bfS{\boldsymbol{S}}
\newcommand\bfT{\boldsymbol{T}}
\newcommand\bfX{\boldsymbol{X}}
\newcommand\bfY{\boldsymbol{Y}}
\newcommand\bfP{\boldsymbol{P}}

\newcommand\bfGamma{\boldsymbol{\Gamma}}

\newcommand\bfmu{{\boldsymbol{\mu}}}
\newcommand\bfalpha{{\boldsymbol{\alpha}}}
\newcommand\bfbeta{{\boldsymbol{\beta}}}

\newcommand\bfpartial{{\boldsymbol{\partial}}}
\newcommand\bfgamma{{\boldsymbol{\gamma}}}
\newcommand\bftheta{{\boldsymbol{\theta}}}
\newcommand\bftau{{\boldsymbol{\tau}}}
\newcommand\bfnu{{\boldsymbol{\nu}}}
\newcommand\bfxi{{\boldsymbol{\xi}}}
\newcommand\bfeta{{\boldsymbol{\eta}}}
\newcommand\bfsigma{{\boldsymbol{\sigma}}}
\newcommand\bflambda{{\boldsymbol{\lambda}}}
\newcommand\bfpi{{\boldsymbol{\pi}}}

\newcommand\bfthetaX{{\boldsymbol{\theta_{\kern -1 pt X}}}}
\newcommand\bfthetaY{{\boldsymbol{\theta_Y}}}
\newcommand\bftildeA{\boldsymbol{{\tilde A}}}
\newcommand\bftildeB{\boldsymbol{{\tilde B}}}
\newcommand\bftildeE{\boldsymbol{{\tilde E}}}
\newcommand\bftildeF{\boldsymbol{{\tilde F}}}
\newcommand\bftildeG{\boldsymbol{{\tilde G}}}
\newcommand\bftildeH{\boldsymbol{{\tilde H}}}

\newcommand\bfhA{\boldsymbol{{\hat A}}}
\newcommand\bfhG{\boldsymbol{{\hat G}}}

\newcommand\bfhF{\boldsymbol{{\hat F}}}

\newcommand\bfheta{\boldsymbol{{\hat \eta}}}

\newcommand\bfhpi{\boldsymbol{{\hat \pi}}}
\newcommand\bfhsigma{\boldsymbol{{\hat \sigma}}}
\newcommand\bfhomega{\boldsymbol{{\hat \omega}}}

\newcommand\bfhtau{\boldsymbol{{\hat \tau}}}

\newcommand\bfcF{\boldsymbol{{\check F}}}
\newcommand\bfcH{\boldsymbol{{\check H}}}

\newcommand\bfceta{\boldsymbol{{\check \eta}}}

\newcommand\bfctau{\boldsymbol{{\check\tau}}}

\newcommand\bfcpi{\boldsymbol{{\check\pi}}}
\newcommand\bfcsigma{\boldsymbol{{\check\sigma}}}
\newcommand\bfcomega{\boldsymbol{{\check\omega}}}

\newcommand\Rtheta{{ \raise 1pt \hbox{$\scriptstyle {\boldsymbol{\theta}}$}}}
\newcommand\Ltheta{{ \lower 1pt \hbox{$\scriptstyle {\boldsymbol{\theta}}$}}}

\newcommand\Rsigma{{ \raise 1pt \hbox{$\scriptstyle \sigma$}}}

\newcommand\Reta{{ \raise 1pt \hbox{$\scriptstyle \eta$}}}

\newcommand\vecthsigma{\partial_{{\displaystyle \hat {\raise 1.3pt \hbox{$\scriptstyle \sigma$}}}^a}}
\newcommand\vectcsigma{\partial_{{\displaystyle \check {\raise 1.3pt \hbox{$\scriptstyle \sigma$}}}^\alpha}}

\setboolean{proofmode}{false}

\begin{document}
\title{B\"acklund Transformations for Darboux Integrable Differential Systems} 
\author{I.M. Anderson}
\address{Department of Mathematics and Statistics, Utah State University, Logan Utah, 84322}
\email{Ian.Anderson@usu.edu}

\author{M.E. Fels}
\address{Department of Mathematics and Statistics, Utah State University, Logan Utah, 84322}
\email{Mark.Fels@usu.edu}

\begin{abstract}
The purpose of this paper is to present a novel group-theoretical method for constructing B\"acklund transformations between
systems of differential equations. Our approach is based upon the definition  of B\"acklund transformations as integrable extensions of  exterior differential systems. The construction of these transformations is obtained using the general theory of symmetry reduction of  differential systems. Our method is then applied to differential systems which are integrable by the  method of Darboux and a detailed  understanding of the 
B\"acklund transformations  for these systems is obtained.
\end{abstract}

\maketitle



\bigskip

\setcounter{page}{1}
{\small
\noindent
\begin{center}
{\sc Table of Contents}
\end{center}
\par\noindent
1. Introduction
\par\noindent
2. Preliminaries
\par\noindent
\quad 2.1. Integrable Extensions of Differential Systems
\par\noindent
\quad  2.2. Reduction of Differential Systems
\par\noindent
\quad  2.3. Symmetry Reduction of Pfaffian Systems by Free Group Actions
\par\noindent
3. B\"acklund Transformations by Symmetry Reduction of  EDS
\par\noindent
\quad 3.1. Reductions of EDS and Commutative Diagrams
\par\noindent
\quad 3.2.  Subgroup Reductions and Integrable Extensions
\par\noindent
4. Darboux Integrable Differential Systems
\par\noindent
5. Integrable Extensions of Darboux Integrable Differential Systems
\par\noindent
\quad 5.1.  Integrable Extensions are Darboux Integrable
\par\noindent
\quad 5.2.  Maximally Compatible Integrable Extensions for Darboux Integrable Differential Systems
\par\noindent
6. Group Theoretic Constructions of  Darboux Integrable Systems
\par\noindent
\quad 6.1. Darboux Integrable by Symmetry Reduction of Sums of EDS
\par\noindent
\quad 6.2. Reduction by Diagonal Actions
\par\noindent
\quad 6.3. Reduction by non-Diagonal Actions
\par\noindent
7.  Vessiot Algebras as Fundamental Invariants for Darboux Integrable Systems
\par\noindent
\quad 7.1. The Vessiot Algebra of a Darboux Integrable Differential System
\par\noindent
\quad 7.2. The Prolongation of Darboux Integrable Systems
\par\noindent
\quad 7.3. 4-Adapted Coframes and the Vessiot Algebras for the Quotient Construction
\par\noindent
8. Quotient Representations for Maximally Compatible Integrable Extensions
\par\noindent
\quad 8.1.Uniqueness of Maximally Compatible Integrable Extensions
\par\noindent
9.  B\"acklund Transformations for Darboux  Integrable Systems
\par\noindent
Appendices
\par\noindent
References
\par\noindent
}

\section{Introduction}
	Let $\Delta_1(z) =0$   and $\Delta_2(u) = 0$  be two systems of  partial differential equations in the  unknown  functions
	$z$ and  $u$.  A B\"acklund transformation between these two systems of equations  is a third system of equations 
	$\Delta_3(u, z) = 0$ with the  following property. If  a solution for $\Delta_1(z)$ =0 is given, then
	$\Delta_3(u, z) = 0$ becomes  a system of total differential equations for $u$  which is  then a solution to  $\Delta_2(u) = 0$. 
	Likewise, if a solution to  $\Delta_2(u)$ =0 is given, then 
	$\Delta_3(u, z) = 0$ becomes  a system of total differential equations for $z$  which is  then  a solution to  $\Delta_1(z) = 0$. 
	One  of the best known  examples of a B\"acklund transformation, which  relates  solutions of the wave equation  to solution of
	Liouville's equation
\begin{equation*}
	\Delta_1(z) = z_{xy} = 0  \quad\text{and}\quad \Delta_2(u) = u_{xy} = e^u 
\end{equation*}
	is 
\begin{equation}
	z_x-u_x=  \sqrt{2} \exp {\dfrac{z+u}{2}} \quad z_y+u_y = -\sqrt{2}\exp {\dfrac{u-z}{2}} 
\EqTag{BT1}
\end{equation}
	While B\"acklund transformations have a long and  distinguished history in  the study of integrable differential equations
	\cite{goursat:1925a}, \cite{rogers-shadwick:1982}, it 
	nevertheless remains a challenging problem to give systematic methods for their construction.
	
The purpose of this paper is  two-fold:

\smallskip
\noindent
{\bf [i]} to show that B\"acklund transformations can easily be constructed using the general theory of 
	symmetry reduction of exterior differential systems; and

\smallskip
\noindent
{\bf[ii]} 
	to prove that  this symmetry approach to B\"acklund transformations  
	leads to a very general, yet precise, understanding of B\"acklund transformations
	for Darboux integrable differential systems.

\smallskip
\noindent
	
	Using these group-theoretic constructions, it is possible 

\smallskip
\noindent
{\bf[i]}  to systematically re-construct many known examples of  B\"acklund transformations; 

\smallskip
\noindent
{\bf[ii]}  to construct new examples of  B\"acklund transformations  which are beyond the capabilities of existing methods; and 

\smallskip
\noindent
{\bf[iii]} to establish  simple, but novel and non-trivial,   geometric obstructions for the existence of B\"acklund transformations
	    for certain classes of differential equations.

\smallskip
\noindent
	These, and other applications,  will appear in a sequel to this article.

	Our work is formulated within the differential-geometric setting of  symmetry  reduction of  exterior differential systems (EDS).
	From this geometric viewpoint,  
	two  differential systems $\CalI_1 \subset \Omega^*(M_1) $ and $\CalI_2  \subset \Omega^*(M_2) $, 
	defined on manifolds $M_1$ and $M_2$, are said to be related by a B\"acklund transformation
	if there exists a  differential system $\CalB  \subset \Omega^*(N) $ on a manifold $N$  and maps
\begin{equation}
\begin{gathered}
\begindc{\commdiag}[3]
\obj(0, 14)[B]{$(\CalB, N)$}
\obj(-12, 0)[I1]{$(\CalI_1, M_1)$}
\obj(12, 0)[I2]{$(\CalI_2, M_2)$\,}
\mor{B}{I1}{$\bfp_1$}[\atright, \solidarrow]
\mor{B}{I2}{$\bfp_2$}[\atleft, \solidarrow]
\enddc
\end{gathered}
\EqTag{BackDef}
\end{equation}
	which  define  $\CalB$  as  integrable extensions for both  $\CalI_1$ and $\CalI_2$. For the explicit  re-formulation 
	of the B\"acklund transformation  \EqRef{BT1} into this EDS setting, we refer the reader to \cite{anderson-fels:2012a}.
	The fact that  $\CalB$ is 
	an integrable extension of $\CalI_1$  means, in essence,  that if $P \subset M_1$ is an integral manifold
	of  $\CalI_1$ then  the restriction of $\CalB$  to $\bfp_1^{-1}(P)$ is a Pfaffian system which is completely integrable
	in the sense of Frobenius, and so $\bfp_1^{-1}(P)$ is foliated by integral manifolds of $\CalB$ \cite{bryant-griffiths:1995a}. 
	In this way integral manifolds of $\CalI_1$ 
	can  be locally  lifted,  by solving  ordinary differential equations, to integral manifolds of $\CalB$ which 
	then project by $\bfp_2$ to integral manifolds of $\CalI_2$ and conversely.  
	Although  the first impression of  this EDS formulation of a B\"acklund transformation  may be that it is rather abstract 
	and difficult  to work with, we shall  see that  it is actually ideally suited for a general  analysis of B\"acklund transformations.

\medskip
	This paper  consists  of 3  parts.  In the first part, which consists of Sections 2 and 3,  
	we give a  very simple, group-theoretic method for 
	constructing B\"acklund transformations. 
	To describe our method in some detail, let us
	first recall a few definitions which we shall  use many times.  Let $G$ be a Lie group which acts on  a 
	manifold $M$.  We say that $G$ acts {\deffont regularly} on $M$ whenever the orbit space $M/G$ admits  
	the structure of a differentiable manifold for which the canonical projection map $\bfq_G : M \to M/G$ is a surjective submersion. 
	We say that $G$ is a  {\deffont symmetry group} of a differential system $\CalI \subset \Omega^*(M)$ 
	if the diffeomorphisms of $M$ defined by $G$ 
	preserve  $\CalI$, and we say the $G$ acts {\deffont transversely} to $\CalI$ if the tangent spaces to the orbits of
	$G$ are pointwise transverse to the space of 1-forms  in $\CalI$ (that is, no non-zero vector in the tangent space to the orbits is annihilated by all the 1-forms in $\CalI$, see equation \EqRef{Itrans}).   As we  shall see below in Theorem A, this transversality 
	condition plays an important role in the construction of integrable extensions.

	Let   $\bfp: M \to N$ be a smooth submersion. Then we define the  {\deffont reduced differential system} $ \CalI/\bfp$ on $N$ by
\begin{equation}
	 \CalI/\bfp = \{\, \theta \in \Omega^*(N) \, |\, \bfp^*(\theta) \in \CalI \, \} .
	\EqTag{Intro1}
\end{equation}
	In the special case where $G$ is a regularly acting symmetry group of $\CalI$,
	we shall write $\CalI/G$ in place of $\CalI/\bfq_G$.  Finally,  we say that the diagram 
\begin{equation*}
\begin{gathered}
\begindc{\commdiag}[3]
\obj(-25, 0)[I]{$(\CalI, M)$}
\obj(0, 0)[H]{$(\CalJ, N)$}
\obj(0, -17)[I1]{$(\CalK, L)$}
\mor{I}{H}{$\bfp$}[\atleft, \solidarrow]
\mor{H}{I1}{$\bfq$}[\atleft, \solidarrow]
\mor{I}{I1}{$\bfr $}[\atright, \solidarrow]
\enddc
\end{gathered}
\end{equation*}
	is a {\deffont commutative diagram} of EDS if  $\bfq \circ \bfp = \bfr$,   $\CalJ= \CalI/\bfp$ and  $\CalK = \CalJ/\bfq =  \CalI/\bfr$.
	Additional information on these definitions can be
	found in Section 2 and also in \cite{anderson-fels:2005a}. 

	With these definitions in place we now state our  main theorem on the  construction of B\"acklund transformations by  symmetry 
	group reduction. The proof  is given in Section 3.

\begin{ThmA}
\StTag{ThIntro1}
	Let $\CalI$ be a differential system on $M$ with symmetry groups $G_1$ and $G_2$. 
	Let $H$ be a  common subgroup of  $G_1$ and $G_2$ and assume that the actions of $G_1$, $G_2$ and $H$ are all
	regular on $M$.  Then the orbit  projection maps $\bfp_1: M/H \to M_1/G_1$ and $\bfp_2: M/H \to M_2/G_2$ are smooth
	surjective submersions  and
\begin{equation}
\begin{gathered}
\begindc{\commdiag}[3]
\obj(0, 32)[I]{$(\CalI,M)$}
\obj(0, 12)[H]{$(\CalI /H,M/H)$}
\obj(-30, 0)[I1]{$(\CalI /G_1,M/G_1)$}
\obj(30,  0)[I2]{$(\CalI /G_2,M/G_2)$}
\mor{I}{H}{$\bfq_H$}[\atleft, \solidarrow]
\mor{H}{I1}{$\bfp_1$}[\atleft, \solidarrow]
\mor{H}{I2}{$\bfp_2$}[\atright, \solidarrow]
\mor{I}{I1}{$\bfq_{G_1}$}[\atright, \solidarrow]
\mor{I}{I2}{$\bfq_{G_2.}$}[\atleft, \solidarrow]
\enddc
\end{gathered}
\EqTag{Intro21}
\end{equation}
	is a commutative diagram of EDS. Furthermore, if  the actions of $G_1 $ and $G_2$ are transverse to $\CalI$, then 
	the  maps  in \EqRef{Intro21} are all  integrable extensions and  the diagram
\begin{equation}
\begin{gathered}
\begindc{\commdiag}[3]
\obj(0, 20)[H]{$(\CalI /H,M/H)$}
\obj(-20, 0)[I1]{$(\CalI /G_1,M/G_1)$}
\obj(20,  0)[I2]{$(\CalI /G_2,M/G_2)$}
\mor{H}{I1}{$\bfp_1$}[\atleft, \solidarrow]
\mor{H}{I2}{$\bfp_2$}[\atright, \solidarrow]
\enddc
\end{gathered}
\EqTag{Intro222}
\end{equation}
	defines  $(\CalI/H, M/H)$ as  a B\"acklund transformation between $(\CalI/G_1, M/G_1)$ and  $(\CalI/G_2, M/G_2)$.
\end{ThmA}
	The double fibration in \EqRef{Intro222} is,  at least in terms of the manifolds which are constructed, a simple generalization
	of the double fibration  construction  for homogeneous spaces presented  by
	Baston and Eastwood  \cite{baston-eastwood:1989a}  (page 69)  in the context of the Penrose transform. 

\medskip

	 There are many  known examples of  B\"acklund transformations, 
	such as the B\"acklund transformation  \EqRef{BT1},
	between systems of 
	differential equations which are integrable by the method of Darboux \cite{Clelland-Ivey:2009a}, \cite{zvyagin:1991a}.
	In the classical literature, which deals only with scalar  partial differential equations in the plane,
	an equation is called Darboux integrable if it admits a sufficient number 
	of  intermediate integrals.  The method of Darboux states that the
	general solution to the  given partial differential equation can be found by integrating  ODE systems constructed from these 
	intermediate integrals. See \cite{bryant-griffiths-hsu:1995a},  \cite{goursat:1897a}, \cite{goursat:1899a}, 
	\cite{ivy-langsberg:2003},  \cite{stormark:2000a}.  
	Vessiot  \cite {vessiot:1939a}, \cite{vessiot:1942a} made the remarkable discovery that these ODE are of Lie type and 
	thereby introduced Lie group-theoretic methods into the study of Darboux integrable equations.  In
	\cite{anderson-fels-vassiliou:2009a} a very general definition of 
	Darboux integrability for exterior differential systems, one which goes well beyond the classical definition,	  
	was introduced and Vessiot's results were re-interpreted and generalized in terms of symmetry reduction of
	differential systems (see also Appendix A in \cite{anderson-fels:2013a}).

\medskip 

 {\it The main  goal for the remainder of this article
	is to show that   B\"acklund transformations between Darboux integrable systems can be   constructed
	by  combining the general group-theoretic approach provided  of Theorem A  with  the realization of Darboux integrable systems by
	symmetry reduction  given in  \cite{anderson-fels-vassiliou:2009a}. In certain   cases, we  will be able to prove that 
	all  B\"acklund transformations between Darboux integrable systems are obtained using  the  group-theoretic method of Theorem A.}

\medskip 

	Sections 4 and 5, which constitute  the second part of this paper,  provide some foundational 
	material on Darboux integrable systems  and their integrable extensions. 
	The  general definition of Darboux integrability  is reviewed in Section 4.  A simplification 
	of the criteria for Darboux integrability is given and used to establish the following fundamental relationship 
	between integrable extensions and Darboux integrability (see Theorem \StRef{DIext} for an additional technical  hypothesis).

\begin{ThmB} 
	Let $\bfp : (\CalE, N) \to (\CalI, M)$ be an integrable extension of differential systems. If   
	$\CalI$ is Darboux integrable, then $\CalE$ is Darboux integrable. 
\end{ThmB}

	This theorem shows, in particular, that if either of the differential systems $\CalI_1$ or $\CalI_2$  in  \EqRef{BackDef} is 
	Darboux integrable, then  the B\"acklund  transformation $\CalB$ must be Darboux integrable.  This important observation 
	has  heretofore not appeared in the  literature; it alone greatly simplifies the study of  B\"acklund transformations for Darboux 
	integrable systems.  Indeed,   in view of Theorem B,   one can  begin to ask how  the various geometric invariants for
	Darboux  integrable systems are related though integrable extensions. One immediately identifies a 
	very natural class of integrable extensions with many important properties --- we refer to these extensions
	as  {\deffont maximally compatible}  integrable extensions for Darboux integrable systems (see  Definition  \StRef{DarbComp}). 
	We find, for example, that if
	$\bfp : (\CalE ,N) \to (\CalI, M)$ is a maximally compatible integrable extension of  a  Darboux integrable system, then the number
	of intermediate integrals for $\CalE$ is  the number of intermediate integrals for $\CalI $ plus the  fibre dimension of $\bfp$. This implies that the number of intermediate integrals for $\CalE$ is maximal among all integrable extensions of $\CalI$ with the same fibre dimension as the extension $\CalE$. 
	Another important property of these extensions is given below in Theorem C. We shall also 
	find that maximally compatible  integrable extensions for Darboux integrable systems are exactly those which can be found by the
	group-theoretic methods of Theorem A (see Theorem D below).
	
\medskip
 
	In  \cite{anderson-fels-vassiliou:2009a} it was shown that the symmetry reduction of a pair of differential systems 
	defined on  product manifolds by   a free diagonal  group action  provides a novel way of constructing Darboux integrable systems. 
	Part 3 of the paper, consisting of Sections 6, 7 and 8, explores the relationship between this group-theoretic 
	method of constructing Darboux integrable systems and the group-theoretic method in Theorem A for constructing B\"acklund
	transformations.  The key result here is  the remarkable fact that  the group-theoretic construction of 
	Darboux integrable systems  given in \cite{anderson-fels-vassiliou:2009a}  can be generalized 
	to the case of reductions by non-diagonal group actions.  The main conclusion of  Section 6   (see the remarks after Theorem 
	\StRef{Dreduce}) 	
	is that one can now construct by purely group-theoretic methods {\it   B\"acklund transformations of the form \EqRef{Intro222}
	where all the differential systems appearing in the diagram are now  Darboux integrable}. Moreover, when one of the 
	differential systems $\CalI/G_1$ or  $\CalI/G_2$ is made using  the original  diagonal  group action 
	construction from  \cite{anderson-fels-vassiliou:2009a},  we find the corresponding extension 
	$\bfp_1$ or $\bfp_2$ is maximally compatible. This observation will play an important role in Sections 7 and 8 where we 
characterize  those B\"acklund transformation for Darboux integrable equations which 
	arise from the group-theoretic construction of Section 6.

	In Section 7  we  introduce the  {\deffont Vessiot    algebra}  $\vess(\CalI)$  of a Darboux integrable system $\CalI$. 
	This a Lie algebra whose structure equations can be calculated from a certain coframe adapted to the geometry of $\CalI$.  The 
	Vessiot algebra is a fundamental diffeomorphism invariant  of $\CalI$.  In Section 7.1
	we return to the study of the interplay between the geometric  invariants  of $\CalI$ and  an extension $\CalE$ 
	where we prove the following.	
\begin{ThmC}
	Let   $\bfp\:  (\CalE, N) \to (\CalI, M)$ be an integrable extension of  Darboux systems $\CalE$ and $\CalI$.  
	If the pair $(\CalE, \CalI)$ is maximally compatible, then the induced map  (see Section 7.1)
\begin{equation}
	{\bf \tilde \bfp}_x \: \vess(\CalE)  \to \vess(\CalI)
\end{equation}
	is  a Lie algebra monomorphism for each point $x\in N$. The fibre dimension of the integrable extension  $\CalE$ is
	the co-dimension of  the image  ${\bf \tilde \bfp}_x (\vess(\CalE))$ in $\vess(\CalI)$.
\end{ThmC}

	Theorem C implies that the  fibre dimension of  a maximally compatible integrable extension for a Darboux integrable system $\CalI$
	is bounded from below by the co-dimension  of  the largest dimensional subalgebra of  $\vess(\CalI)$.
	In  the sequel to  this article  \cite{anderson-fels:2011b} we use this fact to  study  B\"acklund transformations
	for scalar Darboux integrable equations of the form $u_{xy} =f(x,y,u, u_x, u_y)$ (defining the  system $\CalI_2$). 
	It is  shown that if such an equation is non-Monge integrable and  
	admits a B\"acklund transformation $\CalB$, of fibre dimension 1, 
	to the wave equation $v_{xy} = 0$ then the projection map $\bfp_2:\CalB \to \CalI_2$
	must be maximally compatible and the Vessiot algebra of $\CalB$ must  be 2-dimensional. By Theorem C 
	the Vessiot algebra $\vess(\CalB)$  must be a subalgebra of   $\vess(\CalI)$.  Consider, for example, 
	the equation
\begin{equation}
	u_{xy}= \frac{ \sqrt{1-u_x^2} \sqrt{1-u_y^2} }{\sin u }.
\EqTag{DISO3}
\end{equation}
	This equation  is Darboux integrable and has Vessiot algebra $\lies\lieo(3)$.  Since this algebra  has no {\it  real} 
	2-dimensional subalgebras,  we deduce that \EqRef{DISO3}  cannot be  
	related to the wave equation by a  {\it real} B\"acklund transformation with 1-dimensional fibres. 
	This conclusion is at odds with part 2 of Theorem 1 in  \cite{Clelland-Ivey:2009a}.

	In Section 7.3 we calculate the Vessiot algebra for the general quotient constructions of Darboux integrable systems 
	given in Section 6.   In the particular case  of  a diagonal group action by a Lie group $G$, the Vessiot algebra is
	found to be the Lie algebra of $G$ and we call   $G$ the  {\deffont Vessiot group} for the quotient construction.

	At  this point the question now arises as to whether or not every B\"acklund transformation between pairs of Darboux 
	integrable systems necessarily arises (say at least locally)  through the group-theoretic methods developed in parts 1 and 2. 
	Motivated in part by Theorem C, our main result in this direction is the following theorem which we prove in Section 8. In stating this Theorem we assume for simplicity that the local canonical quotient representations for $\CalE$ and $\CalI$ (Theorem \StRef{Intro4}) are global.
	
\begin{ThmD}
 Let   $\bfp\:  (\CalE, N) \to (\CalI, M)$ be an integrable extension of  Darboux integrable systems $\CalE$ and $\CalI$ 
 and suppose the pair $(\CalE, \CalI)$ is maximally compatible. Then there is a Darboux integrable differential system $(\CalK, P)$ with symmetry
 group $G$  and  subgroup $H \subset G$  such that 
\begin{equation}
 (\CalI, M) \cong_{\text{\rm loc}}  (\CalK/G, P/G) \quad  \text{and} \quad (\CalE, N) \cong_{\text{\rm loc}}  (\CalK/H, P/H).
\end{equation}
 The submersion $\bfp$  is therefore locally identified with the orbit projection map from $P/H$ to  $P/G$.
\end{ThmD}

	Given a B\"acklund transformation (1.2) between Darboux integrable systems $\CalI_1$ and $\CalI_2$ , where $\bfp_2$ is 
	maximally compatible,  we can then apply Theorem D to conclude that the right hand side of the B\"acklund transformation 
	can be obtained through group quotients  just as in the right side of diagram (1.4). 
	In particular, the differential system $\CalB$ can be constructed as a  group quotient.  
	In the case of a B\"acklund transformation between  a Darboux integrable Monge-Amp\'ere equation $\CalI_2$ 
	and the wave equation $\CalI_1$, the map $\bfp_2$ is proved to be  maximally compatible. It is shown in \cite{anderson-fels:2012a} 
	that the entire B\"acklund transformation can then be constructed using Theorem A.

	Section 9, the last section of the paper,  summarizes all our main results.

\begin{ThmE} 
	Let $\CalI_2$ be a Darboux integrable differential system with quotient representation  $\CalI_2 = \CalI/G_2$. Assume that $\dim G_2 > 1$.
	Then there always exists a second symmetry group $G_1$ of  $\CalI$ for which the commutative diagram \EqRef{Intro21} can  be used to 
	construct  a B\"acklund transform  between  $\CalI_2$ and the differential system $\CalI_1 = \CalI/G_1$. The differential systems
	$\CalB$ and $\CalI_1$ are Darboux integrable and 
	$\CalI_1$  has more (functionally  independent) Darboux invariants than $\CalI_2$.
\end{ThmE}
	
	Theorem E allows us to quickly determine the  B\"acklund transformations for  all  the Darboux integrable Monge-Amp\`ere
	equations of the  type considered in   \cite{Clelland-Ivey:2005a}, \cite{Clelland-Ivey:2009a} and \cite{zvyagin:1991a}. 
	Theorem E  will be used in  a sequel to this  article \cite{anderson-fels:2011b} 
	to construct new B\"acklund transformations for equations  not of Monge-Amp\`ere type;
	for equations which are Darboux integrable at higher jet levels; for  systems of equations  in several dependent variables; 
	and for over-determined systems  in 3 independent variables.

	For Darboux integrable scalar partial differential equations in the plane it is often asserted that there
	is a chain  of B\"acklund  transformations  taking the given equation  to the  wave equation. In view of Theorems C and E
	a more general and precise formulation of this statement might be that{ \it for any Darboux integrable system 
	with solvable Vessiot algebra, there is chain  of real B\"acklund  transformations  with 1-dimensional  fibres 
	which terminates at a Darboux integrable system  with 1-dimensional Vessiot algebra.}

	It is a pleasure to thank the referee for an  extraordinarily detailed  review which  lead to significant improvements in this article.
	Support for this research was provided by grant DMS-0713830 from the National Science Foundation.

\section{Preliminaries}

	In this section we  gather together a number of definitions and basic results on integrable extensions and 
	reductions of  exterior differential systems. 	Conventions and some results from \cite{anderson-fels:2005a} are used.
	
	We assume that an EDS $\CalI$, 
	defined on a manifold $M$, has constant rank in the sense that each one  
	of its homogeneous components $\CalI^p \subset \Omega^p(M)$  coincides with the sections 
	$\CalS(I^p)$ of a  constant rank subbundle $I^p \subset \Lambda^p(M)$. If $\CalA$ is a subset of $\Omega^*(M)$, 
	we let  $\langle \CalA \rangle_{\text{alg}} $ and $\langle \CalA \rangle_{\text{\rm diff}} $ be the algebraic and differential ideals  generated by $\CalA$. 
	If $\CalI$ and $\CalJ$ are differential systems on the same manifold, we let  
\begin{equation}
	\CalI +\CalJ =  \langle \CalI \cup \CalJ \rangle_{\text{alg}}.
\end{equation} 
	Note that $\CalI + \CalJ$ is differentially closed. 
	As usual, a differential system $\CalI$ is  called a Pfaffian system if there is  a constant rank subbundle $I\subset   \Lambda^1(M)$  
	such that $\CalI = \langle \CalS(I) \rangle_\text{diff}$.  As is customary, the subbundle $I$ shall also  be referred to as 
	a  Pfaffian system.

	Let $I\subset \Lambda^1(M)$ be a Pfaffian system.  A {\deffont local first integral}  of $I$ is a smooth function 
	$f\: U \to \Real$, defined on an open set $U$,  and such that $d f \in \CalI$. For each point $x \in M$ we define
\begin{equation}
	I_x^\infty  = \spn \{\, df _x \, | \,   \text{$f$ is a local first integral, defined about $x$}\,\} .
\end{equation}
	We shall always assume that  $I^\infty   = \cup _{x\in M} I_x^\infty$ is a constant rank bundle on $M$.
	It is easy to check that $I^\infty$ is  the  (unique) maximal,  completely integrable,  Pfaffian subsystem of  $I$.
	
	The bundle $I^\infty$ can be computed algorithmically from the derived  flag of $I$. 
	The derived system $I'\subset I$ is defined pointwise by
$$
I'_p= \spn \{ \, \theta_p \ | \  \theta \in \CalS(I) \text{ such that } d \theta \equiv 0 \mod \ I \, \}.
$$
	Letting $I^{(0)}=I$ and assuming $I^{(k)}$ is constant rank we define $I^{(k+1)}$ inductively by
	$I^{(k+1)} = (I^{(k)})'  \quad\text{for} \quad k=0,1\ldots, N,$
	where $N$ is the smallest integer where $I^{(N+1)} = I^{(N)}$.  Then $I^{(N)}$ is completely integrable and
	$ I^\infty = I^{(N)}.$	More information about the derived flag of a Pfaffian system can be found in
	\cite{bryant-chern-gardner-griffiths-goldschmidt:1991a} and \cite{ivy-langsberg:2003}.  

\subsection{Integrable Extensions of Differential Systems}

	We recall the definition of an integrable extension \cite{bryant-griffiths:1995a}. First, let $\bfp :N \to M$ be a  submersion 
	and $\CalI$ an EDS on $M$.  An EDS  $\CalE$ on $N$
	is called an {\deffont integrable extension}  of  $\CalI$ if there exists a subbundle  $J \subset  \Lambda^1(N)$ such that 
\begin{equation}
	\text{\bf [i] }  \rank  J = \dim N- \dim M, 
	\quad \text{\bf [ii] }  \ann(J) \cap \text{ker } (\bfp_*)=0 
	\quad\text{and \bf [iii] }
	\CalE = \langle \CalS(J) \cup  \bfp^*(\CalI) \rangle_\text{alg} .
\EqTag{IntExtDef}
\end{equation}
	Here $\ann(J)$ is the subbundle of vectors in $TM$ which annihilate the 1-forms in $J$. The second condition in \EqRef{IntExtDef}
	states that $J$ is transverse to  $\bfp$.

	The integrals manifolds of $\CalE$ and $\CalI$ are related as follows. 
	Let  $s:P \to N$ be  an immersed integral manifold of $\CalE$ and let  $ \tilde s = \bfp \circ s:P\to M$. Then condition {\bf [ii]}
	of \EqRef{IntExtDef} implies that $\tilde s$ is an immersion and hence  {\bf [iii]} implies that $\tilde s$
	is an integral manifold of $\CalI$.  
	Conversely, if $\tilde s:P\to M$ is an integral manifold of
	$\CalI$ and $Q =  \bfp^{-1}(\tilde s(P))$,  then  the third condition in  \EqRef{IntExtDef}   implies that
	$\CalE|_Q$ is a  Frobenius system. Consequently integral manifolds of $\CalI$ can be lifted locally to integrable manifolds
	of $\CalE$ by integrating a system of ordinary differential equations. 

	A subbundle $J$  satisfying the three properties \EqRef{IntExtDef} is called an 
	{\deffont admissible} subbundle for the extension $\CalE$. We note that conditions {\bf [i]} and {\bf [ii]} of 
	\EqRef{IntExtDef} imply that
\begin{equation}
	\Lambda^1(N)  = J \oplus  \bfp^* (\Lambda^1(M)) .
\EqTag{IntExt5}
\end{equation}
	If   $\CalE^1 = \CalS(E^1)$  and   $\CalI^1 = \CalS(I^1)$, then  condition {\bf [iii]} implies
\begin{equation}
	E^1 =  J \oplus \bfp^* (I^1). 
\EqTag{IntExt3}
\end{equation}
	Thus,  if $ \{\, \xi^u \,\}$ are 1-forms on $N$
	which define a local basis for $\CalS(J)$, then  condition {\bf [iii]} in \EqRef{IntExtDef}  also states that
\begin{equation}
	d\, \xi^u  \equiv 0  \mod \{\bfp^*(\CalI) ,\, \xi^u \} \quad \text{or} \quad  d\, \xi^u \equiv 0 \mod \{ E^1,\, \bfp^*(I^2) \} .
\EqTag{IntExt4}
\end{equation}
	Any subbundle $J \subset T^*(N)$ which satisfies  \EqRef{IntExt5},  \EqRef{IntExt3} and \EqRef{IntExt4} is  admissible	

A few additional properties of integrable extensions will be needed.

\smallskip
\newcommand\IntExtAdmiss{{\bf IE\,[i]}}
\noindent
{\bf IE\,[i]}
	If  $\bfp \: \CalE \to \CalI$ is an integrable extension, then {\it any}  complementary subbundle $J$ to $\bfp^*(I^1)$ in $E^1$ 
	satisfies \EqRef{IntExtDef} and is therefore admissible.

\smallskip
\newcommand\IntExtDer{{\bf IE[ii]}}
\noindent
{\bf IE\,[ii]} 
	If $\CalI$ is a Pfaffian system, then $\CalE$ is a Pfaffian system with $E=J \oplus \bfp^*(I)$ and, using \EqRef{IntExtDef},
\begin{equation}
	\rank(E) = \rank( I )+ v  \quad \text{ and} \quad \rank (E') = \rank (I') + v, 
\EqTag{IntEv}
\end{equation}
	where  $v = \dim N - \dim M = \rank(J)$.

\smallskip
\newcommand\IntExtInt{{\bf IE\,[iii]}}
\noindent
{\bf IE\,[iii]} 
	If $\CalI$ is a completely integrable Pfaffian system, then equation \EqRef{IntEv} implies that $\CalE$ is completely integrable. 
	Moreover, if $ s:P \to N$ is the maximal integral manifold though $x_0\in N$  for $\CalE$ and $\tilde s:\tilde P \to M$ 
	is the maximal integral manifold through $\bfp(x_0) \in M$ for $\CalI$, then $\bfp\circ s:P \to \tilde s(\tilde P)$ is a local diffeomorphism.

\smallskip
\newcommand\IntExtProl{{\bf IE\,[iv]}}
\noindent
{\bf IE\,[iv]}  
	Let $\bfp : (\CalE, N) \to (\CalI, M)$  be an integrable extension with $J\subset T^*N$ an admissible bundle 
	and let $\pi_N : (\CalE^{[1]}, N^{[1]}) \to (\CalE, N)$ and $\pi_M : (\CalI^{[1]}, M^{[1]}) \to (\CalI, M)$ 
	be the prolongations of $\CalE$ and $\CalI$ to the space $N^{[1]}$ and $M^{[1]}$ 
	of $k$-dimensional integral elements for $\CalE$ and $\CalI$.
	Then $\bfp$ lifts to a unique map $\bfp^{[1]} :(\CalE^{[1]}, N^{[1]}) \to (\CalI^{[1]}, M^{[1]}) $ 
	which covers $\bfp$ and  $\CalE^{[1]}$ is an integrable extension of $\CalI^{[1]}$ with admissible bundle $\pi_N^*(J)$.
 
\smallskip
\newcommand\IntExtInftymodp{{\bf IE\,[v]}}
\noindent
{\bf IE\,[v]} 
	Let  $\bfp:(E, N) \to (I, M)$ be an integrable extension of a Pfaffian system $I$.
	If $J$  is an  admissible bundle for this extension, then  
\begin{equation}
	\rank(I^\infty) \leq \rank(E^\infty/ \bfp) \leq  \rank( I^\infty)  + \rank (J).
\EqTag{IntEviii}
\end{equation}

Properties \IntExtAdmiss, \IntExtDer\ and \IntExtInt\ are easily checked. See Appendix A for the proof of  \IntExtProl\ and  \cite{anderson-fels:2012a} for the proof of \IntExtInftymodp.

\subsection{Reduction of Differential Systems}

	Recall from the introduction that if $\bfp :N \to M$ is a smooth surjective submersion and $\CalI$ an EDS on $N$,  
	then the {\deffont reduction of $\CalI$} with respect to $\bfp$ which is denoted by $\CalI/\bfp$, is the EDS on $M$ defined by
\begin{equation}
	\CalI/\bfp = \{ \ \theta \in \Omega^*(M) \ | \ \bfp^* \theta \in \CalI \ \}.
\EqTag{Imodp}
\end{equation}
	Note that $\CalI/\bfp$ is not necessarily constant rank without additional hypotheses. 
	If the fibres of $\bfp$ are connected, then the reduction $\CalI/\bfp$ can be computed 
	using Corollary II.2.3 of \cite{bryant-chern-gardner-griffiths-goldschmidt:1991a} which states that
	$\theta \in \CalI$ satisfies $\theta = \bfp^* \bar \theta$ for some $\bar \theta \in \CalI/\bfp$ if and only if,
\begin{equation}
X\hook \theta = 0,\ \ X \hook d \theta = 0 \quad {\rm for \ all \ } \ X \in \ker (\bfp_*) .
\EqTag{IredX}
\end{equation}
	Similarly, if  $I \subset \Lambda^p(N)$, then the reduction of the bundle $I$ with respect to $\bfp$ is
\begin{equation}
	I/\bfp = \{ \ \theta \in \Lambda^p(M)\ | \ \bfp^* \theta \in I  \ \}.
\EqTag{bundlemodp}
\end{equation}

	We now specialize \EqRef{Imodp} to the case of  reduction by a Lie  group.
	Let $G$ be a finite dimensional Lie group acting on $M$ with left action $\mu:G \times M \to M$.
	Define maps  $\mu_x :G \to M$ and $\mu_g :M  \to M$ by
\begin{equation*}
	\mu_x(g)  = \mu(g, x) = g \cdot x = \mu_g(x) \quad\text{for $x\in M$ and  $g \in G$}.
\EqTag{actionmaps}
\end{equation*}
For each  $Z_e \in T_eG$, the corresponding infinitesimal generator for the action $\mu$ is 
the vector field $X$ on $M$ defined pointwise by
\begin{equation}
	X_x  = (\mu_x)_*(Z_e)  \quad \text{for all $x \in M$.}
\EqTag{infgen}
\end{equation}
The set of all  infinitesimal generators for the action $\mu$ is a Lie algebra of vector fields on $M$ which we
	denote by  $\Gamma_G$. If $\lieg$ is the Lie algebra of {\it right} invariant vector fields on $G$
	and $Z\in \lieg$, then the map $\tilde \mu: \lieg \to \Gamma_G$ defined by \EqRef{infgen} is a Lie algebra homomorphism.

	Let $\bfGamma_G \subset TM$ be the integrable distribution generated by the point-wise span of $\Gamma_G$. 
	We will assume that all actions are  {\deffont regular} in the sense that the orbit space $M/G$ has a smooth manifold structure
	such that the  canonical projection $\bfq_G \: M \to M/G$ is a smooth submersion 
	and   hence $\bfGamma_G = \ker (\bfq_{G*})$.  Since each orbit of $G$ is the inverse image of a point in $M/G$,
	the implicit function theorem implies that the orbits are  imbedded submanifolds.
	
	The group $G$ acting on $M$ is a {\deffont symmetry group of $\CalI$} if,
	for each $g\in G$ and $\theta \in \CalI$,   $\mu_g^*(\theta) \in \CalI$. Under these circumstances
	we define the  {\deffont symmetry reduction of $\CalI$  by G} to be
\begin{equation}
	\CalI/G = \CalI/\bfq_G=  \{\ \bar \theta \in \Omega^*(M/G) \ | \ \bfq_G^*(\bar \theta ) \in \CalI  \ \}.
\EqTag{IGdef}
\end{equation}
 	In other words $\CalI/G$ is the reduction given by equation \EqRef{Imodp} with respect to the submersion $\bfq_G:M\to M/G$. 
	However, by utilizing the $G$-invariance of $\CalI$, the computation of a local basis of sections for $\CalI/G$ 
	can now be done algebraically \cite{anderson-fels:2005a}. 
	This fact is not necessarily true for the reduction in \EqRef{Imodp} for generic $\bfp$. 
	In analogy with \EqRef{IredX},  a form $\theta \in \Omega^p(M)$ satisfies   $\theta=\bfq_G^*(\bar \theta ) $ for some 
	$\bar \theta \in \Omega^p(M/G)$ if and only if $\theta$ is {\deffont  $G$-basic}, that is, $G$  semi-basic and $G$ invariant so that
\begin{equation}
	X\hook \theta = 0 \quad\text{and} \quad \mu_g^*(\theta) = \theta
\	\quad\text{for all  $X \in \bfGamma_G$ and  $g \in G$.}
\end{equation}

	Likewise,  if $A \subset \Lambda^p(M)$ is a $G$-invariant subbundle, then
\begin{equation}
	A/G = A/\bfq_G=  \{\ \bar \theta \in  \Lambda^p(M/G) \ | \ \bfq_G^*(\bar \theta ) \in A \ \}.
\EqTag{AmodG}
\end{equation} 
As a cautionary remark, we observe that if $\CalI$ is a $G$-invariant Pfaffian differential system with  $\CalI = \langle \CalS(I)\rangle_{\text{diff}}$,
	then it is generally not true that $\CalI/G$ is the Pfaffian system for $I/G$, that is
\begin{equation*}
	\CalI/G  \neq   \langle \CalS(I/G)\rangle_{\text{diff}}.
\end{equation*}
However, it is true that $ \CalS(I/G)\rangle_{\text{diff}}\subset \CalI/G$.	See \cite{anderson-fels:2005a} for examples.
	
	A symmetry group $G$ of an EDS $\CalI$ is said to be {\deffont transverse} to $\CalI$ if
\begin{equation}
	 \ann(I^1) \cap \bfGamma_G = 0.
\EqTag{Itrans}
\end{equation}
	This transversality condition, which holds for all the examples and application we consider, will be an essential hypothesis 
	 for  almost all of the results in this paper.  Therefore,  it is useful to summarize a few relevant facts from  \cite{anderson-fels:2005a} about symmetry reduction and the role of transverse actions.

	The first step in computing the reduction in \EqRef{AmodG} is to define, for any subbundle $A\subset \Lambda^p(M)$, 
 the subset $A_{G, \semibasic}\subset A$ of $G$ semi-basic forms as
\begin{equation}
	A_{G,\semibasic} = \{\,\alpha \in A \, | \, X \hook \alpha = 0 \text{ for all $ X \in \bfGamma_G$} \} .
\end{equation}
It is clear that if $A$ is $G$-invariant, then $A_{G,\semibasic}$ is $G$-invariant.  
	If, in addition, $A_{G,\semibasic}$ is a constant-rank
	bundle then there is a bundle $\bar A \subset \Lambda^p(M/G)$, of the same rank as  $A_{G,\semibasic}$,
	satisfying 
\begin{equation}
	\bfq_G^*(\bar A) = A_{G,\semibasic}, \quad \text{namely} \quad   \bar A = A_{G,\semibasic}/G.
\EqTag{Asb}
\end{equation}
	Furthermore, about each point $x \in M$, there is an $G$-invariant open set  $U$   
	and a local basis  $\{\alpha^i\}$ for $\bar A|_{\bfq_G(U)}$  such that 
\begin{equation}
A_{G,\semibasic}|_U=	\spn  \{ \bfq^*_G(\alpha^i) \} \, .
\EqTag{Glocb}
\end{equation}
Equation \EqRef{Glocb} shows that  $A_{G,\semibasic}$ always admits a local basis of $G$ invariant $p$-forms. 
	
	In the particular case of interest, namely  when $G$ is symmetry group of a differential system $\CalI$ 
	and $G$ acts transversely to $\CalI$, the fact $\ann(\bfGamma_G) = \Lambda^1(M)_{G,\semibasic}$ allows us to re-write condition \EqRef{Itrans} in the equivalent form
\begin{equation}
I^1 +   \Lambda^1(M)_{G,\semibasic} =  \Lambda^1(M).
\EqTag{Itrans2}
\end{equation}
Using equation \EqRef{Itrans2} we showed in \cite{anderson-fels:2005a} that transversality automatically  ensures that the spaces $I^p_{G,\semibasic}$ are constant rank. By applying the above observations to the bundles $I^p$,  the quotient $I^p/G$ is constant rank.  In the special case of the bundle $I^1$ of 1-forms  in $\CalI$, the rank is easily computed from the transversality condition \EqRef{Itrans2} to be
\begin{equation}
	\rank (I^1/G) = \rank(I^1_{G,\semibasic}) =\rank( I^1 \cap   \Lambda^1(M)_{G,\semibasic})= \rank(I^1) - \rank( \bfGamma_G).
\EqTag{rankI1}
\end{equation}

A fundamental consequence of  transversality is the following.

\begin{Theorem}
\StTag{IntEntCor1}
	Let $G$ be  a symmetry group of an EDS $\CalI$ on a manifold $M$. Assume $G$ acts regularly on $M$ and transversely to $\CalI$.  
	Then $\CalI$ is an integrable extension of the reduced differential system $\CalI/G$.
\end{Theorem}
\begin{proof} To prove this theorem we shall construct a set of  algebraic generators for $\CalI$  adapted to the $G$-action.
	Using equation \EqRef{Itrans2} we may choose bundles  $J$ and  a $G$-invariant,  $G$ semi-basic bundle $W$   such  that
\begin{equation}
	I^1_{G, \semibasic} \oplus J = I^1\quad \text{and}\quad  I^1_{G, \semibasic} \oplus W = \Lambda^1(M)_{G, \semibasic}. 
\EqTag{JW}
\end{equation}
	The bundle $J$ is any complement to $I^1_{G, \semibasic}$ in $I^1$, and by equation \EqRef{rankI1} $\rank(J)=\rank(\bfGamma_G)$. To  construct $W$, let 
	$\barW$ be any complement to $I^1/G$ in $\Lambda^1(M/G)$ and set $W =  \bfq_G^*(\barW)$.
	The transversality condition \EqRef{Itrans2} then implies 
\begin{equation*}
 \Lambda^1(M)_{G, \semibasic} \oplus  J = \Lambda^1(M) 
\end{equation*}
in which case
\begin{equation}
	I^1_{G, \semibasic} \oplus J \oplus  W = \Lambda^1(M).
\EqTag{DecompI1}
\end{equation}

Using the decomposition of $\Lambda^1(M)$ in \EqRef{DecompI1}  we may choose an open set $U$ and a set of generators for $\CalI|_U$ such that
\begin{equation*}
\CalI|_{U} = \langle\, \theta^i,\,  \theta^a_{G,\semibasic}, \,  \tau^\alpha_{G,\semibasic} \, \rangle_\text{alg}
\end{equation*}
where $\theta^i \in \Omega^1(J|_U)$, $\theta^a_ {G,\semibasic}\in \Omega^1(I^1_{G,\semibasic}|_U)$, and $ \tau^\alpha_{G,\semibasic}   \in \Omega^{*}(W|_U)$. By applying the argument leading to equations \EqRef{Asb} and \EqRef{Glocb} to the case of $I^p$, the
generators for $\CalI|_U$ above may be refined so that
\begin{equation}
\begin{gathered}
	\CalI|_{U} = \langle\, \theta^i,\,  \theta^a_{G}, \,  \tau^\alpha_{G} \, \rangle_\text{alg} \quad {\rm and}\
(\CalI/G)|_{\bfq_G(U)} = \langle\,   \bar \theta^a, \, \bar  \tau^\alpha \, \rangle_\text{alg},  \quad \text{where}
\\[6pt]
	J|_U  = \spn \{\, \theta^i\, \},
	\quad
 	I^1_{G, \semibasic}|_U  = \spn \{\, \theta^a_G= \bfq_G^* \bar \theta^a \, \} 
	\quad  \text{and}  \quad \tau^\alpha_G  = \bfq_G^* \bar \tau^\alpha   \in \Omega^{*}(W|_U).
\end{gathered}	
\EqTag{IntExt15}	
\end{equation}
We emphasize that the forms $\theta^a_G$ and $\tau^\alpha_G$ are now $G$-basic. For additional details see Appendix A and Appendix B in \cite{anderson-fels:2005a}.

	We call   \EqRef{IntExt15}  a set  of $G$-adapted algebraic generators for $\CalI$.
	By definition, $\CalI$ is differentially closed and hence
\begin{equation}
	d\, \theta^i \equiv 0 \mod  \langle\, \theta^i, \theta^a_{G}, \, \tau^\alpha_{G} \rangle_\text{alg}.
\EqTag{IntExt7}
\end{equation}
	Together, equations  \EqRef{IntExt15}	and  \EqRef{IntExt7} prove that $\CalI$ is an integrable extension
	of $\CalI/G$.
	We shall give a refinement of these structure equations for the case of free group actions in Section 2.3.
\end{proof}
	We shall also need the following.
\begin{Theorem}
\StTag{RedGInfty}
	Let $I$ be a Pfaffian system on a manifold $M$.
	Suppose $G$ acts regularly on $M$ with quotient map $\bfq_G\:M  \to M/G$,  $G$ is a symmetry  group of $I$, and $G$
	acts transversely to  $I$.
	If $(I^\infty)_{G,\semibasic}$ and  $(I/G)^\infty$ are constant rank vector bundles, then 
\begin{equation}
	(I^\infty)_{G, \semibasic} =  \bfq_G^*((I/G)^\infty) \quad \text{and}\quad  (I/G)^\infty = I^\infty/G.  
\EqTag{crank3}
\end{equation}
\end{Theorem}

\begin{proof}  
	We prove the second equation in \EqRef{crank3}. The first  equation then follows from \EqRef{Asb}.
	We first remark that the transversality hypothesis implies that $I/G$  has constant rank. 
	Since the pullback of an integrable Pfaffian system is integrable and $\bfq_G^*$  is injective, 
	it follows that $\bfq^*_G( (I/ G)^\infty) $ is a  constant rank
	integrable subbundle of $I$ and therefore $ \bfq^*_G( (I/ G)^\infty)  \subset  I^\infty$.  By definition \EqRef{IGdef},  this implies that 
\begin{equation*}
	(I/ G)^\infty  \subset  I^\infty/G.
\end{equation*}
	To prove the reverse inclusion we deduce, from the assumption that 
	$(I^\infty)_{G,\semibasic}$ is of constant rank, that is,  $I^\infty/G$ is a constant rank  subbundle of $I/G$. 
	But  $I^\infty/G$ is integrable \cite{fels:2007a}  (Theorem 2.9) and hence $I^\infty/G \subset (I/G)^\infty$ and the lemma is established.
\end{proof}
	
\subsection{Symmetry Reduction of Pfaffian Systems by Free Group Actions}

	It will be important to have a refinement of the foregoing general results on reduction of differential systems
	for the special case where the action of $G$, in addition to be being regular and transverse,  is  also free. In other words,
	we now consider a Pfaffian system $I$ defined on a left principal $G$ bundle $M$.
	Recall that a local trivialization for $\bfq_G:M \to M/G$ consists of a $G$-invariant open set $U \subset M $ and a diffeomorphism 
	$\Phi \:U \to  \bfq_G(U) \times G$ where  $\Phi(x) = (\bfq_G(x),\phi(x))$  
	and $\phi \:U \to G$ is  $G$-equivariant, that is, $\phi(\mu(g,x))= g \,\phi(x)$. 
	
	For the $r$-dimensional Lie group $G$ we make the following choices. 
	Let $Z_i,\, 1\leq i \leq r$ be a basis of right invariant vector fields, and let
	$\tau^i$ be the dual right invariant one-forms to $Z_i$ on $G$.  Let  $[\, Z_i, \, Z_j\,] =  C_{ij}^k Z_k $.
	Define the matrix-valued function $\bflambda:G \to GL(r, \real)$ by
\begin{equation}
\Ad^*(g) (\tau^i) = \lambda(g)^i_j \,\tau^ j,
\EqTag{deflambda}
\end{equation}
	where $\Ad^*$ is the co-adjoint representation of $G$.  Equation \EqRef{deflambda}  implies that
\begin{equation}
L_g^*\, \tau^i = \lambda(g^{-1})^i_j \tau^j,
\EqTag{lambdap1}
\end{equation}
	where $L_g$ is left multiplication on $G$ by $g$.
	Since $\Ad^*(g\,g') = \Ad^*(g) \Ad^*(g')$ we  also deduce that
\begin{equation}
\lambda(g \,g')^i_j= \lambda(g')^i_k \lambda(g) ^k _j \, .
\EqTag{lambdap2}
\end{equation}
\begin{Lemma}
	The exterior derivative of the matrix-valued function $\bflambda:G \to GL(r, \real)$ is
\begin{equation}
d\lambda^i_j = \lambda^i_k C^k_{lj}\, \tau^l .
\EqTag{dlambdaG}
\end{equation}
\end{Lemma}
\begin{proof}
	With $\{ Z_i \}$ the basis of right invariant vector-fields we  first compute $\bflambda_* Z_i(e) $ by
\begin{equation}
\begin{aligned}
	\frac{d \ }{ dt}  \bflambda \bigl( \exp(t Z_i) \bigr) (\tau^j(e))|_{t=0} 
	&= \frac{d \ }{ dt} (L^*_{\exp(-tZ_i)} \circ R_{\exp(t Z_i)}^* (\tau^j(e))|_{t=0} \\
	& = (-\CalL_{Z_i} \tau^j)(e) = C^j_{i k} \tau^k(e)
\end{aligned}
\end{equation}
	and therefore
$$
	d \lambda^i_j(e) = C^i_{kj} \tau^k(e) .
$$
	Then we may use property \EqRef{lambdap2} to compute $\bflambda_*(Z_l(a))$. 
	Since $Z_i$ is a right invariant vector-field we have
\begin{equation}
	\frac{d \ }{ dt} \lambda(\exp(t Z_\ell(a) )^i_j |_{t=0} = \lambda^i_k(a) \frac{d \ }{ dt} \lambda^k_j(\exp( tZ_\ell))|_{t=0} \\
	= \lambda^i_k(a) C^k_{\ell j} 
\end{equation}
which establishes \EqRef{dlambdaG}.
\end{proof}

\begin{Theorem}
\StTag{BigLemma}
	Let $\CalI$ be  a Pfaffian system on a manifold $M$. Let $G$ be a symmetry group of $\CalI$ which acts freely
	and regularly on $M$ and transversely to $\CalI$. Let $\{Z_1, Z_2 , \dots Z_r\}$ be a basis for the right invariant vector fields on $G$	
	with structure equations $[\, Z_i, \, Z_j\,] =  C_{ij}^k Z_k $
	and  let  $\{X_1, X_2 , \dots X_r\}$ be the corresponding basis of infinitesimal generators for the action of $G$ on $M$ 
	(see \EqRef{infgen}).  

	Then about each point $x\in M$ there exists a $G$-invariant open set $U$ and a coframe $ \{\, \theta^i, \eta^a, \sigma^\alpha\, \} $ on $U$ such that:

\noindent
\smallskip
 {\bf[i]} $\CalI = \langle \theta^i, \eta^a \rangle_{\text{\rm diff}}$; 

\noindent
\smallskip
{\bf[ii]} the forms  $\eta^a$, $\sigma^\alpha$ are $G$-basic;

\noindent
\smallskip
	{\bf [iii]} the forms $\theta^i$ satisfy $\theta^i(X_j) = \delta^i_j$,

\noindent
\smallskip
 {\bf [iv]} For any $g \in G$, $\mu_{g}^* \theta^i= \lambda(g^{-1})^i_j \theta^j,$ where $\lambda^i_j(g) $ is the matrix defined in equation \EqRef{deflambda}; and

\noindent
\smallskip
 {\bf[v]} these forms satisfy the  structure equations
\begin{equation}
\begin{aligned}
	d\sigma^\alpha & = 0,  \quad   
\\d\eta^a & =  E^a_{\alpha \beta}\,\sigma^\alpha \wedge \sigma^\beta + F^a_{ b \alpha}\,\eta^b \wedge \sigma^\alpha, \quad \text{and}
\\
	d \theta^i  &=  A^i_{\alpha\beta}\, \sigma^\alpha \wedge \sigma^\beta  +  B^i_{a\beta}\,\eta^a\,\wedge \sigma^\beta - \frac12 C^i_{jk}\,\theta^j \wedge \theta^k .
\end{aligned}
\EqTag{BigStr}
\end{equation}
\end{Theorem}

\begin{proof}
	Let $n = \dim(M)$,  $r = \dim(G)$,  and $p = \rank I^1$. Again, because the action of $G$ is free and regular,
	we may choose a local trivialization with open set $U$ and map $(\bfq_G,\phi):U \to \barU \times G$,  
	where $\barU$ is an open set in $M/G$ which we can  choose to be small enough to be the domain of a coordinate chart.
	Let  $J^1, J^2, \ldots, J^{n-r}$ be a set of coordinates on $\barU$, whose  pullbacks to $U$,  
	which again we denote by $J^1,\ldots, J^{n-r}$, are  $G$-invariant. Since the rank of $\Lambda^1(M)_{G,\semibasic}$ is $n-r$, 
	the differentials of these invariant functions give a basis for $\Omega^1(U)_{G,\semibasic}$.

	On account of  \EqRef{IntExt15}, we may assume  that on the $G$-invariant open set $U$,
	$\{\, \eta^1,\eta^2, \dots, \eta^{p-r} \, \}$ is a basis of $G$-invariant sections for $I^1_{G,\semibasic}$, restricted to $U$.
	Since $I^1_{G,\semibasic} \subset  \Lambda^1(M)_{G,\semibasic}$,
	we can choose a set of differentials $\sigma^\alpha =d J^{k_\alpha}$ complementary
	to the 1-forms  $\eta^a$.  This gives  $\Omega^1(U)_{G,\semibasic} = \spn \{\,\eta^a,\, \sigma^\alpha\, \}$.   
	To simplify the notation in what follows, we re-label the invariant functions as 
	$J^1$,$J^2, \ldots, J^{n-p}$, $K^1$, $K^2\ldots K^{p-r}$ so that now  $\sigma^\alpha =d J^{\alpha}$. 
	The forms $\eta^a$ are $G$-basic 
	and can therefore be written as 
\begin{equation}
	\eta^a = P^a_b \,dK^b + Q^a_\alpha \,\sigma^\alpha,
\EqTag{etadef}
\end{equation}
	where the coefficients are functions of the invariants $J^\alpha, K^a$. Because the 1-forms $\{\,\eta^a,\, \sigma^\alpha\, \}$ 
	are point-wise  linearly independent the matrix $P^a_b$ is invertible and therefore we can replace the forms  in \EqRef{etadef} 
	by 
\begin{equation}
	\eta^a = dK^a + R^a_\alpha \,\sigma^\alpha,  \quad\text{where} \quad  R^a_\alpha  =   [P^{-1}]^a_b Q^b_\alpha
\EqTag{etadef1}
\end{equation}	
	and still preserve the equation $\Omega^1(U)_{G,\semibasic} = \spn \{\,\eta^a,\, \sigma^\alpha\, \}$.  
	The forms $\sigma^\alpha$ are closed by construction.  The forms $\eta^a \in I$ and the
	structure equations for the $\eta^a$  in  \EqRef{BigStr} are an easy consequence of  \EqRef{etadef} and  \EqRef{etadef1}.

	At this point {\bf[i]} and {\bf[ii]}  and the  first 2 structure equations in  \EqRef{BigStr}  are proved. We  turn now to 
	{\bf[iii]} and {\bf [iv]}.  Let $\omega^i = \phi^* (\tau^i)$, where $\tau^i$ are the right-invariant forms on $G$ defined above 
	and dual to $Z_i$. Note that $\omega^i(X_j)=\delta^i_j$, and
$$
	d\omega^i = -\dfrac12 C^i_{jk} \omega^j \wedge \omega^k .
$$
	The forms $\{\omega^i, \eta^a ,\sigma^\alpha\}$ are a basis of sections of $T^*U$.

	Now choose  1-forms $\{\, \tilde \theta^1, \ldots \tilde \theta^p \, \}$ so that  
	$\{\, \tilde \theta^i, \eta^a\,\}$ is a basis of sections for $I^1$,
	restricted to $U$.   We can write
\begin{equation}
	\tilde \theta^i  = \tilde P^i_j  \omega^i + \tilde S^i_a \eta^a +  \tilde T^i_\alpha \sigma ^\alpha.
\end{equation} 
	Transversality implies that the matrix  $\tilde P^i_j $   is invertible and therefore 
	there are  smooth functions  $S^i_a, T^i_\alpha $  on $U$ such that, for each $i = 1\ldots r$,
\begin{equation}
	\omega^i + S^i_a \eta^a + T^i_\alpha \sigma ^\alpha \in \CalI .
\end{equation}
	Since $\eta^a\in \CalI$, this implies that 
 \begin{equation}
	\theta^i = \omega^i +T^i_\alpha \sigma ^\alpha \in \CalI .
\EqTag{deftheta}
\end{equation}
	At this point the 1-forms $\{\, \theta^i, \eta^a, \sigma^\alpha \, \}$ satisfy  parts ${\bf [i],[ii],[iii]}$ of the theorem.

To show {\bf [iv]} first note that by the equivariance of $\phi$ and equation \EqRef{lambdap1},  
\begin{equation}
\mu_g^* (\omega^i) =\mu_g^* \phi^* (\tau^i)= \phi^*(  \lambda(g^{-1})^i_j) \omega^j .
\EqTag{mugom}
\end{equation}
Now using the $G$-invariance of $\sigma ^\alpha$ and equation \EqRef{mugom} we have (where $g\cdot  p = \mu_g(p)$)
\begin{equation}
\mu_g^*\theta^i(g \cdot p) = \lambda(g^{-1}) ^i_j \omega^j(p) +  T^i_\alpha(g\cdot  p) \sigma ^\alpha(p) .
\EqTag{mugth}
\end{equation}
Since $\CalI$ is $G$-invariant $\mu_g^* \theta^i \in \CalS(I^1|_U)$, and equation \EqRef{mugth} leads to
$$
\mu_g^*\theta^i = \lambda(g^{-1})^i_j \theta^j \quad \text{and} \quad T^i_\alpha(g \cdot p) = \lambda(g^{-1}) ^i_j T^j_\alpha(p) .
$$
This proves property {\bf [iv]}.

Finally we complete the proof of part {\bf [v]} of the theorem.  By direct calculation using equation \EqRef{deftheta}, we have
\begin{align*}
	d \theta^i  &	= d\omega^i +  dT^i_\alpha\,\wedge \sigma^\alpha = -\frac12 C^i_{jk}\,\omega^j\wedge \omega^k + dT^i_\alpha \,\wedge \sigma^\alpha
\\
	& = -\frac12 C^i_{jk}\,\theta^j\wedge \theta^k  + C^i_{jk} T^k_\beta\, \theta^j\wedge \sigma^\beta
	- \frac12 C^i_{jk}T^j_\alpha T^k_\beta\, \sigma^\alpha\wedge \sigma^\beta +  d T^i_\alpha\,\wedge\sigma^\alpha.
\end{align*}
	There are no 2-forms of the kind $\theta^i \wedge \eta^a $ in this formula for $d\theta^i$.
	Moreover, since the vector fields $X_\ell$ are infinitesimal  symmetries for $\CalI$,
	it follows that $X_\ell\hook d\,\theta^i \in \CalI$
	and therefore the terms of the kind  $\theta^i \wedge \sigma^\alpha$
	(appearing in  $C^i_{jk}T^k_\beta\, \theta^j\wedge\sigma^\beta$ and $d T^i_\alpha\,\wedge\sigma^\alpha$) must cancel out.
	This leads to the structure equations \EqRef{BigStr}  for the 1-forms $\theta^i$.
\end{proof}

\section{B\"acklund Transformations and Symmetry Reduction of  EDS}

The goal of this section is  to prove Theorem A.  In Section 3.1 we prove  two  lemmas regarding  commutative diagrams  of  differential systems.  In Section 3.2  we  give simple conditions under which an orbit projection map defines an integrable extension from which  Theorem  A will immediately follow.

\subsection{Reductions of EDS and commutative diagrams}

In this section we  prove the first statement in  Theorem A. We start with the following lemma.

\begin{Lemma} \StTag{Gcomd} Let $\CalI$ be a differential system on  a manifold $P$ and let 
\begin{equation}
\begin{gathered}
\begindc{\commdiag}[3]
\obj(-25, 0)[I]{$P$}
\obj(0, 0)[H]{$N$}
\obj(0, -17)[I1]{$M$}
\mor{I}{H}{$\bfp_1$}[\atleft, \solidarrow]
\mor{H}{I1}{$\bfp_2$}[\atleft, \solidarrow]
\mor{I}{I1}{$\bfp_3$}[\atright, \solidarrow]
\enddc
\EqTag{CDsubs}
\end{gathered}
\end{equation}
	be a commutative diagram of manifolds and suppose that $\bfp_1$ and $\bfp_3$ are surjective submersions. 
	Then $\bfp_2$ is a surjective submersion and 
\begin{equation}
	(\CalI/\bfp_1)/\bfp_2 = \CalI/\bfp_3 
\EqTag{CDeds0}
\end{equation}
\end{Lemma}

\begin{proof} 
	The fact that   $\bfp_2$ is  a surjective submersion is easily checked.
	The definitions of $(\CalI/\bfp_1)/\bfp_2$,  $\CalI/\bfp_1$ and $\CalI/\bfp_3$  imply, in turn, that
\begin{equation*}
\begin{aligned}
	(\CalI/\bfp_1)/\bfp_2
&	= \{ \, \theta \in \Omega^*(M) \ | \ \bfp_2^* \theta \in \CalI/\bfp_1 \, \}
 	=  \{ \, \theta \in \Omega^*(M) \ | \ \bfp_1^*(\bfp_2^* \theta) \in \CalI\, \}
\\
& 	= \{ \, \theta \in \Omega^*(M) \ | \ \bfp_3^* \theta \in \CalI\, \} = \CalI/\bfp_3 
\end{aligned}
\end{equation*}
\end{proof}

When ever  \EqRef{CDeds0} holds , we shall say that 
\begin{equation}
\begin{gathered}
\begindc{\commdiag}[3]
\obj(-25, 0)[I]{$\CalI$}
\obj(0, 0)[H]{$\CalI/\bfp_1$}
\obj(0, -17)[I1]{$\CalI/\bfp_3$}
\mor{I}{H}{$\bfp_1$}[\atleft, \solidarrow]
\mor{H}{I1}{$\bfp_2$}[\atleft, \solidarrow]
\mor{I}{I1}{$\bfp_3$}[\atright, \solidarrow]
\enddc
\EqTag{CDeds}
\end{gathered}
\end{equation}
is a {\deffont commutative diagram of EDS}. For future reference, we remark that the commutative diagram \EqRef{CDsubs}  also implies that
\begin{equation}
	\ker \bfp_{2*} = \bfp_{1*}(\ker \bfp_{3*}).
\EqTag{kerp}
\end{equation}

\begin{Lemma}
\StTag{IntExt4}
	Let $G$ be  a symmetry group of an EDS $\CalI$ on a manifold $M$ which acts regularly on $M$.
	Let $H\subset G$  be a subgroup of  $G$ which also acts regularly on $M$.  Then the orbit mapping
\begin{equation}
\bfp : M/H \to  M/G \quad\text{defined by} \quad\bfp(Hx) = Gx
\EqTag{IntExt14}
\end{equation}
	is a surjective submersion  which gives rise to the following commutative diagram of EDS
\begin{equation}
\begin{gathered}
\begindc{\commdiag}[3]
\obj(-25, 0)[I]{$\CalI$}
\obj(0, 0)[H]{$\CalI /H$}
\obj(0, -17)[I1]{$\CalI /G\ ,$}
\mor{I}{H}{$\bfq_H$}[\atleft, \solidarrow]
\mor{H}{I1}{$\bfp$}[\atleft, \solidarrow]
\mor{I}{I1}{$\bfq_{G}$}[\atright, \solidarrow]
\enddc
\EqTag{CDfund}
\end{gathered}
\end{equation}
that is, $(\CalI/H)/\bfp  = \CalI/G$.
\end{Lemma}

\begin{proof} The given hypotheses leads immediately to the commutative diagram
\begin{equation}
\begin{gathered}
\begindc{\commdiag}[3]
\obj(-25, 0)[I]{$M$}
\obj(0, 0)[H]{$M /H$}
\obj(0, -17)[I1]{$M /G\ $}
\mor{I}{H}{$\bfq_H$}[\atleft, \solidarrow]
\mor{H}{I1}{$\bfp$ , }[\atleft, \solidarrow]
\mor{I}{I1}{$\bfq_{G}$}[\atright, \solidarrow]
\enddc
\EqTag{orbitmap}
\end{gathered}
\end{equation}		
	where $\bfq_G$ and $\bfq_H$ are surjective submersions. 
	Theorem \StRef{Gcomd}  then implies  that $\bfp$ is a surjective submersion and that the diagram \EqRef{CDfund} commutes.
\end{proof}

\noindent
{\bf Proof of part [i] of Theorem A}: Apply Corollary \StRef{IntExt4} to each of the two cases $H \subset G_1$ and $H\subset G_2$ and this immediately produces the commutative diagram \EqRef{Intro21}.

\begin{Remark}
\StTag{NReduction}
	 The subgroup $H$ in Lemma  \StRef{IntExt4}  acts regularly if and only if it a closed subgroup
	(see  \cite{dieudonne:1970},  (16.10.3)  on page 58). 	If the action of $G$ in   Lemma \StRef{IntExt4} 
	is free and if $H$ is a closed normal subgroup
	then  the natural action of $G/H$  on $M/H$  is both free and regular.
	If $G$ is a symmetry group of $\CalI$, 
	then $G/H$ is a symmetry group of $\CalI/H$. By Corollary  \StRef{IntExt4}, $$\CalI/G=(\CalI/H)/(G/H), $$ 
	where the projection map $\bfp$ in \EqRef{CDfund} is now  the quotient map for the action of $G/H$ on $M/H$. 
	Finally, if $G$ acts transversely to $\CalI$ then $G/K$ acts transversely to $\CalI/H$.
\end{Remark}

\subsection{Subgroup Reductions and Integrable Extensions}

	The proof of the second statement in Theorem  A  follows directly from the following general result.

\begin{Theorem}
\StTag{IntExt3}
	Let $G$ be  a symmetry group of an EDS $\CalI$ on a manifold $M$. Assume $G$ acts regularly on $M$ and transversely to $\CalI$, and
	let $H\subset G$  be a  closed Lie subgroup of  $G$.  Then the orbit mapping $\bfp$ in
	\EqRef{IntExt14} defines $\CalI/H$ as an integrable extension of  $\CalI/G$.
\end{Theorem}	

\begin{proof} 
	Since $G$ acts transversely to $\CalI$,  $H$ also acts transversely to $\CalI$.
	These transversality conditions  ensure that the quotient bundles $I^1/G$ and $I^1/H$ 
	are constant rank and consequently  we
	may construct a  constant rank bundle $\barK \subset  \Lambda^1(M/H)$ such that 
\begin{equation}
	I^1/H =   \barK \oplus \bfp^*(I^1/G).
\EqTag{barK}
\end{equation}
	We prove that  $\barK$ satisfies the  3 conditions of  definition  \EqRef{IntExtDef}. 
	From \EqRef{barK}  and   \EqRef{rankI1} it follows that 
$$	
	\rank( \barK) = \rank(I^1/H) - \rank(I^1/G)  = \rank(\bfGamma_G) -\rank(\bfGamma_H) = \dim(M/H) - \dim(M/G).
$$
	so that condition {\bf [i]} of definition  \EqRef{IntExtDef} is satisfied. 

	The second condition in   \EqRef{IntExtDef} is  $\ann(\barK) \cap \ker(\bfp_*) = 0$ which we now re-write as a condition on the 
	distributions $\bfGamma_G$ and $\bfGamma_H$. Set  $K = \bfq_{H}^*(\barK)$ so that $K_{H, \semibasic}  = K$.   
	Then   the application of  \EqRef{kerp} to the diagram \EqRef{CDfund} yields
\begin{equation}
	\ker( \bfp_* )= \bfq_{H*} (\bfGamma_G).
\EqTag{IntExt9}
\end{equation}
	Therefore, since $\bfq_{H*}( \ann (K )) = \ann(\barK)$,  condition  {\bf [ii]} in   \EqRef{IntExtDef} is equivalent to 
\begin{equation*}
\bfq_{H*}( \ann (K ) \cap \bfGamma_G) = 0\qquad {\rm or} \qquad 
\bfGamma_G  \cap \ann(K) \subset \bfGamma_H\, . 
\end{equation*}
We will prove the slightly stronger condition
\begin{equation}
\bfGamma_G  \cap \ann(K)= \bfGamma_H,
\EqTag{newAdK1}
\end{equation}
which implies condition {\bf [ii]}.

	To  prove \EqRef{newAdK1} we use  the commutativity of  \EqRef{CDfund}, together with   \EqRef{Asb},  to write
\begin{equation}	
	I^1_{H ,\semibasic} = K \oplus I^1_{G, \semibasic}.
\EqTag{IntExt10}
\end{equation}
	Transversality (see \EqRef{rankI1})  then gives 
\begin{equation*}
	\rank(K)  =  \rank  I^1_{H ,\semibasic}  - \rank I^1_{G, \semibasic}  =  \rank \bfGamma_G  -  \rank \bfGamma_H.
\end{equation*}
	which in turn, implies that 
\begin{equation}
\rank( \Lambda^1(M)_{H,\semibasic}) =  \dim(M) -\rank \bfGamma_H  = \rank( \Lambda^1(M)_{G, \semibasic})  +\rank (K).
\end{equation}
	It follows from \EqRef{IntExt10} that  $ \Lambda^1(M)_{G, \semibasic} \cap K$  =0  and consequently
\begin{equation}
	\Lambda^1(M)_{H,\semibasic}  =   \Lambda^1(M)_{G, \semibasic} \oplus K.
\EqTag{LambdaH}
\end{equation}
	The annihilator of this equation yields   \EqRef{newAdK1} : 
\begin{equation}
	\bfGamma_H = \ann( \Lambda^1(M)_{H,\semibasic}  )  = \ann( \Lambda^1(M)_{G, \semibasic}) \cap  \ann(K)  
	= \bfGamma_G  \cap \ann(K).
\end{equation}

	Condition {\bf [iii]} in   \EqRef{IntExtDef} is $\CalI/H  =    \langle \CalS(\barK)  \cup \bfp^*(\CalI/G) \rangle_\text{alg}. $
	To prove this we shall    use  \EqRef{IntExt10}  and \EqRef{LambdaH}
	to  refine the local coframe   \EqRef{IntExt15} used in  the proof of Theorem \StRef{IntEntCor1}. 
	Choose complementary bundles $W$ and $L$ such 
\begin{equation}
I^1_{H, \semibasic} \oplus L =  I^1
	, \quad I^1_{G, \semibasic} \oplus W
	=  \Lambda^1(M)_{G ,\semibasic}\quad {\rm in \ which \ case } \quad	I^1_{G, \semibasic} \oplus J = I^1
\EqTag{IntExt11}
\end{equation}
where $J=K\oplus L $. We also find
\begin{equation}
	 \Lambda^1(M)_{H,\semibasic}  =   \Lambda^1(M)_{G, \semibasic}  \oplus K  = I^1_{G, \semibasic} \oplus W   \oplus K.
\EqTag{IntExt13}
\end{equation}
	The subbundles
	$J$, $W$ introduced in that proof (see \EqRef{JW}) will be taken to be  the same as those defined here.
	
	Since $I^1_{H, \semibasic} \oplus L =  I^1 = I^1_{G, \semibasic} \oplus J$,
	we can deduce from \EqRef{IntExt10} that
\begin{equation*}
		  I^1_{G, \semibasic} \oplus K \oplus L \oplus W 
		= \Lambda^1(M).
\end{equation*}
	Accordingly, about each point of  $M$,  there is an open set $U$ and  forms
	$\, \theta^u,\,   \theta^r_{H} , \theta^a_{G}, \,  \tau^\alpha_{G}$  with
\begin{gather*}
	L|_U = \spn\{\, \theta^u\,\}, \quad   K|_U = \spn\{\, \theta^r_{H} = \bfq_H^*\tilde \theta^r \, \} , 	\quad
	I^1_{G, \semibasic}|_U =  \spn\{\, \theta^a_{G} = \bfq_G^*\bar \theta^a \, \}, \quad  
\\[2\jot]
	\CalI|_{U} = \langle\, \theta^u,\,   \theta^r_{H} , \theta^a_{G}, \,  \tau^\alpha_{G} \, \rangle_\text{alg}  ,
	\quad 
	\barK|_{ \bfq_H(U)} =   \spn\{\, \tilde \theta^r \,\}
	\quad\text{and}\quad 
	(\CalI/G)|_{\bfq_G(U)}  = \langle\   \bar \theta^a , \,   \bar \tau^\alpha  \, \rangle_{\text{alg}}, 
\end{gather*}
	where $\tau^\alpha_{G} = \bfq_G(\bar \tau^\alpha) \in  \Lambda^*(W_{G, \semibasic} )$.
	Put $\tilde \theta^a = \bfp^*(\bar \theta^a)$ and   $\tilde \tau^\alpha = \bfp^*(\bar \tau^\alpha)$. Then, 
	in terms of these forms, we have
\begin{gather}	
	\langle\bfp^*(\CalI/G)\rangle _{\text{alg}} |_{\bfq_H( U)}  =  \langle\   \tilde \theta^a , \,   \tilde \tau^\alpha  \, \rangle_{\text{alg}} 
	\quad\text{and}\quad  
	(\CalI/H)|_{\bfq_H( U) }  =  \langle\  \tilde \theta^r  ,  \tilde \theta^a, \,   \tilde \tau^\alpha  \, \rangle_{\text{alg}} 
\end{gather}
	and  condition {\bf [iii]} in  \EqRef{IntExtDef} is proved.
\end{proof}

\begin{Corollary}
\StTag{IntEntCor2}
Under the hypotheses of Theorem 	\StRef{IntExt3}, the diagram \EqRef{CDfund} is a commutative diagram of integrable extensions.
\end{Corollary}

\noindent
{\bf Proof of part [ii] of Theorem A}: 
	The transversality hypothesis in part  {\bf [ii]} of  Theorem A allows to apply Theorem
	\StRef{IntEntCor1} and   Corollary \StRef{IntEntCor2} to the  two cases $H\subset G_1$ and $H \subset G_2$.

\medskip
	We remark that $\CalI$ is a rather special type of integrable extension of $\CalI/G$ and perhaps deserving of the 
	designation {\deffont integrable extension of  Lie type}.
	Indeed,   it is shown in \cite{anderson-fels:2005a},
	in the special case where $\CalI$ is Pfaffian and the action of $G$ is free, that the integral manifolds for $\CalI$ can be 
	constructed from the integral manifolds of $\CalI/G$ by the integration of ODE's of Lie type.

	It  should be emphasized that, unless  $H$ is a normal subgroup of G, 
	then {\it  the orbit projection maps $\bfp_1$ and $\bfp_2$  in  \EqRef{Intro21}.
 are  not  group quotients.}
	Note that  the fibre dimensions for the maps $\bfp_1$ and $\bfp_2$ 
	defining the B\"acklund transformation  \EqRef{Intro21} are given by the difference of orbit dimensions
	$\dim \CalO_{G_1} -\dim \CalO_H$ and $\dim \CalO_{G_2} - \dim \CalO_H$, respectively.  Thus, 
	in the particular case  of  free actions, 
	Theorem A  can be used to construct  a B\"acklund transformation  
	with 1-dimensional fibres only when  $G_1$ and $G_2$ have a common subgroup of  co-dimension 1.


\section{Darboux Integrable Differential Systems}
	In this section we review the definition of Darboux integrable differential systems and 
	we state a new result which simplifies the  verification of this definition. We begin with the definition of  a decomposable differential system and the  notion of a Darboux pair as found in \cite{anderson-fels-vassiliou:2009a}.

\begin{Definition}
\StTag{DI1}	
	An exterior differential system  $\CalI$ on $M$ is   {\deffont decomposable of type $[p,\rho]$}, where
	$p,\rho \geq 2$, if about each point $x \in M$ there is a local coframe
\begin{equation}
	\tilde \theta^1,\ \ldots,\ \tilde \theta^r,\
	\hsigma^1,\ \dots,\ \hsigma^{p}, \
	\csigma^1,\ \dots,\ \csigma^{\rho},
\EqTag{Intro4}
\end{equation}
	such that $\CalI$ is algebraically generated by  1-forms and 2-forms
\begin{equation}
	\CalI = \langle \, \tilde \theta^1, \ \dots, \ \tilde \theta^r, \ \hOmega^1,\ \dots,\ \hOmega^{s},\ \cOmega^1,\  \dots, \ \cOmega^{\tau} \, \rangle_\text{\rm alg},
\EqTag{Intro5}	
\end{equation}
	where $s, \tau \geq  1$,
	$\hOmega^a \in \Omega^2(\hsigma^1,\dots,\hsigma^{p})$, and
	$\cOmega^\alpha \in \Omega^2(\csigma^1,\dots,\csigma^{\rho})$.
\end{Definition}	
\par
	Equation \EqRef{Intro5} implies that the 1-forms $\tilde \theta^e$ satisfy structure equations of the form
\begin{equation}
	d \tilde \theta^e \equiv A^e_{ab} \ \hsigma^a \wedge  \hsigma^b +  B^e_{\alpha\beta} \ \csigma^\alpha \wedge  \csigma^\beta \quad \mod \{\, \tilde \theta ^e\, \} 
\EqTag{DecStEq}
\end{equation}
	where 
\begin{equation}
	\text{span} \{A^e_{ab} \ \hsigma^a \wedge  \hsigma^b \} \subset  \text{span} \{\hOmega^a \} 	
	\quad \text{and}\quad
	\text{span} \{ B^e_{\alpha\beta} \ \csigma^\alpha \wedge  \csigma^\beta \} \subset \text{span}\{ \cOmega^a\}.	
\end{equation}
	In the special case that $\CalI$ is a  Pfaffian system, these inclusions then  become equalities.

	In particular, any class $r$ hyperbolic differential system, as
	defined in   \cite{bryant-griffiths-hsu:1995a},  is a decomposable differential system  of type [2, 2].  
	In the special 
	case of  a class $r$ hyperbolic differential system the decomposition is essentially unique but  
	this uniqueness property is not true for general decomposable systems.

\begin{Definition}  Let $\CalI$ be a decomposable differential system. The bundles $\hV, \cV \subset T^*M$ defined by
\begin{equation}
\hV = \text{span} \{\,\tilde \theta^e,\, \hsigma^a  \}
	\quad\text{and}\quad
	\cV = \text{span} \{\, \tilde \theta^e,\, \csigma^\alpha  \}
\EqTag{SPS}
\end{equation}
are called the {\deffont associated singular Pfaffian systems} with respect to the decomposition \EqRef{Intro5}. The 
	 the differential systems  generated by
\begin{equation}
	\hCalV =\langle \,\tilde \theta^e,\, \hsigma^a, \,\cOmega^\alpha\, \rangle_\text{\rm diff}
	\quad\text{and}\quad
	\cCalV = \langle \, \tilde \theta^e,\, \csigma^\alpha, \,\hOmega^a\, \rangle_\text{\rm diff}
\EqTag{Intro6}
\end{equation}
	are called the associated {\deffont singular differential systems} for $\CalI$. 
\end{Definition}
	
	Note that
\begin{equation}
	I^1 = \hV \cap \cV  \quad\text{and}\quad T^*M  =  \hV  + \cV.
\EqTag{hVcapcV}
\end{equation}
	For further information on the relationship between a decomposable differential system and its singular systems 
	see Theorem 2.6 of  \cite{anderson-fels-vassiliou:2009a}.

	 The  characteristic  systems  defined in  
	\cite{bryant-griffiths-hsu:1995a} for a class $r$ hyperbolic differential system  coincide with the singular Pfaffian systems \EqRef{SPS}.
	Two  decomposable differential systems $\CalE$ and $\CalI$, defined on manifolds  $N$ and $M$, and with singular Pfaffian systems
	$\{\hZ, \cZ\}$ and $\{\hV, \cV\}$ are defined to be equivalent if there is a diffeomorphism $\phi :N \to M$ such that $\phi^*(\CalI) = \CalE$,
	$\phi^*(\hV) = \hZ$ and $\phi^*(\cV) = \cZ$.  Under these circumstances, it follows that  $\phi^*(\hCalV) = \hCalZ$ and 
	$\phi^*(\cCalV) = \cCalZ$. 

	The definition of a Darboux integrable differential system is given in terms of its singular Pfaffian systems.
\begin{Definition}
\StTag{Intro2}
	 A pair of Pfaffian systems $\hV$ and $\cV$  on a manifold $M$ define a
	{\deffont Darboux pair} if 
\par
\smallskip
\noindent
{\bf[i]}
\vskip -30pt
\begin{equation}\hV  + \cV^{\infty} = \cTM
	\quad\text{and}\quad
        \cV  + \hV^{\infty} = \cTM ,\quad \text{and}
\EqTag{Intro7}
\end{equation}
\par
\smallskip
\noindent
{\bf[ii]}
\vskip -30pt
\begin{equation}
 \hV^{\infty} \cap \cV^{\infty} = \{\,0\, \}.
\EqTag{Intro8}
\end{equation}
	A decomposable differential system $\CalI$ is {\deffont Darboux integrable}  
	if its singular Pfaffian systems  \EqRef{SPS} define a Darboux pair.	
\end{Definition}

	This definition of Darboux integrability is slightly more general than that given in \cite{anderson-fels-vassiliou:2009a}, 
	where it was assumed that the singular systems $\hCalV$ and $\cCalV$ are Pfaffian systems.
	It is a simple matter to argue (see Appendix C, or Section 2.2 in \cite{anderson-fels-vassiliou:2009a}) 
	that if $\CalI$ is Darboux integrable, then 
	about each point of $M$ there exists a local coframe 
	$\{\, \bftheta,\,  \bfheta,\,  \bfceta,\,  \bfhsigma, \bfcsigma \ \}$ such that
\begin{equation}
	I^1  =   \spn \{\ \bftheta,\,  \bfheta,\,  \bfceta \, \} \quad\text{and}\quad
\begin{aligned}
	\hV =  \spn \{\ \bftheta,\,  \bfheta,\,  \bfceta,\,  \bfhsigma \ \}, & \quad \hV^\infty = \spn \{\ \bfheta,\,  \bfhsigma \ \}, 
\\
	\cV=  \spn \{\ \bftheta,\,  \bfheta,\, \bfceta,\,  \bfcsigma \ \}, & \quad \cV^\infty = \spn \{ \ \bfceta,\,  \bfcsigma  \ \} .
\end{aligned}
\EqTag{ZeroC}
\end{equation}
	Such coframes are said to be {\deffont 0-adapted coframes with respect to the Darboux pair $\{\hV, \cV\}$}. Note that
\begin{equation}
	\hV^\infty \cap \cV =  \spn\{\, \bfheta\, \} \quad   \text{and}  \quad \cV^\infty \cap \hV =  \spn\{\, \bfceta \, \}. 
\EqTag{ZeroC1}
\end{equation} 

	The  structure equations for the forms $\bfhsigma$, $\bfcsigma$ (arising from the fact that  $\hV^\infty$ and 
	$\cV^\infty$ are completely integrable Pfaffian systems) 
	 immediately imply that the  differential generators for $\hCalV$ and $\cCalV$ in \EqRef{Intro6} are now 
	algebraic generators, that is,
\begin{equation}
	\hCalV =\langle \,\tilde \theta^e,\, \hsigma^a, \,\cOmega^\alpha\, \rangle_\text{\rm alg}
	\quad\text{and}\quad
	\cCalV = \langle \, \tilde \theta^e,\, \csigma^\alpha, \,\hOmega^a\, \rangle_\text{\rm alg}.
\EqTag{Intro6DI}
\end{equation}
Equations \EqRef{Intro6DI} are also transparent from the structure equations in Theorem 2.9 in \cite{anderson-fels-vassiliou:2009a} for a 1-adapted coframe. 

\begin{Definition} The Darboux invariants for the Darboux pair $\{\hV, \cV\}$ or for a Darboux integrable differential system  $\CalI$ 
	are the first integrals for $\hV$ or $\cV$,
	that is,   $C^\infty(M)$ functions  $f$ such that $d f \in \CalS( \hV)$ or  $df \in  \CalS(\cV)$. 
\end{Definition}

  We shall make repeated use of the following.
\begin{Lemma}\StTag{dfL}  
	If  $f$ is a locally defined, real-valued function on $M$ such that 
\begin{equation*} 
	df  \in \CalS(\hV) = \spn \{\,\bftheta, \, \bfheta, \, \bfceta,\, \bfhsigma\, \}
\end{equation*} 
	then $df \in  \CalS(\hV^\infty)$ and hence  $df \in  \spn \{ \, \bfheta,\,  \bfhsigma\, \}$. 
	Likewise, if $df \in  \CalS(\cV) = \spn \{ \, \bftheta, \,  \bfheta, \, \bfceta, \, \bfcsigma\, \}$ then
	$df \in  \CalS(\cV^\infty)$  and $df \in \spn \{ \,  \bfceta,\, \bfcsigma \, \}$.
\end{Lemma}

\begin{proof} 
	If $df \in \CalS(\hV)$,
	then,  by  the definition of the terminal derived Pfaffian system,  $df \in \CalS(\hV^\infty)$, and the lemma follows from equation \EqRef{ZeroC}.
\end{proof}

	The number of (functionally) independent Darboux  invariants is therefore given by the sum of the ranks  of  the completely 
	integrable systems  $\hV^\infty$ and $\cV^\infty$.  

	The next theorem greatly simplifies the task of  verifying that a decomposable differential system is Darboux integrable.
	This simplification will be used in Section 5.1 to prove  Theorem B.

\begin{Theorem} 
\StTag{DPD}
	Let $\CalI$ be decomposable differential system with singular Pfaffian systems  $\hV$ and $\cV$ 
	and suppose that $(I^1)^\infty = (\hV \cap \cV)^\infty = 0$.  If  $\hV$ and $\cV$ satisfy  conditions {\bf [i]} in the 
	definition of a Darboux pair, then condition {\bf[ii]} is automatically satisfied and $\CalI$ is Darboux integrable.
\end{Theorem}
	The proof of this theorem is given in Appendix C.

\medskip

	The differential systems which we consider in our examples and applications  arise from partial differential equations 
	and are linear differential systems with independence condition (see \cite{bryant-chern-gardner-griffiths-goldschmidt:1991a}, Chapter IV and especially page 111).  This leads to the following definition.

\begin{Definition}
	\StTag{DCwithIC} 
A linear differential system   $\CalI$  with
	independence  condition $\omega^1\wedge \dots  \wedge \omega^m  \neq  0$ is Darboux integrable if  the following conditions 
	hold.

\smallskip
\noindent
{\bf [i]}
There  exists local coframes
\begin{equation*}
	\tilde \theta^1,\ldots, \tilde \theta^r,\
	\homega^1,\dots, \homega^{m_1}, \htau^1, \dots, \htau^{p_1}, \
	\comega^1,\dots, \comega^{m_2}, \ctau^1, \dots, \ctau^{\rho_2}
\end{equation*}
	with  $m_1 +p_1 \geq 2$, $m_2 + \rho_2\geq 2$,
\begin{equation}
	\omega^1 \wedge \dots \wedge \omega^m
	= \homega^1 \wedge \dots \wedge \homega^{m_1} \wedge \comega^1 \wedge\dots \wedge \comega^{m_2},
\EqTag{omegaIC}
\end{equation}
	where  \EqRef{Intro5} holds,  where the 2-forms  $\hOmega^a$ and $\cOmega^\alpha$ now  assume the form
\begin{equation}
	\hOmega^a = \hL^a_{bc}\,\htau^b \wedge \homega^c 
	\quad\text{and} \quad 
	\cOmega^\alpha = \cL^\alpha_{\beta\gamma}\,\ctau^\beta\wedge \comega^\gamma 
\EqTag{Omegatauomega}
\end{equation}
	and where the  structure equations  are 
\begin{equation}
	d \tilde \theta^e = A^e_{ab}\,\htau^a \wedge \homega^b +
	B^e_{\alpha \beta}\, \ctau^\alpha \wedge \comega^\beta \quad \mod \tilde \theta^e
\EqTag{DecStrEq2}.
\end{equation}

\smallskip
\noindent
{\bf [ii]} The singular Pfaffian systems $ \hV = \text{span}\{\tilde \theta^e,\, \htau^a, \,\homega^i\}$
	and $ \cV =  \text{span}\{\tilde \theta^e, \,\ctau^a, \, \comega^i\}$
define a Darboux pair.

\smallskip
\noindent
{\bf [iii]}  $\hJ = \text{span}\{\homega^1,\dots, \homega^{m_1}\}$ is an integrable sub-system of  $\hV^\infty$ and 
$\cJ = \text{span}\{\comega^1,\dots, \comega^{m_2}\}$ is an integrable sub-system of  $\cV^\infty$

\end{Definition}
		
\section{Integrable Extensions of Darboux Integrable Differential Systems}
	In the previous section we defined Darboux integrable differential systems $\CalI$  as decomposable differential systems 
	whose associated singular systems satisfy the conditions of  Definition \StRef{Intro2}.   Let $\bfp : (\CalE, N) \to (\CalI, M)$ be 
	an integrable extension.  In Section  5.1 we show how the pullbacks  by $\bfp$ of the  singular systems  for $\CalI$ and an admissible 
	subbundle for $\CalE$ can be used to construct singular systems for $\CalE$  which satisfy the conditions of Darboux integrability.
    Then, in Section 5.2, 
	we  introduce a special  class of integrable extensions for Darboux integrable systems
	which we call {\deffont maximally compatible integrable extensions}.  
	These integrable extensions  are maximal in the sense that  
	the number of (functionally independent) Darboux invariants for each singular Pfaffian system is as large as possible. 
	Such extensions arise very naturally in the construction of 
	Darboux integrable systems by symmetry reduction in Section 6, in the proof of the uniqueness of the local quotient
	representation in Section 8, and in the analysis Darboux integrable Monge-Amp\`ere equations in  
	 \cite{anderson-fels:2011b}.

\subsection{Integrable Extensions are Darboux Integrable}
The precise formulation of Theorem B is the following.
\begin{Theorem}
\StTag{DIext} 
	Let $\bfp : (\CalE, N) \to (\CalI, M)$ be an integrable extension with  $J$ an admissible subbundle of $T^*N$ for $(\CalE, \CalI)$.

\smallskip
\noindent
{\bf [i]}  If  $\CalI$ is decomposable of type $[p, \rho]$ with singular Pfaffian systems $\hV$ and $\cV$, 
	then $\CalE$ is decomposable of type $[p, \rho]$ with singular Pfaffian systems  
\begin{equation}
	\hZ = J \oplus \bfp^* ( \hV) \quad\text{and}\quad \cZ = J \oplus \bfp^* (\check V).
\EqTag{singWJV}
\end{equation}
The corresponding singular differential systems satisfy
\begin{equation}
	\hCalZ  = \langle \CalS(J) + \bfp^*( \hCalV) \rangle_{\text{\rm alg}} \quad\text{and}\quad 
	\cCalZ  =  \langle \CalS( J) + \bfp^*( \cCalV ) \rangle_{\text{\rm alg}} \ . 
 \EqTag{IEWs1}
\end{equation}

\smallskip
\noindent
{\bf [ii]} The systems $\hZ$, $\cZ$, $\hCalZ$ and $\cCalZ$ are integrable extensions of  $\hV$, $\cV$, $\hCalV$ and $\cCalV$
respectively.

\smallskip
\noindent
{\bf [iii]}  If $\CalI$ is Darboux integrable and $(E^1)^\infty =0$,
	then $\CalE$ is Darboux integrable.\footnote{If $(E^1)^\infty \neq 0$ ,  
	then one can always let $P$   be an integral  manifold $N'$  of $(E^1)^\infty =0$,  let 
	$\CalF = \CalE_N$  and replace
	$(\CalE, N)$ in the statement of the theorem  by  $(\CalF, P)$, where now $\CalF^\infty$ = 0 }
\end{Theorem}
\begin{proof} 
	Let $\{\tilde\theta^e,\, \hsigma ^a,\, \csigma^\alpha \}$ be a local coframe on $U_0\subset M$ and let
	$\hOmega^c$, $\cOmega^\gamma$ be  2-forms  on $U_0$ such that (see equation \EqRef{Intro5})
\begin{equation*}
	\CalI|_U = \langle\, \tilde\theta ^e,\,  \hOmega^c ,\, \cOmega^\gamma\, \rangle_{ \text{alg}}.
\end{equation*}
	The structure equations for the 1-forms  $\tilde\theta ^e$ are (see \EqRef{DecStEq})
\begin{equation}
	d\tilde\theta^e = A^e_c\, \hOmega^c  + B^e_\gamma\,\cOmega ^\gamma  \quad \mod \{\, \tilde\theta^e \, \}, \EqTag{oeqs}
\end{equation}
	and the singular systems are (see \EqRef{Intro6})
\begin{equation*}
	\hCalV|_U = \langle\, \tilde\theta^ e, \hsigma^a, \cOmega^\gamma\, \rangle_{ \text{diff}} \quad\text{and}\quad 
	\cCalV|_U = \langle\, \tilde\theta^e, \csigma^\alpha, \hOmega^c \,\rangle_{ \text{diff}}  \ .
\end{equation*}	

	Now choose an open set $U \subset \bfp^{-1}(U_0)$ and a local basis of sections $\xi^u$ for $J$.  Allowing for a slight abuse of notation,
	we have that the 1-forms $\{\tilde\theta^e,\,  \xi^u,\,  \hsigma ^a,\, \csigma^\alpha\}$ is a local coframe on $U$.
	By the integrable extension property \EqRef{IntExt4}  the structure equations  for $\CalE$ are \EqRef{oeqs} and 
\begin{equation}
	d \xi ^u = F^u_v \wedge \xi^ v + G^u_e \wedge \tilde\theta^e + H^u _c \, \hOmega^ c + L^u_\gamma \,\cOmega^\gamma,
\EqTag{exteqs}
\end{equation}
	where  the $F^u_v$, $G^u_e$ are 1-forms on $U$  and  $H^u_b$, $L^u_\beta \in C^\infty(U)$. On account of the
	structure equations \EqRef{oeqs} and \EqRef{exteqs}
	we have that
\begin{equation}
	\CalE|_U= \langle \, \tilde\theta^e,\, \xi^u,\, \hOmega^b ,\, \cOmega^ \beta \,\rangle_{\text{alg}}
\end{equation}
	and $\CalE$ is clearly decomposable of type $[p, \rho]$ with singular  Pfaffian systems 
\begin{equation}
	\hZ|_U = \{\, \tilde\theta^e,\,  \xi^u,\, \hsigma^a\, \}
	\quad\text{and}\quad
	\cZ|_U = \{\, \tilde\theta^e ,\,   \xi^u,\, \csigma^\alpha\, \} 
\EqTag{dhW}
\end{equation}
and singular differential systems
\begin{equation}
	\hCalZ|_U = \langle \, \tilde\theta^e, \xi^u, \hsigma ^a, \cOmega^\gamma\,\rangle_{\text{diff}} 
	\quad\text{and}\quad
	\cCalZ|_U = \langle \,\tilde\theta^e, \xi^u, \csigma ^\alpha, \hOmega^c \,\rangle_{\text{diff}}.
\EqTag{dhCalW}
\end{equation}
Equation  \EqRef{dhW} proves  \EqRef{singWJV}.   The systems $\langle \bfp^*\hCalV \rangle_{\text{alg}} $ and $\langle \bfp^*\cCalV \rangle_{\text{alg}} $ are differentially closed. The combination of this with equation \EqRef{exteqs} shows that the systems in equation \EqRef{dhCalW} can be
written as \EqRef{IEWs1} (that is the systems in  \EqRef{IEWs1} are differentially closed). 

	Equations  \EqRef{SPS} and \EqRef{dhW} and also  \EqRef{Intro6} and \EqRef{dhCalW}  immediately prove  part {\bf [ii]}.

	To prove part {\bf [iii]}, we first note that  $\bfp^* (\cV^\infty )$ is a completely integrable Pfaffian system contained in $\cZ$ 
	and therefore 
\begin{equation}
	\bfp^* (\cV^\infty )\subset \cZ^{\infty} .
\EqTag{pVinftyinZ}	
\end{equation}
	Suppose  now that $\CalI$ is Darboux integrable. Then,  since $\hV + \cV^{\infty} = T^*M$,   we deduce  from  \EqRef{pVinftyinZ}
	that
\begin{equation*}
	\hZ + \cZ^\infty  \supset   J + \bfp^*(\hV) +\bfp^* (\cV^\infty) = J \oplus \bfp^*(\hV + \cV^\infty) =  J \oplus  \bfp^*(T^*M) = T^*N.
\label{Wsum}
\end{equation*}
	This proves that the first part of condition {\bf [i]} in equation \EqRef{Intro7} is satisfied. The second equation  is similarly proved.  
	Theorem \StRef{DPD} then shows $\hZ$ and $\cZ$ define a Darboux pair, and so $\CalE$ is Darboux integrable.
\end{proof}

	We remark that if $\bfp\:  (\CalE, N) \to (\CalI, M)$ is  an integrable extension of decomposable systems $\CalE$ and $\CalI$ 
	then  the relations \EqRef{singWJV} between the singular systems for $\CalE$ and $\CalI$ need not automatically hold.  
	This is because  $\CalE$ may be decomposable with respect to different choices of singular systems.

\subsection{ A special class of integrable extensions for Darboux integrable differential systems}

	Motivated by Theorem  \StRef{DIext} and the remark made in the last paragraph in Section 5.1 we now give 
	two  definitions which will play an important  role in all that follows in Sections 7 and 8. In these sections  we shall 
	relate the fundamental invariants for a pair a  Darboux integrable systems $\CalE$ and  $\CalI$ in the case 
	when $\CalE$ is an integrable extension of $\CalI$. 
	The first of  these  definitions formalizes  the relations  in  \EqRef{singWJV} and implies that  the number of Darboux invariants
	for the singular Pfaffian systems $\hZ$ and $\cZ$ for the extension $\CalE$ are bounded by (see \IntExtInftymodp)
\begin{align*}
	\rank(\bfp^*(\hV^\infty))  & \leq  \rank( \hZ^\infty)  \leq \rank (\bfp^*(\hV^\infty))    + \rank (J) \quad \text{and}\\
	\rank(\bfp^*(\cV^\infty))  & \leq   \rank( \cZ^\infty)  \leq \rank (\bfp^*(\cV^\infty))    + \rank (J). 
\end{align*}
	The case that is the most important  to us is the case where the number of Darboux invariants for  $\CalE$
	exceeds the number of Darboux invariants for  $\CalI$ by the maximal number possible, that is, the case 
	where 
\begin{equation}
	\rank( \hZ^\infty)  =  \rank (\bfp^*(\hV^\infty))    + \rank (J)   \quad \text{and}\quad
	\rank( \hZ^\infty)  =  \rank (\bfp^*(\hV^\infty))    + \rank (J) .
\EqTag{rankZinfty}
\end{equation}
	These rank equalities will be implied by part {\bf [ii]} of the following definition.

\begin{Definition}
\StTag{DarbComp}
	Let   $\bfp\:  (\CalE, N) \to (\CalI, M)$ be an integrable extension of decomposable systems $\CalE$ and $\CalI$ with
	singular Pfaffian systems $\{\hZ, \cZ\}$ and  $\{\hV, \cV\}$ respectively.

\smallskip
\noindent
	{\bf [i]}  
	The extension $(\CalE, \CalI)$ is said to be {\deffont compatible with respect to 
	the singular Pfaffian systems} $\{\hZ, \cZ\}$ and  $\{\hV, \cV\}$
	if  
\begin{equation}
	 \bfp^* ( \hV) \subset \hZ  \quad\text{and}\quad \bfp^* (\cV) \subset \cZ .
\EqTag{DarbAdmis0}
\end{equation}

\smallskip
\noindent
	{\bf[ii]}
	The integrable extension $(\CalE, \CalI)$  is  called  {\deffont maximally compatible  with respect to 
	the singular Pfaffian systems} $\{\hZ, \cZ\}$ and  $\{\hV, \cV\}$ if there exists 
	an admissible subbundle $\hJ \subset \hZ^\infty\cap \cZ$  for $(\CalE, \CalI)$  such that 
\begin{equation}
	\cZ =  \hJ \oplus \bfp^* ( \cV)  \quad  \text{and} \quad
	\hZ^\infty = \hJ \oplus \bfp^*(\hV^\infty),  
\EqTag{DarbAdmis1}
\end{equation}
	and, similarly, an admissible 
	subbundle $\cJ \subset \cZ^\infty\cap \hZ$   for $(\CalE, \CalI)$ such that
\begin{equation}
	\hZ =    \cJ  \oplus \bfp^* ( \hV) \quad \text{and} \quad  \cZ^\infty = \cJ \oplus \bfp^*(\cV^\infty).
\EqTag{DarbAdmis2}
\end{equation}
\end{Definition}

	A couple of simple observations regarding these definitions are in order. First, suppose that 
	 $(\CalE, \CalI)$ is a  compatible integrable extension and that $J$ is an 
	admissible subbundle.  The inclusions \EqRef{DarbAdmis0}, the definitions of $\hZ$ and $\cZ$, 
	and the fact $E^1 = J \oplus  \bfp^*(\hV)$
	show that  $J \oplus \bfp^*(\hV) \subset \hZ$ and $J \oplus \bfp^*(\cV) \subset \cZ$.
	Moreover, on account of \EqRef{IntExt3},  \EqRef{IntExt5}  and \EqRef{hVcapcV} 
	we find that 
\begin{gather*}
	\hZ + \cZ = T^*N = J  \oplus \bfp^*(T^*M) = ( J \oplus \bfp^*(\hV) )  + ( J \oplus \bfp^*(\cV) ) 
\\[1\jot]
	\hZ \cap \cZ = E^1 = J  \oplus \bfp^*(I^1) = ( J \oplus \bfp^*(\hV) ) \cap ( J \oplus \bfp^*(\cV) ) .
\end{gather*}
	Together these equations imply that
\begin{equation}
	\dim (\hZ) + \dim (\cZ) =  \dim ( J \oplus \bfp^*(\hV))  + \dim ( J \oplus \bfp^*(\cV) ) 
\end{equation}
	which forces the inclusions $J \oplus \bfp^*(\hV) \subset \hZ$ and $J \oplus \bfp^*(\cV) \subset \cZ$ to be equalities,
\begin{equation}
	\hZ =  J \oplus \bfp^*(\hV) \quad\text{and}\quad \cZ =  J \oplus \bfp^*(\cV).
\EqTag{commonJ}
\end{equation}
	The proves that $J$ is also an admissible bundle for  the Pfaffian integrable extensions $(\hZ,\hV)$ and $(\cZ,\cV)$. 
	However, it should be emphasized that  \EqRef{commonJ} does not the imply the direct sum decomposition formulas for 
	$\hZ^\infty$ and $\cZ^\infty$ in  \EqRef{DarbAdmis1} and \EqRef{DarbAdmis2}.

	With regards to the definition {\bf [ii]}  in  \StRef{DarbComp}, 
	from the  first two equations in   \EqRef{DarbAdmis1} and \EqRef{DarbAdmis2} it is clear that any maximally compatible 
	extension is compatible. One also easily checks that  \EqRef{IntExtDef}, equations \EqRef{DarbAdmis1} 
	and \EqRef{DarbAdmis2} imply that $\rank(\hJ) = \rank(\cJ)$.  
	Consequently, in addition to the direct sum  decompositions in \EqRef{DarbAdmis1} 
	and  \EqRef{DarbAdmis2} we also have
\begin{equation}
	\hZ =    \hJ  \oplus \bfp^* ( \hV)  \quad\text{and} \quad \cZ =    \cJ  \oplus \bfp^* ( \cV).
\EqTag{hZmore}
\end{equation}

	In terms of this definition, Theorem \StRef{DIext} states that if  $\bfp\:  (\CalE, N) \to (\CalI, M)$ is an integrable extension 
	and $\CalI$ is decomposable, then  there always exists singular Pfaffian systems for $\CalE$  which are compatible.
	We shall  see in Theorem \StRef{DiagDarbComp1}  that maximally compatible
	extensions naturally arise when integrable  extensions are constructed by group reductions using diagonal actions.   


	The next theorem gives  a simple criteria for  an  integrable extension to be  maximally compatible.  We emphasize that
	the Darboux integrability of the differential system $\CalI$ is not required.
	This result will be used in Sections  6 and  11.

\begin{Theorem}\StTag{H}
	Let   $\bfp\:  (\CalE, N) \to (\CalI, M)$ be an integrable extension of  decomposable systems  $\CalE$ and $\CalI$ with the compatible singular systems $\{\hZ, \cZ\}$ and  $\{\hV, \cV\}$ respectively.  Assume that $\CalE$ is Darboux integrable. Then $(\CalE, \CalI )$ is maximally compatible if and only if 
\begin{equation*}
\begin{aligned}
& {\rm {\bf [i]}} & \ker( \bfp_*) \cap \ann(\hZ^\infty) =0 ,\quad & \quad  \ker (\bfp_*) \cap \ann(\cZ^\infty) = 0 \quad {\text \rm and} \\
& {\rm {\bf [ii]}} &  \rank (\hZ^\infty) = \rank (\ker(\bfp_*)) + \rank (\hV^\infty),  \quad &\quad \rank (\cZ^\infty) = \rank (\ker(\bfp_*)) + \rank (\cV^\infty). &
\end{aligned}
\end{equation*}
\end{Theorem}

	The rank equalities {\bf [ii]} imply that the number of  independent Darboux invariants for  $\hZ$ and $\cZ$ 
	exceeds the number of  independent Darboux invariants for  $\hV$ and $\cV$ by exactly 
	the fibre dimension of the submersion $\bfp$.

\begin{proof}
	It is easy to check that \EqRef{DarbAdmis1} and  \EqRef{DarbAdmis2} implies conditions  {\bf [i]} and {\bf[ii]} of Theorem \StRef{H}. 
	
	We prove that compatibility, specifically the second equation in  \EqRef{DarbAdmis0}, 
	and  the first equations in {\bf [i]} and {\bf[ii]} in the statement of the theorem give  \EqRef{DarbAdmis1}.  
	We shall need the following  the rank equalities
\begin{align}
\rank (\hZ^\infty \cap \cZ)& =  \rank(\ker \bfp_*) + \rank (\hV^\infty \cap \cV) \quad \text{\rm and} 
\EqTag{R1}
\\[1\jot]
\rank\bigl(\hZ^\infty \cap \ann(\ker \bfp_*)\bigr)& = \rank (\hV^\infty) = \rank(\bfp^*(\hV^\infty)).
\EqTag{R2}
\end{align}
	which we now derive from  the hypotheses of  Theorem \StRef{H}.

	Since $\CalE$ is Darboux integrable,  \EqRef{Intro7}  implies that $\hZ^\infty + \cZ = T^*N$ and thus
\begin{equation}
	\dim N  = \rank \hZ^\infty + \rank \cZ - \rank ( \hZ^\infty \cap \cZ).
\EqTag{dimWW}
\end{equation}
	Therefore, on account of equation \EqRef{dimWW}, condition {\bf [ii]} in  Theorem \StRef{H},  \EqRef{DarbAdmis0}, 
	 the definition of integrable extension  (to compute $\dim N$) and equation \EqRef{Intro7}, we obtain
\begin{equation*}
\begin{aligned}
	&\rank (\hZ^\infty \cap \cZ)
	 =  \rank \hZ^\infty + \rank \cZ - \dim N
\\
&  \quad = [\rank (\ker \bfp_*) + \rank \hV^\infty] +  [\rank(\ker \bfp_*) + \rank \cV] - [\rank(\ker\bfp_*) + \dim M]
\\
& \quad = \rank(\ker(\bfp_*)) + \rank \hV^\infty +  \rank \cV - \dim M   = \rank(\ker \bfp_*) + \rank (\hV^\infty \cap \cV) .
\end{aligned}
\end{equation*}

	We calculate the dimension of the subbundle $ \hZ^\infty \cap \ann(\ker \bfp_*)$,  first using the transversality condition {\bf[i]} in Theorem \StRef{H} and  then {\bf[ii]},
	to be
\begin{align*}
	&\rank\bigl(\hZ^\infty \cap \ann(\ker \bfp_*)\bigr)
	= \dim(N) - \rank\bigl(\ann(\hZ^\infty \cap \ann(\ker \bfp_*))\bigr)
\\
	&\quad= \dim(N)  - \rank\bigl(\ann(\hZ^\infty)  +\ker \bfp_*\bigr)
	= \dim(N) -  \rank\bigl(\ann(\hZ^\infty)\bigr)   - \rank(\ker \bfp_*)
\\	
	&\quad= \rank(\hZ^\infty) -  \rank(\ker \bfp_*)  = \rank (\hV^\infty) = \rank(\bfp^*(\hV^\infty)).
\end{align*}
	Equations \EqRef{R1} and \EqRef{R2} are now established. 

	The inclusions (see \EqRef{DarbAdmis0})
\begin{equation}
	\bfp^*(\hV) \subset \hZ \quad \text{and}\quad 	\bfp^*(\hV^\infty) \subset \hZ^\infty , 
\EqTag{hVInc}
\end{equation}
	imply  that  we can choose  a subbundle $\hJ \subset \hZ^\infty \cap \cZ$   such that
\begin{equation}
	\hZ^\infty \cap \cZ = \hJ \oplus \bfp^*(\hV^\infty \cap \cV)
\EqTag{cap1}
\end{equation}
	We claim that the bundle $\hJ$  is an admissible subbundle which satisfies \EqRef{DarbAdmis1}.  

	The next step in the proof is  to check that the complementary bundle  $\hJ$  is transverse to $\bfp$, that is, 
	$\hJ \cap  \ann(\ker \bfp_*) = 0$.
	Equation  \EqRef{R1} implies that  $\rank \hJ = \rank\,(\ker \bfp_*)$.
	Note also that  the inclusions \EqRef {hVInc} and  rank equality \EqRef{R2} imply that
\begin{equation}
	\hZ^\infty \cap \ann(\ker \bfp_*) = \bfp^*(\hV^\infty).
\EqTag{BunEq}
\end{equation}
	Since $\rank(\hJ) = \rank(\ker \bfp_*)$,
	transversality is equivalent to the non-degeneracy of the canonical pairing, restricted to $ (\ker \bfp_*) \times \hJ$.
	Let $\alpha \in \hJ$ and suppose  $\alpha(X) = 0$ for all $X \in \ker \bfp_*$, that is, suppose  $\alpha \in \hJ_{\bfp,\semibasic}$.
	Since $\hJ \subset \hZ^\infty$ and $\hJ\subset \cZ$ 
	we have $\alpha \in \hZ^\infty_{\bfp,\semibasic}$ and $\alpha \in \cZ_{\bfp,\semibasic}$.
	Equation \EqRef{BunEq} implies $\alpha \in \bfp^*( \hV^\infty)$ while equation  \EqRef{commonJ} implies that $\cZ_{\bfp,\semibasic}= \bfp^* (\cV) $ and so $\alpha \in \bfp^* (\cV)$.
	Therefore $\alpha \in \bfp^*(\hV^\infty\cap \cV)$, and so by \EqRef{cap1}, $\alpha =0$.
	
	Transversality shows that  $\hJ \cap \bfp^* (\cV) = \hJ \cap \bfp^* (\hV^\infty) = 0$ and therefore, by dimensional considerations,
\begin{equation*}
	 \cZ = \hJ \oplus \bfp^*( \cV), \quad \text{and}\quad 
	\hZ^\infty  = \hJ \oplus \bfp^*(\hV^\infty).
\end{equation*}
	Finally, from the definition of Darboux pair, we have   $\hZ^\infty \cap \cZ \subset E^1$ and $ \cZ^\infty \cap \hZ \subset E^1$
	and therefore $\hJ \subset E^1$.  The assumption that  the singular systems $\{\hZ, \cZ\}$ and  $\{\hV, \cV\}$ are compatible implies that 
\begin{equation}
	\rank E^1-  \rank \bfp^*( I^1) = 	\rank \hZ - \rank \bfp^*( \hV)  = \rank \cZ - \rank \bfp^*( \cV ).
\end{equation}
	Thus, again by dimensional considerations,  we deduce $E^1 = \hJ  \oplus \bfp^* (I^1) $ and $\hJ$ is admissible. 
	
	The proof of  \EqRef{DarbAdmis2} is obtained from  the second equations in {\bf [i]} and {\bf[ii]} 
	by interchanging the roles of  the accents $\check{}$ and $\hat{}$ in the above.
\end{proof}

\section{\kern -2pt Group Theoretic Constructions of Darboux Integrable Systems}
	Let $G$ be a Lie group acting on manifolds $M_1$ and $M_2$.
	In  Section 3 of \cite{anderson-fels-vassiliou:2009a}  a general group-theoretic construction 
	of Darboux integrable systems was given, based upon symmetry reduction with respect to 
	the diagonal action $G_\diag$ on the product manifold  $M_1\times M_2$.
	In this section we extend this result to the more general case  of symmetry reduction by subgroups $L$
	of the product group  $G\times G$. This leads to a purely  group-theoretic construction of  B\"acklund transformations
	between Darboux integrable systems.
 	We also  show,  in the special case of the diagonal action $G_\diag$, that symmetry reduction
	 leads to pairs of  Darboux integrable systems 
	which are  always maximally compatible in the sense described in Section 5.2.  This fact will play an key
	role in Sections 7 and 8  when we 
	characterize those integrable extensions of Darboux integrable equations which arise as group quotients.
		
\subsection{Darboux Integrable by Symmetry Reduction of Sums of EDS}

	To begin, let $\CalK_1$ and $\CalK_2$ be EDS  on manifolds $M_1$ and $M_2$.
	The {\deffont direct sum } $\CalK_1 + \CalK_2$ is the EDS on $M_1 \times M_2$  which is
	algebraically generated by the pullbacks of $\CalK_1$ and $\CalK_2$  to $M_1 \times M_2$ 
	by the canonical projection maps $\pi_a\colon M_1\times M_2 \to M_a$. 
	 {\it Throughout this section 
	we will assume that $\CalK_1$ and $\CalK_2$ are algebraically generated by
	1- and 2-forms.} This granted, 
	the system $\CalK_1+\CalK_2$ is then clearly decomposable  with singular Pfaffian systems
\begin{equation}
	\hW = K^1_1 \oplus T^*M_2
	\quad\text{and}\quad
	\cW = T^*M_1 \oplus K^1_2 .
\EqTag{Wdefs}
\end{equation}
	It is a simple matter to check that   $\{\, \hW, \cW\,\}$ form a Darboux pair if and only if $(K^1_a)^{\infty} = 0 $. 

	Let  $L$ be a subgroup of  the product  group $G \times G$, 
	let $\rho_a : G\times G \to G$, $a=1,2$ be the projection onto the $a^{\text{th}}$ factor and let 
	$L_a = \rho_a(L) \subset  G$.  The 
	action of $L$ on $M_1\times M_2$  is  then given in terms of these projection maps and the actions of $G$ on 
	$M_1$  and $M_2$ by
\begin{equation}
	  \ell \cdot ( x_1, x_2) = (\rho_1(\ell) \cdot x_1, \,\rho_2(\ell) \cdot x_2) \quad \text{for $\ell \in L$}.
\EqTag{Laction}
\end{equation}
	The projection maps $\pi_a: M_1 \times M_2 \to M_a$ are equivariant  with respect the actions of $L$ and $L_a$, that  is,  
	for any $\ell \in L$ and  $x\in M_1 \times M_2$ 
\begin{equation}
	\pi_a(\ell \cdot x) = \rho_a(\ell) \cdot \pi_a(x).
\EqTag{Lequiv}
\end{equation}
	In the special case of the diagonal action   $G_\diag$, each $L_a \cong G$ and the 
	actions of $L_a$ on $M_a$ coincide with the original actions of $G$ on  $M_a$. No assumptions will be made regarding  the dimensions of the Lie groups $L_a$.
 
	We also suppose that $L$ acts regularly on $M_1\times M_2$ and let $\bfq_L\:M_1\times M_2 \to (M_1 \times M_2)/L$
	be the canonical submersion.
	For a  given point $ (p_1, p_2) \in M_1 \times M_2$, let $\iota_{M_1} :M_1 \to  M_1\times M_2 $ 
	and  $\iota_{M_2} :M_2 \to  M_1\times M_2 $  be the inclusion maps. 
\begin{equation}
\begin{gathered}
	\iota_{M_1}(x_1) =  (x_1, p_2) \
	\quad\text{and}\quad   
	\iota_{M_2}(x_2) =  (p_1, x_2) 
	\quad\text{and set}\quad 
\\[2\jot]
	\bfq_{M_a} = \bfq_{L} \circ \iota_{M_a} :  M_a \to (M_1\times M_2)/L
\end{gathered}
\EqTag{qMdef}
\end{equation}
	Let 
	$\bfq_{L_a} : M_a \to M_a/L_a$ be the canonical quotient maps.  By the $L$-equivariance of the
	projection maps $\pi_a$ (see \EqRef{Lequiv}), we can then define  maps 
\begin{equation}
	\bfp_a:(M_1\times M_2)/L \to M_a/L_a
\end{equation}
	such that the diagrams 
\begin{equation}
\begin{gathered}
\begindc{\commdiag}[3]
\obj(-35, 0)[I1]{$M_1\times M_2$}
\obj(0, 0)[I1G]{$(M_1\times M_2)/L$}
\obj(-35, -13)[I1C]{$M_1$}
\obj(0, -13)[I1Z]{$M_1/L_1$}
\mor{I1}{I1G}{$\bfq_L$}[\atleft, \solidarrow]
\mor{I1G}{I1Z}{$\bfp_1$}[\atleft, \solidarrow]
\mor{I1C}{I1Z}{$\bfq_{L_1}$}[\atright, \solidarrow]
\mor{I1}{I1C}{$\pi_1$}[\atright, \solidarrow]

\enddc
\end{gathered}
\quad\text{and }\quad 
\begin{gathered}
\begindc{\commdiag}[3]
\obj(-35, 0)[I1]{$M_1\times M_2$}
\obj(0, 0)[I1G]{$(M_1\times M_2)/L$}
\obj(-35, -13)[I1C]{$M_2$}
\obj(0, -13)[I1Z]{$M_2/L_2$}
\mor{I1}{I1G}{$\bfq_L$}[\atleft, \solidarrow]
\mor{I1G}{I1Z}{$\bfp_2$}[\atleft, \solidarrow]
\mor{I1C}{I1Z}{$\bfq_{L_2}$}[\atright, \solidarrow]
\mor{I1}{I1C}{$\pi_2$}[\atright, \solidarrow]

\enddc
\end{gathered}
\EqTag{H2action}
\end{equation}
	commute. 
	If we make the slightly stronger  assumptions that  the groups  $L_a$ act regularly on $M_a$,  then the maps $\bfq_{L_a}$ 	
	and  $\bfq_{L_a} \circ \pi_a$  are smooth submersions
	and hence,  by  Theorem \StRef{Gcomd},    the maps $\bfp_a$ are smooth submersions. 

	The following theorem summarizes the essential facts  regarding the construction of Darboux integrable systems by 
	symmetry reduction.

  \begin{Theorem}
\StTag{Dreduce}  
	Let $\CalK_a, a=1,2$ be exterior differential systems on $M_a$, $a = 1,2$ which are algebraically generated by
	1- and 2-forms. Assume that $( K_a^1)^\infty  = 0$ and let
\begin{equation}
	\hW= K_1^1 + T^*M_2, \quad \cW= T^*M_1 +K_2^1
\EqTag{Wsing}
\end{equation}
	be the corresponding Darboux pair on  $M_1\times M_2$. 
	Consider a  Lie group  $G$  which  acts freely on $M_a$,  is a {\it common\/} symmetry group of  both $\CalK_1$  
	and $\CalK_2$, and  acts transversely to $\CalK_1$  and $\CalK_2$.
	Assume also that  the actions of $L \subset G \times G$ on $M_1\times M_2$ and  $L_a$ on $M_a$ 
	are regular and set
\begin{equation}
	M = (M_1 \times M_2)/L, \quad \hV = (K^1_1+T^*M_2)/L, \quad \cV=(T^*M_1+K^1_2)/L.
\EqTag{WMODG}
\end{equation}
	Finally, assume that $\hV^\infty$ and $\cV^\infty$ are constant rank bundles.

\medskip
\noindent
{\bf [i]}
	Then  $\{\hV, \cV \}$
	define  a Darboux pair on  $M$  with \footnote{The notation  $0 +TM^*_2$ indicates that this bundle is
	to be viewed as subbundle of $T^*M_1 +T^*M_2$.}
\begin{gather}
	\hV^\infty = (0 +T^*M_2)/ L  = \hW^\infty/L,  \quad 
	\cV^\infty = (T^*M_1 + 0)/ L = \cW^\infty/L  \quad\text{and}
\EqTag{WMODG3}
\\[2\jot]
	\rank(\hV^\infty)  = \dim M_2  - \dim L_2, \quad \rank(\cV^\infty)  = \dim M_1 - \ \dim L_1.
\EqTag{WMODG4}
\end{gather}

\medskip
\noindent
{\bf [ii]}
	The quotient differential system  $ \CalI = (\CalK_1 + \CalK_2)/L$ on  $M$
	is Darboux integrable with singular Pfaffian systems $\{\hV, \cV\}$ as defined in  \EqRef{WMODG}.  

\medskip
\noindent
{\bf [iii]} The   
	 bundles $\hV^\infty$ and $\cV^\infty$  are given by 
\begin{equation}
\hV^{\infty} = \bfp_2^*(T^*(M_2/L_2)),\quad\text{and}\quad \cV^{\infty} = \bfp_1^*(T^*(M_1/L_1)).
\EqTag{QINVS}
\end{equation}

\noindent
For parts {\bf[iv]}  and {\bf [v]}, we restrict to the special case where $L = G_\diag$.

\medskip
\noindent
{\bf[iv]}  The quotient map $\bfq_{G_\diag}:(\CalK_1\times \CalK_2, M_1 \times M_2) \to (\CalI, M) $ defines a maximally compatible integrable extension with respect to the singular Pfaffian systems  \EqRef{Wsing} and \EqRef{WMODG}.

\medskip
\noindent
{\bf [v]} Assume that $M_1$ and $M_2$ are connected. Then  for any point $(p_1,p_2) \in M_1 \times M_2$, the maps 
	(see \EqRef{qMdef})
\begin{equation}
	\bfq_{M_1}\:  M_1 \  \to M
\quad\text{and}\quad
	\bfq_{M_2}\:  M_2  \to M 
\EqTag{qMa} 
\end{equation}
	define imbedded, maximal integral manifolds for  $\hV^\infty$ and $\cV^\infty$ respectively.
\end{Theorem}

	We prove all five parts of this theorem for the special case of the diagonal action $G_\text{diag}$ in  Section  6.1. 
	The generalizations of parts {\bf[i]}  -- {\bf [iii]} for non-diagonal actions is shown to reduce to the diagonal case in Section 6.2. 
	To prove {\bf[iv]} and {\bf [v]} we will make repeated use of the simple observation that for free, diagonal group actions the distribution
	$ \bfGamma_{G_\diag}$ of infinitesimal generators for the action of $G_\diag$ satisfies
\begin{equation}
	 \bfGamma_{G_\diag} \cap  (TM_1+0) =  \bfGamma_{G_\diag} \cap (0+TM_2) = \{0\}.
\EqTag{FreeIntersection}
\end{equation}
	For non-diagonal actions equations \EqRef{FreeIntersection} fail to hold and statements  {\bf[iv]} and {\bf [v]} are not true.

\begin{Remark}
	Theorem    \StRef{Dreduce}  can be combined with  Theorem A to construct  B\"acklund transformations for Darboux integrable systems.  Indeed, starting from the hypotheses of Theorem  \StRef{Dreduce} 
	pick two different choices, say $G_1 \subset   G\times G$ and $G_2 \subset G\times G$  for the subgroup $L$  in Theorem  
	\StRef{Dreduce} and  let $H =G_1\cap G_2$.  Then three applications  
	of Theorem  \StRef{Dreduce}  prove that $(\CalK_1 +\CalK_2)/G_1$,  $(\CalK_1 +\CalK_2)/G_2$ and $(\CalK_1 +\CalK_2)/H$ are all Darboux integrable systems. Theorem A  then proves that 
	 $(\CalK_1 +\CalK_2)/H$ is a B\"acklund transformation  between  $(\CalK_1 +\CalK_2)/G_1$ and  $(\CalK_1 +\CalK_2)/G_2$.   See also Theorem  \StRef{BigBack}  and the discussion that follows.

\end{Remark}

\begin{Remark}	
	Equations \EqRef{QINVS} clearly imply that
	if $g$ is a smooth function on $M_2/L_2$, then  $g \circ \bfp_2$ is a first integral for $\hV^\infty$ and hence a Darboux invariant
	for $\CalI$. 
	Of course,  smooth functions on $M_2/L_2$ are in one-to-one correspondence with $L_2$-invariant functions on $M_2$. 
	Accordingly,  {\it there is a one-to-one correspondence between the Darboux invariants for the Darboux pair $\{\hV, \cV\}$  on $(M_1 \times M_2)/G$
	and the $L_a$-invariant functions on $M_a$.}
\end{Remark}

\begin{Remark}
\StTag{ClosedDiag}
	Note  that  $G_{\diag}$  is a closed subgroup of $G\times G$. Therefore, 
	if $G$ acts freely and regularly on $M_1$ and $M_2$ then,  by Remark \StRef{NReduction},
	the action of $G_{\diag}$ is free and regular on $M_1\times M_2$.
\end{Remark}

\subsection{Reduction by Diagonal Actions}
	In this section we prove  Theorem \StRef{Dreduce} for diagonal actions $ L = G_\diag$.
	Part {\bf[i]} and equation  \EqRef{WMODG3}  of   Theorem \StRef{Dreduce} for
	diagonal actions is proved in Theorem 3.2 in  \cite{anderson-fels-vassiliou:2009a}. 
	Equation  \EqRef{WMODG4}, in the case of diagonal actions, follows from   \EqRef{WMODG3},  the transversality condition
	\EqRef{FreeIntersection}, \EqRef{rankI1}, and the assumption that  $G_\diag$ acts freely.

	To prove  part {\bf [ii]} of Theorem \StRef{Dreduce}   we first observe, by applying  \EqRef{IntExt15} to $\CalK_1$ and
	$\CalK_2$, 
	that the 2-forms generators obtained in this way for $\CalK_1+\CalK_2$ are $G$-basic and reduce to 2-form generators for  $\CalI$ (see also Section 7.3). Using this observation it is not difficult to 
	show that $\CalI$ is decomposable where the Pfaffian systems $\{\hV, \cV\}$ defined by  
	\EqRef{WMODG} are actually singular systems for $\CalI$. The details are given in Corollary 3.4 in 
	\cite{anderson-fels-vassiliou:2009a} when $\CalK_1$ and $\CalK_2$ are Pfaffian systems. 
	A simple extension of these arguments using the above observation regarding the 2-form generators proves the more general case. 

	It remains to prove parts  {\bf [iii]}, 
	{\bf [iv]} and {\bf [v]} (for $ L = G_{\diag}$). We start with part {\bf [iv]} which 
	is a simple corollary to the following more general 
	result which we shall  also need in Section 9.

\begin{Theorem} 
\StTag{DiagDarbComp1} 
	Let $(\CalK_1, M_1)$ and $(\CalK_2, M_2)$ be exterior differential systems (generated algebraically by $1$- and  2-forms) 
	with a common symmetry group 
	$G$ and with $(K_a^1)^\infty = 0$.  Suppose that the action of $G$ 
	satisfies all the hypotheses of Theorem \StRef{Dreduce} and  let  
	$\CalI$ be the Darboux integrable differential system $\CalI = (\CalK_1\times \CalK_2)/G_\diag $  (see Theorem  \StRef{Dreduce} 
	part {\bf[i]})
	with singular  Pfaffian systems  
\begin{equation*}
	\hV =(K_1^1+T^*M_2)/G_\diag \quad\text{and}\quad \cV=(T^*M_1+K_2^1)/G_\diag.
\end{equation*}

	Now let  H be a Lie subgroup of G which  acts regularly on $M_1$, $M_2$.   Then 

\medskip
\noindent
{\bf [i]} $\CalE = (\CalK_1\times \CalK_2)/H_\diag$   is a Darboux integrable  differential system with singular Pfaffian systems 
\begin{equation*}
	\hZ =(K_1^1+T^*M_2)/H_\diag \quad\text{and}\quad \cZ=(T^*M_1+K_2^1)/H_\diag; 
\end{equation*}

\medskip
\noindent
{\bf [ii]}  $\CalE$ is an integrable extension of  $\CalI$  with respect  to the orbit projection map 
\begin{equation}
	\bfp \: (M_1 \times M_2)/H_\diag \to  (M_1 \times M_2)/G_\diag; \quad \text{and}
\end{equation}
\medskip
\noindent
{\bf [iii]}  
	$(\CalE, \CalI)$ is maximally compatible with respect to the singular Pfaffian systems $\{\hZ, \cZ\}$ and $\{\hV, \cV\}$.
\end{Theorem}

\begin{proof} 
	We first remark that the hypotheses that $G$ acts freely on $M_a$, is a common symmetry group of both $\CalK_1$ and $\CalK_2$,
	and acts transversely to $\CalK_1$ and $\CalK_2$ implies that the same is true of $H$. Remark \StRef{ClosedDiag} shows that
	the action of $H_\text{diag}$ is  regular on $M_1\times M_2$.   Part {\bf [i]} of  Theorem \StRef{DiagDarbComp1}  now follows
	from  part {\bf [i]} of Theorem  \StRef{Dreduce}  with $L = H_{\diag}$. 
	Note that by Remark \StRef{NReduction}, the hypothesis that $H$ acts regularly on $M_1$ and $M_2$ is unnecessary in the
	case  where $H$ is a closed subgroup of $G$.

	Part {\bf [ii]} of  Theorem \StRef{DiagDarbComp1} follows directly from  Theorem \StRef{IntExt3}.
	
	To prove  part {\bf [iii]} of  Theorem \StRef{DiagDarbComp1} we shall check conditions   {\bf [i]}  and {\bf[ii]} of Theorem \StRef{H}. To verify condition {\bf [i]}  in Theorem \StRef{H} we first note that  (see \EqRef{WMODG3} with $L = H_\diag$)
\begin{equation}
	\hZ^\infty=(0+T^*M_2)/H_\diag\quad\text{and}\quad\cZ^\infty=(T^*M_1+0)/H_\diag .
	\EqTag{NCE1}
\end{equation}
	Also, by  \EqRef{kerp},  we have 
\begin{equation}
	\ker (\bfp_*)= \bfq_{H_\diag*}(\bfGamma_{G_{\text{diag}}}) .
\EqTag{NCE2}
\end{equation}
	Let $ R \in \ker (\bfp_*) \cap \ann(\hZ^\infty)$.   Then the annihilator of the first equation in \EqRef{NCE1} and \EqRef{NCE2} imply
	there exists vectors $R_1 \in TM_1 + 0$ and  $R_2 \in \bfGamma_{G_{\text{diag}}}$ such that
\begin{equation}
	R = \bfq_{H_\diag*}(R_1) = \bfq_{H_\diag*}(R_2)
\EqTag{NCE3}
\end{equation}
	in which case $R_1-R_2 \in \bfGamma_{H_\diag}$.  Since $H_\diag$ is  the restriction of the diagonal action $G_\diag$, this implies that
	$R_1 - R_2 \in  \bfGamma_{G_{\text{diag}}}$  and therefore $R_1 \in \bfGamma_{G_{\text{diag}}}$.   
	Equation \EqRef{FreeIntersection} then implies that $R_1=0$  and consequently $R=0$. 
	This proves the first equation in  {\bf [i]} in Theorem \StRef{H} and the proof of the  second equation is similar.
\smallskip

	To check part {\bf [ii]} of Theorem  \StRef{H}, we note that equations  \EqRef{WMODG4}  give
\begin{equation}
	\rank (\hZ^\infty) =  \dim M_2 - \dim H
\quad\text{and}\quad\
	\rank (\hV^\infty) =  \dim M_2 - \dim G.
\EqTag{NS1}
\end{equation}
	Moreover, by equation \EqRef{NCE2}, it follows that
\begin{equation}
\rank (\ker(\bfp_*)) = \dim G - \dim H.
\EqTag{NS3}
\end{equation}
	Equations \EqRef{NS1} and \EqRef{NS3} lead to  the first equation in condition {\bf [ii]} of Theorem \StRef{H}. 
	The proof of the second equation is similar.
\end{proof}

We can now complete the proof of Theorem  \StRef{Dreduce} for diagonal actions.

\begin{proof} [Proof of Theorem  \StRef{Dreduce}, Part {\bf [iv]} for diagonal actions]
	We simply apply  Theorem  \StRef{DiagDarbComp1} with $H$ as the identity group. 
\end{proof}

\begin{proof}[Proof  of Theorem \StRef{Dreduce}, Part {\bf [iii]} for diagonal actions $L = G_\diag$]
	We prove the first formula in  \EqRef{QINVS}. Starting from \EqRef{WMODG3},  we first use Theorem \StRef{RedGInfty}   to write
\begin{equation}
	\bfq_{G_\diag}^*(\hV^\infty) =  (\hW^\infty)_{G_\diag, \semibasic} 
\EqTag{hVH}
\end{equation}
	From the assumption that  $(K^1_1)^{\infty} =0$, we immediate conclude that 
\begin{equation}
	\hW^\infty  =  0+  T^*M_2
\end{equation}
	from which it then follows that 
\begin{equation}
	(\hW^\infty)_{G_{\diag},\semibasic} =  0+ (T^*M_2)_{G,\semibasic} =  \pi^*_2((T^*M_2)_{G, \semibasic})
\EqTag{Winftysb}
\end{equation}
	Also, the assumption that  $G$ acts regularly on $M_2$ gives,  by \EqRef{Asb}, 
\begin{equation}
	(T^*M_2)_{G, \semibasic} =  \bfq^*_{G_2}\bigl(T^*(M_2/ G) \bigr).
\end{equation}
	The combination of  this equation, equation \EqRef{Winftysb}, 
	and the commutativity of the second diagram in \EqRef{H2action} yields
\begin{equation}
	(\hW^\infty)_{G_{\diag}, \semibasic} =  \pi^*_2((T^*M_2)_{G, \semibasic})
	=  \pi^*_2\bigl(\bfq_{G_2}^*(T^*(M_2/ G))\bigr) =  \bfq^*_{G_\diag}\bigl(\bfp^*_2(T^*(M_2/ G))\bigr). 		
\end{equation} 
	 This equation, together with  \EqRef{hVH}, now leads to the desired result  \EqRef{QINVS} . 
	The formula for $\cV^\infty$ is similarly established.	
\end{proof}
 	
\begin{proof}[Proof  of  Theorem \StRef{Dreduce}, Part {\bf [v]}] 
	The fact that the action of $G$ is free on each $M_a$   immediately implies that  the maps $\bfq_{M_a}\: M_a \to M$ are one-to-one.
	Equation \EqRef{FreeIntersection} shows that
	$\ker  \bfq_{M_1, *}  =  0 $ and  $\ker  \bfq_{M_2, *} = 0$
	and therefore the maps  $\bfq_{M_a}$ are immersions.
	
	Let  
\begin{equation}
	N_1 =  \bfq_{M_2}(M_2)   = \bfq_{G_\diag}(\{\,p_1\,\}  \times M _2) 
	\quad\text{and} \quad 
	N_2 =  \bfq_{M_1}(M_1) =\bfq_{G_\diag}(M_1 \times \{\,p_2\,\}),
\EqTag{N1N2def2}
\end{equation}
	and let $x_0 =\bfq_{G_\diag}(p_1, p_2) $   and $[p_a] = \bfq_{G_a}(p_a)$.  We  will show that 
\begin{equation}
	N_a = \bfp_a^{-1}([p_a])  \subset  (M_1 \times M_2)/G_\diag.
\EqTag{N1N2def1}
\end{equation}
	Then, in view of  Theorem \StRef{Dreduce}, part {\bf[iii]}, $N_1$ is  
	an integral manifold of $\cV^\infty$  and  $N_2$ is an integral manifold of $\hV^\infty$, each through the point $x_0$.
	Since the maps $\bfp_a$ are smooth submersions,  the  level sets $N_a$ are smooth,  imbedded submanifolds of  
	$ M = (M_1 \times M_2)/G_\diag$ with 
\begin{align*}
	\dim N_1& = \dim (M) - \dim(M_1/G) = \dim(M) - \rank(\cV^\infty)  
	\quad\text{and}\quad  
\\	
	\dim  N_2 &= \dim (M) - \dim(M_2/G) = \dim(M) - \rank(\hV^\infty).
\end{align*}  
	This shows that  $N_1$ and  $N_2$ are integral manifolds of maximal dimension. 
	By definition,  $ N_1$  and $N_2$ are  connected whenever  $M_2$ and $M_1$ are connected. 
	To  show that $N_1$ is maximal,  let    $\psi :P \to M$  be a connected integral manifold of 
	$\hV^\infty$ through $x_0$. Then by \EqRef{QINVS}  $(\bfp_1\circ \psi)^* = 0$ and so $ \psi(P) \subset  \bfp_1^{-1}([p_1]) \subset N_1$.  

	It remains only  to check \EqRef{N1N2def1}.  We use the  commutativity of   \EqRef{H2action} to calculate
\begin{equation}
	\bfp_1^{-1}([p_1])  = \bfq_{G_\diag} \bigl( (\bfq_{G_1} \circ \pi_1)^{(-1)}([p_1])\bigr) 
	= \bfq_{G_\diag} \bigl ( \pi_ 1^{(-1)}(G\cdot p_1) \bigr) =  \bfq_{G_\diag}\bigl(( G \cdot p_1) \times M_2\bigr).
\EqTag{N1ver2}
\end{equation}
	But for each $g\in G$ and $x_2\in M_2$ we have
\begin{equation*}
	( g \cdot p_1, \, x_2) =  g \cdot_{\text{\diag}}(p_1, \, (g^{-1}) \cdot x_2) .
\end{equation*}
	Therefore $G\cdot p_1 \times M_2 = G_{\diag}(\{p_1\} \times M_2)$ 
	and consequently,  by \EqRef{qMdef} and \EqRef{N1ver2}
\begin{equation}
	\bfp_1^{-1}([p_1])   = \bfq_{G_{\diag}}\bigl ((G \cdot p_1) \times M_2 \bigr)= \bfq_{G_{\diag}}(\{p_1\} \times M_2) = \bfq_{M_2}(M_2) = N_1.
\end{equation}
	The second equation in \EqRef{N1N2def1} is similarly proved.
\end{proof}
	
\begin{Remark} Integrable Pfaffian systems whose maximal integral manifolds are all imbedded submanifolds are 
	called {\deffont regular integrable Pfaffian systems} \cite{palais:1957a} . Theorem \StRef{Dreduce}, part {\bf [v]}  therefore implies that $\hV^\infty$ and $\cV^\infty$
	are regular integrable Pfaffian systems.
\end{Remark}

\begin{Remark}
	The hypothesis that the $M_a$ are connected  in part {\bf [v]} of Theorem \StRef{Dreduce}  implies 
	that the inclusions $\iota_{M_a}:M_a \to M_1 \times M_2$ are maximal integral manifolds of $\hW^\infty $ and $\cW^\infty$. 
	Therefore, by the Darboux compatibility condition established in part {\bf [iii]} of Theorem \StRef{Dreduce}, the manifolds $N_1$ and $N_2$ in equation \EqRef{N1N2def2} are integral manifolds of maximal dimension for $\hV^\infty$ and $\cV^\infty$. 
	Part {\bf [v]} shows, in addition, that these integral manifolds are maximal and imbedded.
\end{Remark}

\subsection{Reduction by  non-Diagonal Actions}

	We return to the general case of  non-diagonal actions $L \subset G\times G$ with the objective of 
	proving  parts {\bf [i]}--{\bf [iii]}  of Theorem  \StRef{Dreduce} in their full generality.
	Let
\begin{equation}
	A_1 = \rho_1(\ker \rho_2) = \{ g_1\in G \,|\, (g_1,e) \in L\} \quad\text{and} \quad 
	A_2 = \rho_2(\ker \rho_1)  =  \{ g_2\in G \,|\, (e, g_2) \in L\} .
\EqTag{Asubgroup}
\end{equation} 
	These are closed normal subgroups of $L_1$ and $L_2$.  
	The product $A_1\times A_2 \subset L  $  is therefore also normal and closed.  Let  $\tilde L = L/(A_1\times A_2)$. Then,  
	from the composition of  the epimorphisms $L \to L_a  \to L_a/A_a$ (with  kernel $A_1 \times A_2$) one  obtains  the isomorphisms
\begin{equation}
\psi_a: \tilde L  \to L_a/A_a \quad\text{defined by }\quad \psi_a( \ell \, A_1 \times A_2  ) = \rho_a(\ell) A_a \quad \text{ for} \ \ell \in L \ .
\EqTag{tildeL}
\end{equation}
	
	We assume the groups $L_a$ act freely and regularly on $M_a$.
	Then,  by Remark \StRef{NReduction},   the groups $A_a$ act freely and regularly on $M_a$, 
	and the groups $L_a/A_a$  act freely and regularly on $M/A_a$.		
	Likewise, $A_1\times A_2$ acts freely and regularly on $M_1\times M_2$, 
	and the quotient group $\tilde L = L/(A_1\times A_2)$  acts freely and regularly on $(M_1\times M_2)/(A_1\times A_2)$.
	The isomorphism \EqRef{tildeL} implies that $\tilde L$ acts freely and regularly  on  $M_1/A_1$  and  $M_2/A_2$. 
	Remark \StRef{ClosedDiag}  shows that the diagonal action $\tilde L_{\diag}$ on   $M_1/A_1 \times M_2/A_2$  is  free and  regular.

	The following lemma is the key to reducing the non-diagonal versions  of parts {\bf [i]} and {\bf [ii]}  of Theorem  \StRef{Dreduce}
	to the corresponding diagonal versions - it gives a canonical identification of the quotient space $(M_1\times M_2)/ L$
	with a quotient space constructed using the aforementioned diagonal action of $\tilde L$.

\begin{Lemma}\StTag{CDPSI} There exists a canonically defined diffeomorphism $\Phi$  such that the diagram
\begin{equation}
\begin{gathered}
\begindc{\commdiag}[3]
\obj(0, 0)[I]{$M_1\times M_2$}
\obj(50, -15)[H]{$(M_1/A_1 \times M_2/A_2)/\tilde L_{\diag}$}
\obj(50, 0)[H2]{$M_1/A_1 \times M_2/A_2$}
\obj(0, -15)[I1]{$(M_1\times M_2)/L$}
\mor{I}{H2}{$\bfq_{A_1 \times A_2}$}[\atleft, \solidarrow]
\mor{H2}{H}{$\bfq_{\tilde L_{\diag}}$}[\atleft, \solidarrow]
\mor{I1}{H}{$\Phi$}[\atright, \solidarrow]
\mor{I}{I1}{$\bfq_L$}[\atright, \solidarrow]
\enddc
\EqTag{nondiag1}
\end{gathered}
\end{equation}
commutes.
\end{Lemma}

\begin{proof}  
	The canonical diffeomorphism $ \Phi_2:  (M_1\times M_2)/(A_1\times A_2) \to M_1/A_1\times M_2/A_2$
	is $\tilde L - \tilde L_\diag$ equivariant  and  
	hence induces the right square in the commutative diagram
	\EqRef{nondiag2} (below), where $\tilde \Phi_2$ is a diffeomorphism. 
	Because $A_1\times A_2\subset L$ is normal we may also construct  (see Remark	  \StRef{NReduction})
	the left-hand square in  \EqRef{nondiag2}, where $\tilde \Phi_1$ is a diffeomorphism.
\begin{equation}
\begin{gathered}
\begindc{\commdiag}[3]
\obj(0, 0)[I]{$\ M_1\times M_2\ $}
\obj(40, 0)[J]{$\ (M_1\times M_2)/(A_1\times A_2)\ $}
\obj(40, -15)[J2]{$\ (M_1\times M_2)/(A_1\times A_2)/\tilde L\ $}
\obj(87, -15)[H]{$\ (M_1/A_1 \times M_2/A_2)/\tilde L_{\diag}\ . $}
\obj(87, 0)[H2]{$\ M_1/A_1 \times M_2/A_2\ $}
\obj(0, -15)[I1]{$\ (M_1\times M_2)/L\ $}
\mor{I}{J}{$\bfq_{A_1\times A_2}$}[\atleft, \solidarrow]
\mor{J}{H2}{$\Phi_2$}[\atleft, \solidarrow]
\mor{H2}{H}{$\bfq_{\tilde L_{{\rm \diag}}}$}[\atleft, \solidarrow]
\mor{I1}{J2}{$\tilde\Phi_1$}[\atright, \solidarrow]
\mor{J2}{H}{$\tilde\Phi_2$}[\atright, \solidarrow]
\mor{I}{I1}{$\bfq_L$}[\atright, \solidarrow]
\mor{J}{J2}{$\bfq_{\tilde L}$}[\atleft, \solidarrow]
\enddc
\EqTag{nondiag2}
\end{gathered}
\end{equation}
This diagram proves the lemma, with $\Phi =  \tilde \Phi_2 \circ \tilde \Phi_1$.
\end{proof}

	Our goal is to prove that  $(\CalK_1 \times \CalK_2)/L$ is Darboux integrable. We shall do this by first using
	Theorem \StRef{Dreduce} to  prove that	$\bigl( \dfrac{\CalK_1 \times \CalK_2}{A_1 \times A_2}\bigr)/ \tilde L_\diag$
	 is Darboux integrable and then using the diagram \EqRef{nondiag2}	to identify these two differential systems.

\begin{proof}[Proof  of Theorem \StRef{Dreduce}, part {\bf [i]}]   
	Let $\tilde \CalK_a = \CalK_a/A_a$ be the reduced differential systems  
	on $M_a/A_a$.   We shall  check that these systems and the actions of $\tilde L$ on $M_a/A_a$
	satisfy all the hypotheses of Theorem   \StRef{Dreduce}.   
	First we need to show $\CalK_a/A_a$  are generated by $1$-forms and $2$-forms. 
	This follows from the fact that $\CalK_a\to \CalK_a/A_a$ are integrable extensions,
	see Theorem 2.1 and equation  \EqRef{IntExt15} in particular.

	Theorem  \StRef{RedGInfty} 
	shows that  $(\tilde K^1_a)^\infty=0$.  We have already noted that $\tilde L$  acts freely and regularly 
	on $ M_a/A_a$.  It is easy to check that $\tilde L$ is a symmetry group of  $\tilde \CalK_a$ and is
	transverse.  Then,  in accordance with equations \EqRef{WMODG} we set
\begin{equation}
	 \tilde M = \bigl(M_1/A_1 \times M_2/A_2 \bigr)/\tilde L_{\diag}, 
	\quad\hZ=\bigl(\tilde K_1^1 + T^*(M_2/A_2) \bigr)/\tilde L_{\diag}
	\text{ \ and \ } 
	\cZ= \bigl(T^*(M_1/A_1) + \tilde K^1_2\bigr)/\tilde L_{\diag}. 
\EqTag{newZhat}
\end{equation}
	The last hypotheses of  Theorem   \StRef{Dreduce}  requires that we  
	check that $\cZ^\infty$ and $\hZ^\infty$ are constant rank bundles.  By the definition of the product action  of $A_1\times A_2$ 
	on $M_1\times M_2$ we have
\begin{equation*}
	\tilde K^1_1  + T^*(M_2/A_2) = (K^1_1 + T^*M_2)/(A_1 \times A_2) \quad \text{and} \quad 
	T^*(M_1/A_1) + \tilde  K^1_2  =  (T^*M_1 + K^1_2)/(A_1 \times A_2). 
\end{equation*}
	Equations \EqRef{WMODG}, the  application of Theorem \StRef{Gcomd} to the commutative diagram \EqRef{nondiag1}, 
	and these equations  lead to
\begin{equation}
\begin{aligned}
	\hV	& = (K_1^1+T^*M_2)/L  = \Phi^*\bigl((K_1^1+T^*M_2)/(A_1 \times A_2)/\tilde L_\diag \bigr) = \Phi^*(\hZ) \quad\text{and} \quad 
\\[2\jot]
	\cV  &  = (T^*M_1+K_2^1)/L = \Phi^*\bigl((T^*M_1+K_2^1)/(A_1 \times A_2)/\tilde L_\diag \bigr) =  \Phi^*(\cZ).
\end{aligned}
\EqTag{Phipb}
\end{equation}
	The hypotheses that $\hV^\infty$ and $\cV^\infty$ are constant rank bundles now  implies that $\cZ^\infty$ and $\hZ^\infty$ 
	are constant rank bundles.

	The application of Theorem \StRef{Dreduce}  part {\bf[i]}, for diagonal actions,  then implies
	that  $\{\hZ, \cZ\}$ is a Darboux pair with
\begin{equation}
	\hZ^\infty = \bigl( 0  + T^*(M_2/A_2)\bigr)/\tilde L_{\diag} \quad\text{and}\quad 
	\cZ^\infty = \bigl(T^*(M_1/A_1) +  0\bigr)/\tilde L_{\diag}.
\end{equation}
	  We conclude, again by \EqRef{Phipb},  that $\{\hV, \cV\}$ is a Darboux pair with  
	$\hV^\infty$,  $\cV^\infty$ given by \EqRef{WMODG3}.
\end{proof}

\begin{proof}[Proof  of Theorem \StRef{Dreduce}, part {\bf[ii]}]    
	We again utilize the commutative diagram \EqRef{nondiag1} in Lemma 6.8 and apply part [ii] of Theorem 6.1 in the case of diagonal actions 
	(which was established in Section 6.2)  to  conclude that $(\CalK_1/A_1\times \CalK_2/A_2 )/\tilde L_\diag$ is 
	Darboux integrable with $\hZ$ and $\cZ$ in equation \EqRef{newZhat}
	as the singular Pfaffian systems. Since $(\CalK_1+\CalK_2)/L =\Phi^*( \CalK_1/A_1\times \CalK_2/A_2 )/\tilde L_\diag$,
	it is Darboux integrable, while  \EqRef{Phipb}  and the commutative diagram  \EqRef{nondiag1} shows that $\hV$ and $\cV$ in  
	\EqRef{Phipb} are the
	singular Pfaffian systems.
\end{proof}

\begin{proof}[Proof  of Theorem \StRef{Dreduce}, part {\bf[iii]}]  
	The counter-parts of the commutative diagrams in \EqRef{H2action}, as applied to the diagonal action  of $\tilde L_\diag$, are
	the commutative diagrams ($a = 1, 2$)
\begin{equation}
\begindc{\commdiag}[3]
\obj(-35, 0)[I1]{$M_1/A_1\times M_2/A_1$}
\obj(10, 0)[I1G]{$(M_1/A_1\times M_2/A_2)/\tilde L_\diag$}
\obj(-35, -13)[I1C]{$M_a/A_a$}
\obj(10, -13)[I1Z]{$(M_a/A_a)/ \tilde L$\ .}
\mor{I1}{I1G}{$\bfq_L$}[\atleft, \solidarrow]
\mor{I1G}{I1Z}{$\tilde \pi_a$}[\atleft, \solidarrow]
\mor{I1C}{I1Z}{$\bfq_{\tilde L}$}[\atright, \solidarrow]
\mor{I1}{I1C}{$ \pi'_a$}[\atright, \solidarrow]
\enddc
\end{equation}
	Therefore,  by part {\bf [iii]} of Theorem  \StRef{Dreduce} (which we have already verified for diagonal actions) it follows that
\begin{equation}
	\hZ^\infty= \tilde \pi_2 ^* \bigl( T^*((M_2/A_2)/\tilde L)\bigr) \quad \text{and} 
	\quad  \cZ^\infty= \tilde \pi_1 ^*\bigl( T^*((M_1/A_1)/\tilde L)\bigr).
\EqTag{newhZ}
\end{equation}
	
	Let $\Psi_a \colon  (M_a/A_a)/\tilde L \to M_a/L_a$ be the canonical smooth diffeomorphisms, 
	let  $\tilde \bfp_a \colon (M_1/A_1\times M_2/A_2)/\tilde L_{\diag} \to M_a/L_a$ be the smooth projection maps defined 
	by $\tilde \bfp_a= \Psi_a \circ \tilde \pi_a  $ 
	and note, on account of \EqRef{nondiag1}, that the projection maps $\bfp_a \colon (M_1\times M_2)/L \to M_a/L_a$ satisfy $\bfp_a = \tilde \bfp_a \circ \Phi $. 
	Equation \EqRef{newhZ} then yields
\begin{equation*}
	\hZ^\infty= \tilde \pi_2 ^* \bigl( T^*((M_2/A_2)/\tilde L)\bigr) 
	= \tilde \bfp_2^* \circ \bigl(\Psi_2^{-1})^*(T^* ((M_2/A_2)/\tilde L) \bigl) =  \tilde \bfp_2^*\bigl(T^* (M_2/L_2) \bigr)
\end{equation*}
	and hence
\begin{equation}
	\hV^\infty =\Phi^* \hZ^\infty= \Phi^* \circ \tilde \bfp_2^*\bigl(T^*(M_2/L_2)\bigr)= \bfp_2^*\bigl(T^*(M_2/L_2)\bigr),
\end{equation}
	as required. The formula for $\cV^\infty$ is similarly derived.
\end{proof}
\section{ Vessiot Algebras as Fundamental Invariants for Darboux Integrable Systems}

	The fundamental invariant for any Darboux integrable differential system is the Vessiot algebra.
	In  Section 7.1  we  recall the  definition of  this Lie algebra  and we prove Theorem C. 
	Section 7.2 shows that the prolongation  of a Darboux integrable system is again 
	Darboux integrable  and that the two Vessiot algebras are isomorphic.
	In Section 7.3  we calculate the Vessiot algebra for  the Darboux integrable
	systems constructed  in Section 6.2.  Many of  results will be used in Section 8 where we 
	characterize those integrable extensions of Darboux integrable equations which arise as group quotients through Theorem A.

\subsection{The Vessiot Algebra of a Darboux Integrable Differential System}
	We begin by recalling the  fundamental  technical result of   \cite{anderson-fels-vassiliou:2009a} .
	\begin{Theorem}
\StTag{AFV}
	Let  $(\CalI, M)$ be  a Darboux integrable system with singular systems $\{\hV, \cV\}$. Then there  
	exists,  about each point  of $M$, 0-adapted  (see \EqRef{ZeroC})  local coframes $ \{\,\bftheta_X, \,\bfheta,\, \bfhsigma,\,  \bfceta,\, \bfcsigma\,\}$   and
	$\{\,\bftheta_Y,\, \bfheta,\, \bfhsigma,\,  \bfceta,\, \bfcsigma\,\} $ satisfying the structure equations 
\begin{equation}
\begin{aligned}
	d \bfhsigma& = 0, \quad 
	d\bfheta =\bfhA \, \bfhsigma \wedge \bfhsigma + \bfhG\, \bfheta\wedge \bfhsigma , \
	\quad d \bfcsigma = 0, 
	\quad d\bfceta = \bfcF\, \bfcsigma \wedge \bfcsigma+\bfcH\, \bfceta\wedge \bfcsigma,
\\	
\extd \bftheta_X
&
	=  \frac12  \bftildeA  \, \bfhpi  \wedge \bfhpi   + \frac12 \bftildeB\, \bfcpi  \wedge \bfcpi \
          +\frac12 \bfC\,\bftheta_X\wedge \bftheta_X   +  \tilde \bfM \, \bfhpi \wedge \bftheta_X, \quad \text{and}
\\
	\extd \bftheta_Y&	=  \frac12  \bftildeE \, \bfhpi  \wedge \bfhpi  +  \frac12 \bftildeF \, \bfcpi\wedge \bfcpi  
        - \frac12 \bfC\, \bftheta_Y \wedge \bftheta_Y  +  \tilde \bfN \, \bfcpi \wedge \bftheta_Y,
\end{aligned}
\EqTag{4Adapted3}
\end{equation}
where $\bfhpi=(\bfhsigma, \bfheta)$ and $ \bfcpi=(\bfcsigma, \bfceta)$.
	The coefficients  $\bfC = [C^k_{ij}]$ are constants and  the corresponding dual frames 
	$\{ \bfpartial_{\bftheta_X}, \bfpartial_{\bfheta},\, \bfpartial_{\bfhsigma},\,  \bfpartial_{\bfceta},\, \bfpartial_{\bfcsigma}\}$ 
	and $\{ \bfpartial_{\bftheta_Y}, \bfpartial_{\bfheta},\, \bfpartial_{\bfhsigma},\,  \bfpartial_{\bfceta},\, \bfpartial_{\bfcsigma}\}$ satisfy
\begin{equation}
[\bfpartial_{\bftheta_{X_i}},\bfpartial_{\bftheta_{Y_j}}] = 0.
\EqTag{4Adapted4}
\end{equation}
\end{Theorem}

	Equations \EqRef{4Adapted3} imply that
\begin{equation}
	[\bfpartial_{\bftheta_{X_i}},\bfpartial_{\bftheta_{X_j}}]= -C^k_{ij} \bfpartial_{\bftheta_{X_k}}, 
	\quad [\bfpartial_{\bftheta_{Y_i}},\bfpartial_{\bftheta_{Y_j}} ]= C^k_{ij} \bfpartial_{\bftheta_{Y_k}} .
\EqTag{VessiotBracket}
\end{equation}
	and therefore the constants $\bfC$ are the structure constants for a real Lie algebra.  Any  pair of  
	0-adapted coframes satisfying  \EqRef{4Adapted3} and \EqRef{4Adapted4} is  said to be  a pair of {\deffont 4-adapted coframes}. 

\begin{Remark} 
\StTag{Rem4Adapted} 
	If  $ \{\,\bftheta, \,\bfheta,\, \bfhsigma,\,  \bfceta,\, \bfcsigma\,\}$  is any 0-adapted coframe, then 4-adapted coframes are
	constructed  by taking $\bfthetaX, \bfthetaY \in \spn \{\, \bftheta,\bfheta, \bfceta \, \}$ and keeping the $\bfheta,\, \bfhsigma,\,  \bfceta,\, \bfcsigma\,$ unaltered.
	Consequently, if  $\CalI$ is Darboux integrable linear Pfaffian system with independence condition \EqRef{omegaIC} 
	then the structure equations \EqRef{4Adapted3} for a pair of 4-adapted coframes can be written as
\begin{equation}
\begin{aligned}
d \bfhsigma& = 0, \quad d\bfheta =  \bfhA \, \bfhtau \wedge \bfhomega + \bfhG \, \bfheta\wedge \bfhsigma  , \
	\quad d \bfcsigma = 0, \quad d\bfceta = \bfcF \, \bfctau \wedge \bfcomega+ \bfcH\,  \bfceta\wedge \bfcsigma,
\\	
\extd \bftheta_X
&
	=  \frac12 \bftildeA \, \bfhtau  \wedge \bfhomega   +  \frac12 \bftildeB\, \bfctau  \wedge \bfcomega  + \bftildeG_1 \bfheta \wedge \bfhpi\,  +\bftildeG _2\, \bfceta \wedge \bfcpi         +\frac12 \bfC\,\bfthetaX\wedge \bfthetaX  +  \tilde \bfM \, \bfhpi \wedge \bfthetaX, \\
\extd \bftheta_Y&	=  \frac12\bftildeE \, \bfhtau  \wedge \bfhomega  +  \frac12 \bftildeF \, \bfctau \wedge \bfcomega  
 + \bftildeH_1 \, \bfheta \wedge \bfhpi  +\bftildeH _2 \, \bfceta \wedge \bfcpi
        - \frac12 \bfC\, \bftheta_Y \wedge \bftheta_Y  + \tilde \bfN \, \bfcpi \wedge \bftheta_Y,
  \end{aligned}
\EqTag{4Adaptedic}
\end{equation}
	where $\bfhsigma=(\bfhtau,\bfhomega)$,  $\bfcsigma=(\bfhtau, \bfhomega)$, $\bfhpi=(\bfhtau,\bfhomega, \bfheta)$, and 
	$\bfcpi=(\bfctau,\bfcomega, \bfceta)$.
\end{Remark}

\medskip

	Let $(\CalE, N)$ be  another Darboux integrable differential system  with singular Pfaffian systems 
	$\{\, \hZ, \cZ\, \}$ and 4-adapted coframes
	$\{\,\bftheta_X', \,\bfheta',\, \bfhsigma',\,  \bfceta',\, \bfcsigma'\,\}$
	and $\{\,\bftheta'_Y,\, \bfheta',\, \bfhsigma',\,  \bfceta',\, \bfcsigma'\,\}$.
	If $\phi\:N\to M$ is a smooth constant rank map satisfying 
	$\phi^*(\CalI) \subset \CalE$  and $\phi^*(\hV) \subset \hZ$ and $\phi^*(\cV) \subset \cZ$ then,  by virtue of \EqRef{ZeroC},
\begin{alignat*}{2}
	\phi^*(\{\bfheta,\, \bfhsigma \} ) &\subset  \spn \{\bfheta', \, \bfhsigma' \}  ,
	&\quad \phi^*(\{\bfthetaX , \,\bfheta, \,\bfceta,  \, \bfhsigma \}) &\subset \spn \{\boldsymbol{\theta_{X'}}, \, \bfheta', \, \bfceta', \, \bfhsigma'  \}  ,
\\
	\phi^*(\{\bfceta, \, \bfcsigma\} )  &\subset  \spn \{\bfceta',\,  \bfcsigma' \} ,  
	&\phi^*(\{\bfthetaY,  \, \bfheta,  \, \bfceta,  \, \bfcsigma  \})  &\subset \spn \{ \boldsymbol{\theta_{Y'}}, \,  \bfheta', \, \bfceta', \, \bfcsigma'  \}.
\end{alignat*}
	In  particular, there are matrix-valued functions $\bfR$  and $\bfS$  on $N$ such that  
\begin{equation}
	\phi^*(\bftheta_X)  = \bfR\, \boldsymbol{\theta_{X'}} \mod \, \{ \bfheta', \, \bfceta', \, \bfhsigma' \}
	\quad \text{and}\quad
	\phi^*(\bftheta_Y)  = \bfS\, \boldsymbol{\theta_{Y'}} \mod \, \{ \bfheta', \, \bfceta', \, \bfcsigma' \}.
\EqTag{Pullback4Adapted}
\end{equation}
	The dual vector fields satisfy
\begin{equation}
	\phi_*(\bfpartial_{\bftheta_{X'}}) = \bfR \, \bfpartial_{\bftheta_X} 
	\quad\text{and}\quad
	\phi_*(\bfpartial_{\bftheta_{Y'}}) = \bfS \, \bfpartial_{\bftheta_Y}
\EqTag{VessHomo}
\end{equation}
	and, on account of  the  4-adapted structure equations,  the functions $\bfR$ and $\bfS$ satisfy (see Corollary 4.6 in
	\cite{anderson-fels-vassiliou:2009a}) 
\begin{equation}
	\bfR \, \bfC' = \bfC \bfR\,\bfR,\quad  d\bfR \in \hZ^\infty, \quad \bfS \, \bfC' = \bfC  \bfS\,\bfS, \quad  
	\quad\text{ and}  \quad  d\bfS \in \cZ^\infty.
\EqTag{dRdS}
\end{equation} 
	The first equation in \EqRef{dRdS} is given in components by $R^i_\ell {C'}^\ell_{jk}  = C^i_{\ell m} R^\ell_j  R^m_k$.

	These equations prove that $\phi_{x,*}$,  at each point $x\in N$,  induces a homomorphism 
	from the Lie algebra of vector-fields $\{ \bfpartial_{\bftheta_X'}  \}$  to the Lie algebra of vector fields 
	$\{ \bfpartial_{\bftheta_X}\}$. In particular,  when $\phi$ is  a diffeomorphism, 
	this demonstrates  that  {\it the Lie algebra defined by the  structure constants $\bfC$ in a pair of 4-adapted coframes is
	an invariant of the Darboux integrable system $\CalI$.} The corresponding  abstract Lie algebra we call the {\deffont Vessiot algebra} 
	which we denote by  $\vess(\CalI)$.   We write the Lie algebra homomorphism induced by $\phi$ as
\begin{equation}
	\tilde \phi_x \:\vess(\CalE) \to  \vess(\CalI).
\EqTag{VessHom}
\end{equation}
	Equations \EqRef{ZeroC}   and \EqRef{VessiotBracket} show that the dimension of the Vessiot algebra is 
\begin{equation}
	\dim \vess(\CalI)  = \dim \spn \{\, \bftheta_X\}  =   \dim \spn \{\, \bftheta_Y\}  = \dim M - \rank(\hV^\infty) - \rank(\cV^\infty) .
\EqTag{VessDim}
\end{equation}

	To calculate the Vessiot algebra of a Darboux integrable system, one must  calculate a 4-adapted coframe. Regrettably,
	there is not  at present a more geometric or intrinsic description of this invariant.  We have seen that   integrable extensions and group
	quotients of Darboux integrable systems are Darboux  integrable but in general it seems
	difficult to calculate the Vessiot algebra of
	the extension or quotient  differential system in terms of the Vessiot algebra of the original one.  We return to this issue in Sections 7.3
	and 8.
\begin{Theorem}
\StTag{Fourth142} Let  $(\CalE, N)$ and $(\CalI, M)$ be Darboux integrable differential systems with 
	singular Pfaffian systems  $\{\,\hZ,\cZ\,\}$ and $\{\,\hV, \cV \,\}$  and suppose
	that $\phi :N \to M$  is a smooth map satisfying
\begin{equation}
	\phi^*(\CalI) \subset \CalE, \quad \phi^*(\hV) \subset \hZ \quad\text{and}\quad \phi^*(\cV) \subset \cZ.
\EqTag{Fourth152}
\end{equation}
        Then the induced homomorphism $\tilde \phi _x\colon \vess(\CalE) \to \vess(\CalI)$
        is injective at each point $x$ if and only if
\begin{equation}
	(T^*N)_{\phi,\semibasic} + ( \hZ^\infty \oplus \cZ^\infty ) = T^*N.
\EqTag{Qinj}
\end{equation}
\end{Theorem}

\begin{proof} 
	In order that $\tilde \phi_x$ be injective it is necessary and sufficient that 
	$\phi_{x,*} $,  restricted  to $\text{span}  \{\, \partial_{\bftheta_X'}\}$  be injective. Since  the coframe
	$ \{\,\bftheta_X', \,\bfheta',\, \bfhsigma',\,  \bfceta',\, \bfcsigma'\,\}$ is 0-adapted, we deduce from \EqRef{ZeroC} that 
\begin{equation*}
	\text{span}  \{\, \partial_{\bftheta_X'}\} =  \ann (\hZ^\infty \oplus \cZ^\infty)
\end{equation*}
	and therefore $\tilde \phi_x$ is injective at each point  if and only if
\begin{equation}
	\ker (\phi_{*}) \cap  \ann( \hZ^\infty \oplus \cZ^\infty)  =0.
\EqTag{preinj}
\end{equation}
The dual  of equation \EqRef{preinj} produces equation  \EqRef{Qinj}.
\end{proof}

	The following corollary of Theorem \StRef{H} proves Theorem C in the introduction.

\begin{Corollary}
\StTag{VessiotIntExt}
	Let   $\bfp\:  (\CalE, N) \to (\CalI, M)$ be an integrable extension of  Darboux systems $\CalE$ and $\CalI$.  
	If the pair $(\CalE, \CalI)$ is maximally compatible, then the induced map 
\begin{equation}
	{\bf \tilde \bfp}_x \: \vess(\CalE)  \to \vess(\CalI)
\end{equation}
	is  a Lie algebra monomorphism for each point $x$. The fibre dimension of the integrable extension  $\CalE$ is
	the co-dimension of  the image  ${\bf \tilde \bfp}_x (\vess(\CalE))$ in $\vess(\CalI)$.
\end{Corollary}
\begin{proof}  The fact that $\tilde \bfp_x$ is a monomorphism follows immediately from   {\bf [i]} of Theorem \StRef{H}  and \EqRef{preinj}. For the second part of the corollary we
note,  by equations  \EqRef{IntExtDef}{\bf [i]}, \EqRef{rankZinfty} and \EqRef{VessDim},  that
\begin{equation*}
\begin{aligned}
\dim \vess(\CalE) &= \dim N - \rank(\hZ^\infty)-\rank(\cZ^\infty) = \dim M  - \rank(\hV^\infty) -\rank(\cV^\infty) -  \rank J
\\
	& = \dim \vess(\CalI) - \rank J.
\end{aligned}
\end{equation*}
Since $\tilde \bfp_x$ is a monomorphism,  this last equation  proves that
the co-dimension of the image $\tilde \bfp_x$ is $\rank J$ which is  the fibre dimension of $\bfp:N \to M$.
\end{proof}

\subsection{The Prolongation of Darboux Integrable Systems}

  In this section we show that the prolongation $\CalI^{[1]}$ of a Darboux integrable system $\CalI$ is again Darboux integrable 
and that the two Vessiot algebras are isomorphic. This becomes transparent by writing a pair of $4$-adapted 
coframes for the prolongation $\CalI^{[1]}$  in terms of a pair of $4$-adapted coframe for the original $\CalI$.

\begin{Theorem} \StTag{ProDI}
	Let $\CalI$ be a linear Pfaffian system with independence condition which is Darboux integrable (see \StRef{DCwithIC}).
	Suppose that  $\CalI$ is involutive. Then the prolongation $\CalI^{[1]}$  is also a linear Pfaffian system 
	with independence condition, $\CalI^{[1]}$ is Darboux integrable, and the
	Vessiot algebras $\vess(\CalI^{[1]})$ and $\vess(\CalI)$ are isomorphic.
\end{Theorem}

\begin{proof} 
	Let $\{\bftheta_{\bullet},\,  \bfheta,\,  \bfceta,\,  \bfhtau, \bfhomega, \, \bfctau, \bfcomega\ \}$, where 
	$\bftheta_\bullet = \bftheta_X$ or  $\bftheta_\bullet = \bftheta_Y$,  denote a pair of  $4$-adapted coframes  for 
	$\CalI$ on an open set $U$. The structure equations are  \EqRef{4Adaptedic} which, combined with the  fact that  $\CalI$
	is decomposable, imply that
\begin{equation*}
	\CalI = \langle \bftheta_X, \,  \bfheta,\,  \bfceta,\,   \bftildeA \, \bfhtau  \wedge \bfhomega, \bftildeB \, \bfctau  \wedge \bfcomega, \, 
	\bfhA  \bfhtau \wedge \bfhomega, \,  \bfcF \bfctau \wedge \bfcomega \rangle_{\text{\rm alg}} 
\end{equation*}
and also
\begin{equation*}
	\CalI = \langle \bftheta_Y, \,   \bfheta,\,  \bfceta, \,  \bftildeE \, \bfhtau  \wedge \bfhomega, \,  \bftildeF \, \bfctau  \wedge \bfcomega , \, 
	\bfhA  \bfhtau \wedge \bfhomega, \,  \bfcF \bfctau \wedge \bfcomega \rangle_{\text{\rm alg}} .
\end{equation*}
	We therefore have  
\begin{equation}
 \text{span} \{ \, \bftildeB \, \bfctau  \wedge \bfcomega  \, \}   \subset 
	\text{span} \{ \,  \bftildeF \, \bfctau  \wedge \bfcomega, \,  \bfcF \bfctau \wedge \bfcomega  \, \}
\EqTag{spanB}
\end{equation}
and hence
$$
	\CalI = \langle \bftheta_{\bullet},  \, \bfheta,\,  \bfceta,   \, \bftildeA \, \bfhtau  \wedge \bfhomega,  \, \bftildeF \, \bfctau  \wedge \bfcomega  ,  \, 
	\bfhA  \bfhtau \wedge \bfhomega,  \,  \bfcF \bfctau \wedge \bfcomega \rangle_{\text{\rm alg}}.
$$
The singular Pfaffian systems for $\CalI$  with respect to these decompositions are 
$$
	\hV = \{  \bftheta_{\bullet},  \bfheta,\,  \bfceta, \bfhtau, \bfhomega \} 
	\quad \text{and}\quad
	\cV = \{  \bftheta_{\bullet},  \bfheta,\,  \bfceta, \bfctau, \bfcomega \}.
$$
	
	In order to compute the prolongation (see  \cite{griffiths-jensen:1987a}, \cite{ivy-langsberg:2003}), 
	we first consider  the linear homogeneous equations (with unknowns  $\hSigma^a_c $ and  $\cSigma^\alpha_{ \gamma}$ )
\begin{equation}
	 A^d_{ab} \,\hSigma^a_c\, \homega^b \wedge \homega^c =0,
	\quad \text{and}\quad
	 F ^\delta_{\alpha\beta}\,  \cSigma^\alpha_{ \gamma} \, \comega^\beta \wedge \comega^\gamma =0, 
\EqTag{proeqSimga}
\end{equation}
	where $\bfA=(\bftildeA, \bfhA)$ and $\bfF=(\bftildeF, \bfhF)$.  Implicit in the construction of the  prolongation of $\CalI$ 
	is the assumption  that the solution spaces to these linear systems have constant dimension. 
	Let ${\hS}{}^a_{c, v}$    and ${\cS}{}^\alpha_{ \gamma, \nu}$
	be bases for the solution spaces to \EqRef{proeqSimga},  that is, 
 \begin{equation}
	 A^d_{ab} \hS^a_{c,v}\, \homega^b \wedge \homega^c =0
	\quad \text{and}\quad
	 F ^\delta_{\alpha\beta}\, \cS^\alpha_{\gamma, \nu} \, \comega^\beta \wedge \comega^\gamma =0 
\EqTag{DefSs}
\end{equation}
	so that the general solutions to  \EqRef{proeqSimga} are  $\hSigma^a_c = \hat s^v {\hS}{}^a_{c, v}$ and  
	$\cSigma^\alpha_\gamma =  \check s^\nu {\cS}{}^\alpha_{ \gamma, \nu}$.   The parameters  $\hat s^v$ and $\check s^\nu$ will
	be used as local fibre coordinates for the prolongation manifold $M^{[1]}$ .
	For future use, we note that 
	\EqRef{spanB} and \EqRef{DefSs} imply
\begin{equation}
	 B^a_{\alpha\beta}\, \cS^\alpha_{\gamma, \nu} \, \comega^\beta \wedge \comega^\gamma =0 .
\EqTag{BcS}
\end{equation}

	The chart  for the prolongation space  of $\CalI^{[1]}$  over $U$  is then 
	$U^{[1]} = U \times \real[\hat s^v] \times \real[\check s^\nu]$ with $\pi:U^{[1]}\to U$ and 
\begin{equation}
	\CalI^{[1]} = \langle\,  \bfheta,\, \bfceta,\, \bftheta_{\bullet},\,  \bftheta_1, \,  \bftheta_2 \,\rangle_{\text{\rm diff}},
\EqTag{PRI}
\end{equation}
where
\begin{equation}
	\theta^a_1=\htau^a- \hat s^v\hS^a_{v,c}\, \homega^c \quad\text{and}\quad
	\theta_2^\alpha = \ctau^\alpha -  \check s^\nu   \cS^\alpha_{\gamma, \nu}\, \comega^\gamma.
\EqTag{T12}
\end{equation}

	In order to show that $\CalI^{[1]}$ is decomposable we begin by taking the exterior derivative of the 
	equations in \EqRef{4Adaptedic}  to get
\begin{equation*}
	d  \bftildeA \equiv 0 \mod \hV,\quad d \bfhA \equiv 0 \mod \hV, \quad 
	d  \bftildeF\equiv 0  \mod \cV ,\quad d \bfhF\equiv 0 \mod \cV .
\end{equation*}
	Accordingly,  we may choose the  bases  $\hS^a_{c, v} $ and $\cS^\alpha_{\gamma, \mu}$ for the solutions to  \EqRef{DefSs} 
	so that
\begin{equation}
	d \hS^a_{c, v} \equiv 0  \mod \hV , \quad\text{and}\quad  
	d  \cS^\alpha_{\gamma, \nu} \equiv 0  \mod \cV .
\EqTag{dSs}
\end{equation}
	It  then follows from  \EqRef{T12} and equations \EqRef{4Adaptedic}  that  the structure equations for $\CalI^{[1]}$ are
\begin{equation}
\begin{aligned}
	& d\bfheta \equiv 0 \mod I^{[1]},\quad d\bfceta \equiv 0  \mod I^{[1]},\quad d\bftheta_{\bullet} \equiv 0  \mod I^{[1]}, 
\\[2\jot]
	& d \theta^a_1 =  - ( \hS^a_{c, v} d \hat s^v 	+ s^v  d\hS^a_{c, v})\wedge \homega^c,  
	 \quad d \check \theta^\alpha_2 = -
	(    \cS^\alpha_{\gamma, \nu} d  \check s^\nu + \check s^\nu   d\cS^\alpha_{\gamma, \nu} ) \wedge \comega^ \gamma  .
\end{aligned}
\EqTag{prStrEq}
\end{equation}
	By condition \EqRef{dSs} the structure equations \EqRef{prStrEq} show that system $\CalI^{[1]}$ is  
	decomposable with singular Pfaffian systems
\begin{equation}
	\hV_1 =\{ \bfheta, \bfceta, \bftheta_{\bullet}, \bftheta_1, \bftheta_2,  \bfhomega, d\hat s ^v  \}
	\quad\text{and}\quad
	\cV_1 =\{ \bfheta, \bfceta,\bftheta_{\bullet}, \bftheta_1, \bftheta_2, \bfcomega, d \check s ^\nu \} .
\EqTag{DPI}
\end{equation}
	
	We  next show that  $I^{[1]}$ is Darboux integrable. 
	According to  Theorem \StRef{kerS},   if $k_a \hS^a_{vc} =0 $ for all $v$ and $c$, then $k_a=0$ and 
	likewise if $k_\alpha \cS^\alpha_{\nu \gamma}  =0 $ for all $\nu$ and $\gamma$ then $k_\alpha = 0$. The first consequence of this is 
	(see also Corollary  \StRef{I1D})
$$
	{\CalI^{[1]}}'|_{U^{[1]}} = \langle \bfheta,\bfceta,\bftheta_{\bullet}\rangle_{\text{\rm diff}},
$$
	from which we may conclude  that ${\CalI^{{[1]}}}^\infty  = 0 $. 
	The second consequence  of    Theorem \StRef{kerS}  (see also the argument given in  Corollary  \StRef{I1D})  
	is that it  allows us to compute the derived systems $\hV_1'$ and $\cV_1'$. 
	Using equations \EqRef{T12}, \EqRef{dSs},  \EqRef{prStrEq}, the fact that $\bfhomega$ and $\bfcomega$ are closed, 
	and  Remark B.3,  we deduce that
\begin{equation}
	\hV_1'= \{ \bfheta, \bfceta, \bftheta_{\bullet}, \bftheta_1,  \bfhomega, d\hat s ^v  \} = \pi^*(\hV) + \{ d\hat s^v\} 
	\quad\text{and}\quad 
	\cV_1'=\{ \bfheta, \bfceta,\bftheta_{\bullet}, \bftheta_2, \bfcomega, d \check s ^\nu \}=\pi^* (\cV)+\{ d \check s ^\nu \}
\EqTag{PDIderived}
\end{equation}
	from which we see that $\hV_1'/\pi = \hV$  and $\cV_1'/\pi = \cV$.  Finally, equation \EqRef{PDIderived} implies that (see also \cite{anderson-fels:2012a}) 
\begin{equation}
\begin{aligned}
	\hV_1^\infty &= \pi^*(\hV^\infty) +  \spn\{ d\hat s^v\} =  \spn \{ \bfheta, \bftheta_1, \bfhomega, d \hat s^v  \}
	\quad \text{and}\quad  
\\
	\cV_1^\infty  &= \pi^*(\cV^\infty) + \spn\{ d\check s^\nu \} = \spn\{  \bfceta, \bftheta_2, \bfcomega, d \check s ^\nu \}.
\EqTag{V1infty}
\end{aligned}
\end{equation}
	Together   \EqRef{DPI} and \EqRef{V1infty} imply that  the hypotheses of Theorem \StRef{DPD} are satisfied  
	and so $\CalI^{[1]}$ is Darboux integrable. 

	To show 	that $\vess(\CalI^{[1]})$ and $\vess(\CalI)$ are isomorphic we begin by noting that by \EqRef{VessDim},  
	equation  \EqRef{V1infty} shows that $\dim \vess(\CalI^{[1]}) = \dim \vess(\CalI)$.
	Moreover, on the open set $U^{[1]} $
\begin{equation}
	T^*U^{[1]}_{\pi,sb} = \spn \{ \bftheta_{\bullet}, \bfheta, \bfceta,\bfhtau,\bfctau,\bfhomega, \bfcomega\}
\end{equation}
	and $d\hat s^v , d \check s^\nu \in \hV_1^\infty\oplus \cV_1^\infty $.
	Thus $\pi^* (\CalI) \subset \CalI^{[1]},\  \pi^* (\hCalV) \subset \hCalV_1,\ \pi^* (\cCalV) \subset \cCalV_1$ and condition
	\EqRef{Qinj} of Theorem   \StRef{Fourth142} is satisfied. Since $\vess(\CalI^{[1]})$ and $\vess(\CalI)$ are the same dimension,
	 Theorem \StRef{Fourth142} therefore 
	implies these Lie algebras are isomorphic.
\end{proof}

\begin{Remark}
	The 4-adapted coframes for a Darboux integrable, linear involutive Pfaffian system $\CalI$  with independence condition naturally  
	lift to give 4-adapted coframes for $\CalI^{[1]}$. Indeed, from \EqRef{DPI}  and  \EqRef{V1infty}, we see that 
\begin{equation*}
	\hV^\infty_1 \cap \cV_1 = \spn \{ \bftheta_1, \bfheta \} \quad\text{and}\quad \cV^\infty_1 \cap \hV_1 = \spn \{ \bftheta_2, \bfceta \}
\end{equation*}
	so that, in accordance with the notation introduced in \EqRef{ZeroC} (see also \EqRef{ZeroC1}) we let
\begin{equation}
	 \bfheta_1  = ( \bftheta_1, \bfheta) \quad \text{and} \quad  \bfceta_1  = (\bftheta_2, \bfceta).
\EqTag{bfteta1}
\end{equation}
	Also, set
\begin{equation}
	\bfhtau_1 = (d\hat s ^v) \quad \text{and}\quad \bfctau_1= (d \check s^\nu )  
\end{equation}
	so that the  coframes constructed in the proof of Theorem  \StRef{ProDI} become
\begin{equation}
	\{\, \bftheta_{\bullet},\,\bftheta_1,\, \bftheta_2,\,   \bfheta,\,  \bfceta,\,  \bfhomega, \, \bfcomega, \,  d\hat s ^v, \, d \check s^\nu \}
= 
	\{\, \bftheta_{\bullet},   \bfheta_1,\,  \bfceta_1,\,  \bfhomega, \, \bfcomega, \, \bfhtau_1, \, \bfctau_1 \}
\EqTag{Pr4Adapted4}
\end{equation}
and, in view of  \EqRef{DPI} and \EqRef{V1infty},
\begin{equation}
	I^{[1]}  =   \spn \{\ \bftheta_{\bullet},\,  \bfheta_1,\,  \bfceta_1 \, \} \quad\text{and}\quad
\begin{aligned}
	\hV =  \spn \{\ \bftheta_{\bullet},\,  \bfheta_1,\,  \bfceta_1,\,  \bfhtau_1, \,\bfhomega \}, & \quad \hV^\infty = \spn \{\ \bfheta_1,\,  \bfhtau_1,\, \bfhomega \ \}, 
\\
	\cV=  \spn \{\ \bftheta_{\bullet},\,  \bfheta_1,\, \bfceta_1,\,  \bfctau_1, \, \bfcomega  \}, & \quad \cV^\infty = \spn \{ \ \bfceta_1,\,  \bfctau_1,\, \bfcomega \ \} .
\end{aligned}
\EqTag{ZeroC3}
\end{equation}
	This clearly shows that both coframes  in \EqRef{Pr4Adapted4} are 0-adapted 
	(see  \EqRef{ZeroC}) for the Darboux integrable system $\CalI^{[1]}$ 
	on the prolongation space $U^{[1]}$.

	It is not difficult to calculate the structure equations for \EqRef{Pr4Adapted4} and to match these structure equations 
	against those in \EqRef{4Adaptedic} which characterize 4-adapted coframes for Darboux integrable linear Pfaffian systems
	with independence condition.
	In keeping with the notation in \EqRef{4Adaptedic}  we write
\begin{equation}
	\bfhsigma_1 =(\, \bfhtau_1, \,  \bfhomega\, ),  \quad \bfcsigma_1 =(\, \bfctau_1, \, \bfcomega \, ), 
	\quad \bfhpi_1=(\, \bfhtau_1, \,  \bfhomega, \, \bfheta_1\, ), \quad \bfcpi_1=(\, \bfctau_1,\,  \bfcomega, \, \bfceta_1 \, ).
\end{equation}
	The forms $\bfhtau_1$ and  $\bfhomega$ are closed so that   $d(\bfhsigma_1 ) = 0$ 
	and the first equation in \EqRef{4Adaptedic} is established.
	Since $\bfheta_1  = ( \bftheta_1, \bfheta)$ and  $\bfeta$ satisfies the structure equation in  \EqRef{4Adaptedic}
	 we need only calculate $d(\bftheta_1)$   in terms of the coframe \EqRef{Pr4Adapted4} in order to prove the second structure equation in the first line of  \EqRef{4Adaptedic}.
	From equation  \EqRef{dSs}, we have $d \hS^a_{c, v}  \in \hV^\infty_1$ and  therefore $d \hS^a_{c, v}$ can
	be expressed as a linear combination of the
	forms  $ \{ \bfheta_1,\,  \bfhtau_1, \bfhomega \}$. It then follows from \EqRef{prStrEq} that	
\begin{equation}
	d \theta^a_1 \in 
	\spn\{ \bfhtau_1 \wedge \bfhomega, \,  \bfheta_1  \wedge    \bfhomega,  \,    \bfhomega \wedge \bfhomega  \}.
\end{equation}
	But the terms  $ \bfhomega \wedge \bfhomega$  are torsion terms which cannot appear by the assumption of involutivity and thus
\begin{equation*}
	d \theta^a_1 \in  \spn\{ \,  \bfhtau \wedge \bfhomega,  \,  \bfheta_1  \wedge    \bfhomega \}.
\end{equation*}
	This proves the second structure equation in the first line of  \EqRef{4Adaptedic}. The remaining 
	structure equations in the first line of  \EqRef{4Adaptedic} are similarly proved.
	
	To check the  structure equations in the second line of  \EqRef{4Adaptedic}, we simply have to re-write the structure equation
\begin{equation}
\extd \bftheta_X
	=  \frac12 \bftildeA \, \bfhtau  \wedge \bfhomega   +  \frac12 \bftildeB\, \bfctau  \wedge \bfcomega  + \bftildeG_1 \bfheta \wedge \bfhpi\,  	+	\bftildeG _2\, \bfceta \wedge \bfcpi         
	+	\frac12 \bfC\,\bfthetaX\wedge \bfthetaX  +  \tilde \bfM \, \bfhpi \wedge \bfthetaX
\EqTag{dthetaX4}
\end{equation}
	on  $U$ in terms of the  coframe  \EqRef{Pr4Adapted4} on $U^{[1]}$.  On account of  \EqRef{DefSs} the first terms in  
	\EqRef{dthetaX4} are
\begin{equation*} 
	\bftildeA \, \bfhtau  \wedge \bfhomega   
	=  \bftildeA \, ( \bftheta_1   +\bfS \bfhomega)  \wedge \bfhomega = \bftildeA \,  \bftheta_1\wedge \bfhomega \\
	\in \spn \{ \bfheta_1 \wedge   \bfhomega\}  \subset  \spn \{ \bfheta_1 \wedge   \bfhpi_1\}. 
\end{equation*}
	Similarly, by \EqRef {BcS},  $\bftildeB\, \bfctau  \wedge \bfcomega  \in  \spn \{ \bfceta_1 \wedge   \bfcpi_1\}$.  It then
	follows that  
\begin{equation}
\extd \bftheta_X
	=   \bftildeG_{11} \bfheta_1 \wedge \bfhpi_1\,  +\bftildeG _{21}\, \bfceta_1 \wedge \bfcpi_1 +\frac12 \bfC\,\bfthetaX\wedge \bfthetaX 
	+  \tilde \bfM \, \bfhpi \wedge \bfthetaX
\EqTag{dthetaXlifted}
\end{equation}
	which clearly  matches with the structures equations for a 4-adapted coframe.  The structure equations for
	$\extd \bftheta_Y$ in terms of the coframe on $U^{[1]}$ are similar  and we hence conclude  that  the 
	coframes in \EqRef{Pr4Adapted4} 
	are 4-adapted. Since the coefficients $\bfC$ in \EqRef{4Adaptedic} and  \EqRef{dthetaXlifted} are identical, these calculations 
	also show directly that the Vessiot algebras  $\vess(\CalI^{[1]})$ and $\vess(\CalI)$ are isomorphic.\qed
\end{Remark}


\subsection{ 4-Adapted Coframes and Vessiot Algebras  for  the Quotient Construction}

	In this section we show how to obtain the 4-adapted coframes for the Darboux integrable system 
	$\CalI = (\CalK_1 +\CalK_2)/ G_\diag$ directly from the adapted coframes for each $\CalK_a$  
	introduced in Theorem \StRef{BigLemma}. 
	This construction proves that  $\vess(\CalI) = \lieg$, where $\lieg$ is the Lie algebra of $G$, and 
	will be needed for the proof of  Theorem  \StRef{UIE}.  We assume the hypotheses of  Theorem \StRef{Dreduce} and, in particular, 
	that  the group $G$ acts regularly on $M_a$. We break the construction of these coframes into a number of steps.

\medskip
\noindent
{\bf Step 1.}  Fix a point $p \in M$ and  pick points $p_a \in M_a$ with $\bfq_{G_{\diag}}(p_1,p_2)=p$. 
	Then pick  $G$-invariant open sets $U_a \subset M_a$, with $p_a \in U_a$, and open sets $\barU_a \subset M_a/G $ on which 
	we have local trivializations 
\begin{equation}
		\Phi_a:U_a\to \barU_a\times G, \quad\text{with}\quad \Phi_a=(\bfq_G^a,\phi_a)  
		\quad\text{and}\quad \phi_a(p_a)=e,
\EqTag{Triv}
\end{equation}
	 where  $\bfq^a_G:M_a\to M_a/G$  are the canonical quotient maps and the $\phi_a$  are $G$-equivariant. Let
	 $\{ {\boldsymbol  \bftheta_1,\,  \bfeta_1, \, \bfsigma_1 }\}$ on $U_1$
	and $\{ \,{\boldsymbol \bftheta_2,\,  \bfeta_2,\, \bfsigma_2} \}$ on $U_2$ be coframes  
	satisfying conditions {\bf [i]-[v]} of Theorem \StRef{BigLemma}. Then we have
\begin{equation}
	K_ 1^1|_{U_1}=  \spn \{\,  \bftheta_1,\,  \bfeta_1\, \}
	\quad\text{and}\quad
	K_ 2^1|_{U_2}= \spn  \{\, \bftheta_2,\,  \bfeta_2\, \} \
\EqTag{K1K2basis}
\end{equation}
	and, by  parts {\bf [ii]-[iii]} of Theorem  \StRef{BigLemma},
\begin{equation}
\begin{alignedat}{4}
	\bftheta_1(\bfX) &= \bf{1},
	&\quad   \bfeta_1(\bfX)  &= \bf{0},
	&\quad	\bfsigma_1(\bfX) &= \bf{0},
\\
	{\bftheta}_2(\bfY) & = \bf{1} ,
	&\quad  {\bfeta}_2(\bfY)   & = \bf{0},
	&\quad  {\bfsigma}_2(\bfY) & = \bf{0} .
\end{alignedat}
\EqTag{LLAdapted}
\end{equation}
	Here $\bfX$ denotes the basis of infinitesimal generators for the action of  $G$ on $U_1$ 
	and $\bfY$ the basis of infinitesimal generators for $G$ on $U_2$. 
	The structure equations for $\bfX$ and $\bfY$ are the same.

	For ease of notation we will identify the forms  $\{ {\boldsymbol  \bftheta_1,\,  \bfeta_1, \, \bfsigma_1 }\}$
	on $U_1$ with their pullbacks by $\pi_1$ to $U_1 \times U_2$ and the forms
	 $\{ \,{\boldsymbol \bftheta_2,\,  \bfeta_2,\, \bfsigma_2} \}$ on $U_2$ with their pullbacks by $\pi_2$ to $U_1 \times U_2$.
	Thus 
\begin{equation}
	\{ {\boldsymbol  \bftheta_1,\,  \bfeta_1, \, \bfsigma_1,\, \bftheta_2,\,  \bfeta_2,\, \bfsigma_2} \}
\end{equation}
	defines  a coframe on  $U_1\times U_2$. Let  $\delta:G\times(M_1\times M_2) \to M_1 \times M_2$ 
	denote the diagonal action of $G$ on $M_1\times M_2$, that is, $\delta_g(x_1,x_2)=(\mu_1(g,x_1),\mu_2(g,x_2))$.  
	The vector fields $(\bfX+\bfY)_i = X_i + Y_i$ are then a  basis for the infinitesimal generators for the diagonal action $\delta$.
	Set  $U=(U_1\times U_2)/G_{\diag}$. 
 
\newcommand\bfbareta{\boldsymbol{{\bar \eta}}}
\newcommand\bfbarsigma{\boldsymbol{{\bar \sigma}}}

\bigskip
\noindent
	{\bf Step 2.} Theorem \StRef{BigLemma}, part {\bf[ii]}, states that
	the forms  $ \bfpi_1 = \{\, \bfeta_1, \,\bfsigma_1\, \}$  and  $\bfpi_2 = \{ \bfeta_2,\, \bfsigma_2\,\}$,  
	defined on  $U_1$ and $U_2$ respectively,  are  $G$-basic.
	Accordingly, we can define forms
        $ \bfbareta_1$ and $\bfbarsigma_1$ on $\barU_1 =U_1/G$   and      $\bfbareta_2$ and
	$\bfbarsigma_2$ on $\barU_2 = U_2/G$
	such that
\begin{equation}
	\bfq_{G}^{1,*}(\bfbareta_1) = \bfeta_1,  \quad \bfq_{G}^{1,*}(\bfbarsigma_1) = \bfsigma_1, \quad
	 \bfq_{G}^{2,*}(\bfbareta_2) = \bfeta_2,  \quad  \bfq_{G}^{2,*}(\bfbarsigma_2) = \bfsigma_2.
\end{equation}
	We then use the maps $\bfp_1$ and $\bfp_2$,
	defined by the commutative diagrams \EqRef{H2action} (with $L = G_\diag$), 
	to define 1-forms
	$\bfhpi = \{\bfheta, \bfhsigma\}$ and  $\bfcpi = \{\bfceta, \bfcsigma\}$  on $U$ by
\begin{equation}
	 \bfheta  = \bfp_2^*(\bfbareta_2),   \quad \bfhsigma = \bfp_2^*(\bfbarsigma_2), \quad
	 \bfceta  = \bfp_1^*(\bfbareta_1),  \quad \bfcsigma = \bfp_1^*(\bfbarsigma_1).	
\EqTag{p_pullback}
\end{equation}
	The diagrams \EqRef{H2action} also show that
\begin{equation}
	\quad
	 \bfq_{G_{\diag}}^*(\bfheta)  = \bfeta_2 ,
	\quad
	 \bfq_{G_{\diag}}^*(\bfhsigma) =  \bfsigma_2,
	\quad
	\bfq_{G_{\diag}}^*(\bfceta) =  \bfeta_1,
	\quad
	\bfq_{G_{\diag}}^*(\bfcsigma) = \bfsigma_1.
\EqTag{qL1}
\end{equation}

\medskip
\noindent
{\bf Step 3.}  We now  use \EqRef{deflambda} and the trivializations  \EqRef{Triv} to 
	define matrix-valued functions $ \bflambda_1,\bflambda_2: U_1 \times U_2\to GL(r,\real)$ by
\begin{equation}
\bflambda_1  = \bflambda\circ \phi_1 \circ \pi_1  \quad\text{and}\quad
\bflambda_2  = \bflambda\circ \phi_2 \circ \pi_2.
\EqTag{defl12}
\end{equation}
	By virtue of  the $G$-equivariance of the $\phi_a$,   \EqRef{lambdap2}   and   \EqRef{Triv} one easily checks that
\begin{equation}
\lambda^i_{a, j}(\delta_g(x_1, x_2)) =  \lambda^k_j(g)\lambda^i_k(\phi_a(x_a)) \quad \text{and} \quad \lambda^i_{a, j}(p_1, p_2) = \delta^i_j.
\EqTag{lambdaprop}
\end{equation}
	Next, we use  \EqRef{defl12} to define 1-forms on $U_1\times U_2$ by
\begin{equation}
\vartheta^i_1  =  \lambda^i_{1,j} (\theta^j_2 -\theta^j_1)\quad\text{and}\quad
\vartheta^i_2  =  \lambda^i_{2,j} (\theta^j_2 -\theta^j_1).
\EqTag{defvart}
\end{equation}
	Clearly, equations \EqRef{LLAdapted} imply that  the forms $\vartheta_1^i$ and $\vartheta^i_2$ are $G_{\diag}$ semi-basic.
	By  part {\bf [iv]} of Theorem \StRef{BigLemma}  and  \EqRef{lambdaprop},  it follows that
	the  forms $\vartheta_a^i$  are $G$-invariant and hence $G$-basic. Consequently
	we can define 1-forms  $\bftheta_X$ and $\bftheta_Y$ on $U$ such that
\begin{equation}
	 \bfq_{G_\diag}^*(\bftheta_X) = {\boldsymbol \vartheta_2}
	\quad\text{and}\quad
	 \bfq_{G_\diag}^*(\bftheta_Y) =  {\boldsymbol \vartheta_1}.
\EqTag{qL2}
\end{equation}
	 For future reference we note,  again by  \EqRef{lambdaprop},   that  $\vartheta^i_1(p_1, p_2) =   \vartheta^i_2(p_1, p_2)$
	and hence
\begin{equation}
\bftheta_X(p) = \bftheta_Y(p) .
\EqTag{matchTheta}
\end{equation}
	The 1-forms $\{\, \bftheta_\bullet, \,\bfheta,\, \bfhsigma,\,  \bfceta,\, \bfcsigma \,\}$,
	 define coframes  on $U$.  We claim that these define 
	the sought-after  4-adapted coframes for  $\CalI = (\CalK_1 + \CalK_2)/G$.

\medskip
\noindent
	{\bf Step 4.}  First we check that the above coframes are properly aligned with the singular Pfaffian systems  $\{\hV, \cV \}$ 
	and their derived flags, in other words,  that they define 0-adapted coframes  (see \EqRef{ZeroC}).  
	To begin,  we immediately deduce, using equations \EqRef{LLAdapted}, that on $U = U_1\times U_2$
\begin{gather} 
\begin{aligned}
	\kern -10pt
	(K_1^1+T^*M_2)_{G_{\diag}, \semibasic} &=  \spn \{\, \bftheta_2 - \bftheta_1,\,  \bfeta_1,\, \bfeta_2, \bfsigma_2 \},\
\\
	(T^*M_1+K_2^1)_{G_{\diag}, \semibasic} &=  \spn\{\, \bftheta_2 - \bftheta_1,\,  \bfeta_1,\,  \bfsigma_1,   \, \bfeta_2 \},
\EqTag{Ksb1}
\end{aligned}
\\
\kern -10pt
(T^*M_2)_{G, \semibasic} = \spn  \{\, \bfeta_2, \,\bfsigma_2 \, \} ,
	\ 
	(T^*M_1)_{G, \semibasic} = \spn \{\,  \bfeta_1,\,  \bfsigma_1   \, \}, \ \quad\text{and}
\EqTag{Ksb2}
\\
(K_1^1 + K_2^1)_{G_{\diag}, \semibasic}  = \spn \{\, \bftheta_2 - \bftheta_1,\,  \bfeta_1,\, \bfeta_2\,\}.
\EqTag{Ksb3}
\end{gather}
	Then, since
\begin{equation}
\begin{aligned}
	\bfq_{G_\diag}^*\ \{\bftheta_\bullet,\,  \bfheta, \, \bfceta,\, \bfhsigma\, \}  &=   \spn\{\, \bftheta_2 - \bftheta_1,\,  \bfeta_1,\, \bfeta_2, \bfsigma_2 \}
	\quad\text{and}
\\
	\bfq_{G_\diag}^*\{\bftheta_\bullet,\,  \bfheta, \, \bfceta ,\, \bfcsigma\, \} &=   \spn \{\, \bftheta_2 - \bftheta_1,\,  \bfeta_1,\,  \bfsigma_1,   \, \bfeta_2 \}
\end{aligned}
\end{equation}
	(see \EqRef{qL1} and \EqRef{qL2}), it  follows from   \EqRef{Asb} and definition \EqRef{WMODG} that 
\begin{equation}
	\hV|_U =  \spn  \{\bftheta_\bullet,\,  \bfheta, \, \bfceta,\, \bfhsigma\, \}
	\quad\text{and}\quad
	\cV|_U = \spn  \{\bftheta_\bullet,\,  \bfheta, \, \bfceta ,\, \bfcsigma\, \}.
\EqTag{ZeroAdapted1}
\end{equation}
	Moreover, the combination of \EqRef{p_pullback},  \EqRef{Ksb2} and  Theorem \StRef{Dreduce}{\bf[iii]} gives
\begin{equation}
	\hV^{\infty}|_U = \bfp_2^*(T^*(U_2/G)) = \spn \{\, \bfheta,\,  \bfhsigma\, \}
	\quad\text{and}\quad
  \cV^{\infty}|_U = \bfp_1^*(T^*(U_1/G)) =  \spn\{\, \bfceta,\, \bfcsigma\,\} .
\EqTag{ZeroAdapted2}
\end{equation}
	Equations \EqRef{ZeroAdapted1}  and \EqRef{ZeroAdapted2}  show that the  coframes
	$\{\, \bftheta_\bullet, \,\bfheta,\, \bfhsigma,\,  \bfceta,\, \bfcsigma \,\}$  are
	0-adapted coframes  for the Darboux pair $\{\, \hV, \cV\, \}$ on $M$. 
	We also  remark  that (see \EqRef{qMdef})
\begin{equation}
\begin{aligned}
	\bfq_{M_1}^*(\hV) &=   (\iota_{M_1})^*(\spn \{\, \bftheta_2 - \bftheta_1,\,  \bfeta_1,\, \bfeta_2, \bfsigma_2 \})   = \spn \{\, \bftheta_1,\,  \bfeta_1\, \}  = K^1_1 
	\quad\text{and}\quad 
\\
	\bfq_{M_2}^*(\cV)  &=  (\iota_{M_2})^*(\spn \{\, \bftheta_2 - \bftheta_1,\,  \bfeta_1,\, \bfeta_2, \bfsigma_1\}) = \spn \{\, \bftheta_2,\,  \bfeta_2\, \}  = K^1_2.
\end{aligned}
\EqTag{Kdiag}
\end{equation}

\bigskip
\noindent
{\bf Step 5.} We  now verify the structure equations in \EqRef{4Adapted3}.
	The  structure equations on the  first line of  \EqRef{4Adapted3} follow directly from \EqRef{BigStr}. 
	To prove the remaining structure equations
	we first recall, from equation \EqRef{BigStr},  that the  structure equations for the 1-forms  $ \bftheta_1$  and $ \bftheta_2$  are
\begin{equation}
	d \bftheta_1 = \bfA\,\bfpi_1 \wedge \bfpi_1 - \frac{1}{2}\bfC\,  \bftheta_1  \wedge \bftheta_1
	\quad \text{and}\quad
	d \bftheta_2 = \bfB\, \bfpi_2 \wedge \bfpi_2 - \frac{1}{2}\bfC\,  \bftheta_2  \wedge \bftheta_2.
\EqTag{StrEq2}
\end{equation}
	where  $\bfpi_1  = \{ \bfeta_1, \bfsigma_1\}$  and $\bfpi_2  = \{ \bfeta_2, \bfsigma_2\}$. The coefficients $\bfA$ 
	are defined on $U_1$, the coefficients $\bfB$ are defined on $U_2$ 
	and the $\bfC$ are the structure constants for the Lie algebra of vector fields  $\bfX$ (or $\bfY$).
	Secondly, on account of    \EqRef{dlambdaG} and  \EqRef{deftheta}, we have that
\begin{equation*}
	 d \lambda^i_{1,j} =  \pi_1^* \bigl(\phi_1^*( d \lambda^i_{j})\bigr)  =   \pi_1^* \bigl (\phi_1^*(\lambda^i_kC^k_{\ell j} \, \tau^\ell) \bigr) 
	= \lambda^i_{1,k} C^k_{\ell j}\, \theta^\ell_1  + M^i_{1\alpha j} \sigma_1^\alpha,
\end{equation*}
	where $M^i_{1\alpha j} \in C^\infty(U_1)$. By virtue  of \EqRef{BigStr} and  this last equation 
	we calculate the structure equations for the forms ${\boldsymbol \vartheta_1}$ in \EqRef{defvart} to be
\begin{align*}
	d{\boldsymbol \vartheta_1}  & = d \bflambda_1 \wedge (\bftheta_2 - \bftheta_1) +  \bflambda_1 \, d(\bftheta_2 - \bftheta_1)
\\
	&=	(\bflambda_1\bfC \bftheta_1 +  \bfM \bfpi_1)\wedge ( \bftheta_2 - \bftheta_1)
	+ \bflambda_1  ( - \frac{1}{2}\bfC\,  \bftheta_2  \wedge \bftheta_2
	+   \frac{1}{2}\bfC\,  \bftheta_1  \wedge \bftheta_1  -  \bfA\,\bfpi_1 \wedge \bfpi_1 +\bfB\, \bfpi_2 \wedge \bfpi_2).
\end{align*}
	All the terms on the right-hand side of this equation which  involve the 2-forms $\bftheta_a\ \wedge\bftheta_b$ combine to
	give $ -\frac12 \bfC\,{\boldsymbol \vartheta_1}  \wedge{\boldsymbol \vartheta_1}$  and
	thus
\begin{equation}
	d{\boldsymbol \vartheta_1}  =   -\bflambda_1\bfA\,\bfpi_1 \wedge \bfpi_1  +  \bflambda_1\bfB\, \bfpi_2 \wedge \bfpi_2
	-\frac12 \bfC\,  {\boldsymbol \vartheta_1} \wedge {\boldsymbol \vartheta_1}   +  \bfM\, \bflambda_1^{-1}\bfpi_1\wedge {\boldsymbol \vartheta_1}.
\EqTag{domega1}
\end{equation}
	A similar equation for  $d {\boldsymbol \vartheta_2}$ holds (except for a change in the sign of the term containing
	the structure constants $\bfC$). Finally equations \EqRef{domega1} (and the counterpart for  $d{\boldsymbol \vartheta_2}$) 
	yield the structure equations
	\EqRef{4Adapted3}.

	Equations \EqRef{4Adapted3} immediately imply that the dual vector fields $\partial_{\theta^i_X}$ and $\partial_{\theta^j_Y}$
satisfy 
\begin{equation}
	[\,\partial_{\theta^i_X},\, \partial_{\theta^j_X}\,] =  -C^k_{ij} \partial_{\theta^k_X}  \quad \text{and}\quad
	[\,\partial_{\theta^i_Y},\, \partial_{\theta^j_Y}\,] =   C^k_{ij} \partial_{\theta^k_Y} .
\end{equation}
 
\bigskip
\noindent
{\bf Step 6.} Lastly,   we check  that  the dual vector fields   $\partial_{\theta^i_{X}}$ to 
	the coframe $\{\bftheta_X, \,\bfheta,\, \bfhsigma,\,  \bfceta,\, \bfcsigma\}$  on $M$   and  the dual vector fields 
	$\partial_{\theta^j_{Y}}$  to $\{\bftheta_Y,\, \bfheta,\, \bfhsigma,\,  \bfceta,\, \bfcsigma\,\} $ are related to 
	the infinitesimal generators  $X_i$  and $Y_j$  for the action of $G$  on $M_1$ and $M_2$ by
\begin{equation}
	\bfq_{M_1 *}(X_i)  = -\partial_{\theta^i_{X}} \quad\text{and}\quad  	\bfq_{M_2*}(Y_j)  =  \partial_{\theta^j_{Y}}
\EqTag{XYpush}
\end{equation}
and satisfy
\begin{equation}
	 [\, \partial_{\theta^i_{X}},\  \partial_{\theta^j_{Y}} \, ] = 0.
\EqTag{thetaCom}
\end{equation}
\
	For  \EqRef{XYpush} it suffices to note, by \EqRef{LLAdapted},
	that  $X_i = \partial_{\theta_1^i}$ and $Y_j= \partial_{\theta_2^j}$ and  then to use
	equations  \EqRef{qL1} and \EqRef{qL2} to calculate  the Jacobians  $\bfq_{M_1*}$ and $\bfq_{M_2*}$  in terms 
	of the dual bases  
	$\{ \bfpartial_{\bftheta_1},\,  \bfpartial_{\bfeta_1}, \,\bfpartial_{\bfsigma_1}\}$ 
	on $M_1$ and 
	$\{\bfpartial_{\bftheta_2},\,  \bfpartial_{\bfeta_2}, \,\bfpartial_{\bfsigma_2}\}$ on $M_2$. 

	To prove  \EqRef{thetaCom},  let  $\mu^k_{1j} = (\lambda^{-1}_1)^k_j$   and  $\mu^k_{2j} = (\lambda^{-1}_2)^k_j$.  
	We then  use   \EqRef{qL2}
	to calculate 
\begin{align*}
	\theta_X^i\bigl(\bfq_{G_\diag*}( \mu^k_{2j}\, \partial_{\theta^k_2})\bigr) 
	&= \bigl(\bfq_{G_\diag}^*(\theta_X^i)\bigr) ( \mu^k_{2j}\, \partial_{\theta^k_2})= 	\delta^i_j 
	 \quad \text{and} \
\\[1\jot]
	\theta_Y^i\bigl(\bfq_{G_\diag*}( \mu^k_{1j}\, \partial_{\theta^k_1})\bigr)  
	&= \bigl(\bfq_{G_\diag}^*(\theta_Y^i)\bigr)( \mu^k_{1j}\, \partial_{\theta^k_1}) = -	\delta^i_j.
\end{align*}
	These equations, together with  \EqRef{qL1},  lead to 
\begin{equation}
	\bfq_{G_\diag*}( \mu^k_{2j}\, \partial_{\theta^k_2}) = \partial_{\theta^j_X} 
	\quad\text{and}\quad 
	\bfq_{G_\diag*}( \mu^k_{1j}\, \partial_{\theta^k_1}) = -\partial_{\theta^j_Y} 
\end{equation}
	from which it then follows, because the $\mu^k_{1j}$ are functions on $M_1$ and the  $\mu^k_{2j}$ are functions on 
	$M_2$,  that
\begin{equation}
	 [\, \partial_{\theta^i_{X}},\  \partial_{\theta^j_{Y}} \, ] 
	=  -\bfq_{G_\diag*}\bigr([\, \mu^k_{2i}\, \partial_{\theta^k_2}, \,  \mu^\ell_{1j}\, \partial_{\theta^\ell_1} \,]\bigr) =  0.
\end{equation}	
	This completes the proof that the coframes $\{\, \bftheta_\bullet, \,\bfheta,\, \bfhsigma,\,  \bfceta,\, \bfcsigma \,\}$
	are 4-adapted  coframes  for the Darboux integrable system $\CalI  = (\CalK_1 + \CalK_2)/G$.

	Since the structure constants  $\bfC$ in \EqRef{domega1} are the structure constants for the Lie algebra of $G$, this theorem 
	immediately implies that {\it  the Vessiot algebra for the Darboux integrable system $(\CalK_1 +\CalK_2)/G_\diag$
	is  the Lie algebra of $G$.}
	More generally, Theorems \StRef{Dreduce} and Lemma  \StRef{CDPSI} combine to give the following result.
\begin{Theorem}
\StTag{VessAlgForQuot}
	 Let  $L \subset G\times G$ and let $A_1$ and $A_2$ be the subgroups of $L$ defined by \EqRef{Asubgroup}.
	Under the hypotheses of Theorem \StRef{Dreduce}, the Vessiot algebra of the Darboux integrable system $(\CalK_1 +\CalK_2)/L$
	is the Lie algebra of Lie group  $\tilde L = L/(A_1 \times A_2)$.
\end{Theorem}

%

\section{Quotient Representations for Maximally Compatible Integrable Extensions}
	In  Theorem \StRef{Dreduce} 
	we showed how Darboux integrable systems can be constructed using the group quotient of pairs of differential systems. 
	It is a  remarkable fact,  established in \cite{anderson-fels-vassiliou:2009a}, that the converse is true  locally, that is,   
	every Darboux integrable system can be locally realized as the quotient of a pair of differential systems 
	with a common symmetry group.  The precise formulation of this result is as follows.
	
\begin{Theorem} 
\StTag{Intro4} 
	Let $(\CalI, M)$ be a Darboux integrable differential system with singular Pfaffian systems $\hV$ and  $\cV$.
	Fix a  point $x_0 $ in $ M$ and let

\noindent
{\bf [i]}  $M_1$ and $M_2$ be the maximal integral manifolds of  $\hV^\infty$ and $\cV^\infty$ through $x_0$, and

\smallskip
\noindent
{\bf [ii]}  $\CalK_1$ and $\CalK_2$ be the restrictions of the singular systems $\hCalV$ and  $\cCalV$ to $M_1$ and $M_2$.

\smallskip
\noindent
	Then there are open sets $U \subset M$, $U_1 \subset M_1$, $U_2 \subset M_2$, each containing $x_0$, 
	and a local action of a  Lie group $G$ on $U_1$ and $U_2$ which
	satisfies the hypotheses of  Theorem \StRef{Dreduce} 
	and such that
\begin{equation}
	U = (U_1 \times U_2)/ G_\diag \quad \text{and} \quad 
	\CalI|_U   = (\CalK_1|_{U_1} + \CalK_2|_{U_2})/G_\diag.
\EqTag{LocQuotRep}
\end{equation}
\end{Theorem}
	The Lie algebra of $G$ coincides with the  Vessiot  algebra of $\CalI$ and the
	local actions of $G$ on $U_1$ and $U_2$ are given  by the  restrictions of the  Lie algebras of  vectors fields  
	$\partial_{\bfthetaX}$ and $\partial_{\bfthetaY}$, dual to the 4-adapted coframes defined by Theorem \StRef{AFV}.  
	When these local actions can be extended to global actions of $G$ on $M_1$ and $M_2$ and when  $M = (M_1 \times M_2)/G_\diag$,
	it then follows that 
\begin{equation}
	\CalI  = (\CalK_1 + \CalK_2)/G_\diag.
\EqTag{GloQuotRep}
\end{equation}

	As a consequence of our proof of Theorem D we shall prove that the quotient representation  \EqRef {GloQuotRep}  of  $\CalI$ 
	is unique. We shall refer to \EqRef{LocQuotRep}  and  \EqRef {GloQuotRep} as 
	{\deffont the local/global canonical quotient representation for a Darboux integrable differential system $\CalI$}.

\subsection{Uniqueness of Maximally Compatible Integrable Extensions}
	In this section we tie together results on integrable extensions of  Darboux integrable systems  in Section 5 
	with the group reduction results of Section 6 to prove Theorem  D. 
	We start with a mapping  $\bfp\:  (\CalE, N) \to (\CalI, M)$  which defines  $ (\CalE, N)$  
	as  an integrable extension of  a Darboux integrable system $(\CalI, M)$.    
By Theorem \StRef{DIext}, we know that $\CalE$
	is always Darboux integrable.   We
	postulate that both $\CalE$ and $\CalI$ admit global canonical quotient representations (see Theorem \StRef{Intro4}),
\begin{equation}
\begin{aligned}
	(\CalE, N) &= \bigl((\CalL_1 +  \CalL_2)/{H_\diag} ,\  (Q_1  \times Q_2)/H_\diag \bigr)   
	\ \text{ and} 
\\ 
	(\CalI, M)  &= \bigl((\CalK_1 + \CalK_2)/ {G_\diag} ,\ (P_1  \times P_2)/G_\diag \bigr),
\end{aligned}    
\EqTag{globCon}
\end{equation}
 	where  $(\CalL_a, Q_a, H)$ and $(\CalK_a, P_a, G)$ satisfy the hypotheses of Theorem  \StRef{Dreduce}.
	The singular systems for these differential systems are (see  \EqRef {Wsing} and \EqRef{WMODG})
\begin{equation}
\begin{alignedat}{4}
	\hU &= L_1^1  \oplus  T^*Q_2,
 &\quad  
	\cU &= T^*Q_1 \oplus L_2^1, 
&\quad  
	\hZ &= (L_1^1  \oplus  T^*Q_2)/H_\diag,
&\quad
	\cZ &= ( T^*Q_1 \oplus L_2^1 )/H_\diag,
\\
	\hW &=  K_1^1  \oplus  T^*P_2,
&\quad 
	\cW &= T^*P_1 \oplus K_2^1, 
&\quad 
	\hV &= (K_1^1  \oplus  T^*P_2)/G_\diag,
&\quad
	\cV &=  (T^*P_1 \oplus K_2^1 )/G_\diag
\end{alignedat}
\EqTag{SS}
\end{equation}
	with derived systems   \EqRef{WMODG3}.  In addition, we postulate that  the pair $(\CalE, \CalI)$ is  maximally 
	compatible with respect to the singular systems $\hZ$, $\cZ$, $\hV$, $\cV$ in \EqRef{SS}.

	From \EqRef{globCon}  we have the  following diagram
\begin{equation}
\begin{gathered}
\begindc{\commdiag}[3]
\obj(0, 0)[L]{$(\CalL_1 +\CalL_2,\  Q_1 \times  Q_2)$}
\obj(40, 0)[P]{$(\CalK_1 +\CalK_2,\  P_1 \times P_2)$}
\obj(0, -17)[N]{$(\CalE, N)$}
\obj(40, -17)[M]{$(\CalI, M)\, .$}
\mor{L}{N}{$\bfq_{H_\diag}$}[\atright, \solidarrow]
\mor{P}{M}{$\bfq_{G_\diag}$}[\atleft, \solidarrow]
\mor{N}{M}{$\bfp$}[\atleft, \solidarrow]
\enddc
\end{gathered}
\EqTag{qHqG1}
\end{equation}
	and our goal is to prove that  there is a subgroup  $H' \subset G$  such that  
\begin{equation}
	(\CalE, N) \cong _{\text{\rm loc}}   (\CalK_1 \times \CalK_2, P_1\times P_2) /H'_\diag .
\EqTag{ElocKH}
\end{equation}
	Specifically, we shall prove that there is a globally defined,  local diffeomorphism   $\Psi : N \to  (P_1\times P_2)/H'_\diag$  
	such that $ \Psi^*((\CalK_1 \times \CalK_2)/H'_\diag) = \CalE $.
	This then proves Theorem D. The key to proving \EqRef{ElocKH}   from \EqRef{qHqG1} is to construct 
	a Lie group monomorphism $\phi :H \to G$  and
	$H$ equivariant maps $\psi_a:Q_a\to P_a$ such that $\psi_{a}^*(\CalK_a) = \CalL_a$ from which $\Psi$ is easily defined.

\begin{Theorem}  
\StTag{UIE}
	Let  \EqRef{globCon}, \EqRef{SS} and \EqRef{qHqG1} be given.  Suppose that  the manifolds $Q_a$,  $P_a$  in \EqRef{qHqG1} are connected, that the
	Lie group $H$ is connected  and that the actions of the groups $H$ and $G$ are free and regular 
	on $Q_1$, $Q_2$ and $P_1$, $P_2$ respectively and transverse to  $\CalL_a$ and $\CalK_a$.  Suppose that the  
	differential systems $( \CalE, \CalI )$ are maximally compatible with respect to the singular systems \EqRef{SS}.
	Pick points $(q_1, q_2) \in Q_1\times Q_2$ and $(p_1, p_2) \in  P_1\times P_2$ 
	such that $\bfp \circ \bfq_{H_\diag}(q_1, q_2) =  \bfq_{G_\diag}(p_1, p_2)$ and set (see \EqRef{qMdef}) 
\begin{equation}
	\bfq_{Q_a}  = \bfq_{H_\diag} \circ \iota_{Q_a} : Q_a \to N
	\quad\text{and}\quad
	\bfq_{P_a}  = \bfq_{G_\diag} \circ \iota_{P_a} : P_a \to M.
\end{equation}

	Then there are  globally defined local diffeomorphisms $\psi_a:Q_a\to P_a$ and a Lie group homomorphism $\phi :H \to G$ 
	with injective differential $\phi_*$ such that:

\smallskip
\noindent
{\bf[i]} the following two diagrams 
\begin{equation}
\begin{gathered}
\begindc{\commdiag}[3]
\obj(0, 0)[L]{$Q_1$}
\obj(28, 0)[P]{$P_1 $}
\obj(0, -17)[N]{$N$}
\obj(28, -17)[M]{$M$}
\mor{L}{N}{$\bfq_{Q_1}$}[\atright, \solidarrow]
\mor{P}{M}{$\bfq_{P_1}$}[\atleft, \solidarrow]
\mor{N}{M}{$\bfp$}[\atleft, \solidarrow]
\mor{L}{P}{$\psi_1$}[\atleft, \solidarrow]
\enddc
\end{gathered}
\quad\text{and}\quad
\begin{gathered}
\begindc{\commdiag}[3]
\obj(0, 0)[L]{$Q_2$}
\obj(28, 0)[P]{$P_2$}
\obj(0, -17)[N]{$N$}
\obj(28, -17)[M]{$M$}
\mor{L}{N}{$\bfq_{Q_2}$}[\atright, \solidarrow]
\mor{P}{M}{$\bfq_{P_2}$}[\atleft, \solidarrow]
\mor{N}{M}{$\bfp$}[\atleft, \solidarrow]
\mor{L}{P}{$\psi_2$}[\atleft, \solidarrow]
\enddc
\end{gathered}
\EqTag{psi1psi2}
\end{equation}
	are commutative;

\smallskip
 \noindent
{\bf[ii]}  $\psi_{a}^*(\CalK_a) = \CalL_a$ so that  the differential systems  $\CalK_ a$ and $\CalL_a$  are locally equivalent;

\smallskip
\noindent
{\bf[iii]} the maps $\psi_a$ are $H$ equivariant in the sense that  $\psi_a(h\cdot q) = \phi(h) \cdot \psi_a(q)$ for all $q \in Q_a$ and $h \in H$.

\smallskip
\noindent {\bf[iv]}  Assuming that $H' = \phi(H) \subset G$  is a closed subgroup,   {\bf [i]}, {\bf [ii]} and {\bf [iii]}  together imply that
\begin{equation}
	(\CalE, N) \cong _{\text{\rm loc}}   (\CalK_1 \times \CalK_2, P_1\times P_2) /H'_\diag .
\end{equation}
\end{Theorem}

\begin{proof} The first step is to construct the map $\psi_1$.
	By Theorem \StRef{Dreduce}{\bf[v]} the map $\bfq_{Q_1}\: Q_1 \to N$ is a maximal integral manifold for $\hZ^\infty$ 
	through the point $\bfq_{H_\diag}(q_1,q_2)$. 
	The hypothesis of maximal compatibility implies that  $\bfp :N \to M$ defines $\hZ^\infty$  as an integrable extension of $\hV^\infty$
        (see Definition \StRef{DarbComp}  part {\bf [ii]} and equation \EqRef{DarbAdmis1}).  
	Thus $\bfp$ maps integral manifolds of  $\hZ^\infty$  to integral manifolds  of $\hV^\infty$ and moreover, by \IntExtInt,  
	$\bfp$ maps 	integral manifolds of  $\hZ^\infty$ of maximal dimension
	to integral manifolds  of $\hV^\infty$ of maximal dimension. 
	Therefore  $\bfp \circ \bfq_{Q_1}\: Q_1 \to N$ is an integral manifold   of 
       $\hV^\infty$  through $\bfp(\bfq_{H_\diag}(q_1,q_2))$ of maximum dimension $\dim Q_1$.
	Again, by   Theorem  \StRef{Dreduce}{\bf[v]}, 
	the map $\bfq_{P_1}\: P_1 \to M$ is a maximal integral manifold for $\hV^\infty$  through  $\bfq_{G_\diag}(p_1,p_2)$ and  therefore 
\begin{equation}
	\bfp \circ \bfq_{Q_1}(Q_1) \subset \bfq_{P_1}( P_1)  \quad \text{and hence}\quad  \dim Q_1 = \dim P_1.
\end{equation}
	The fact that $ \bfq_{P_1}$ is one-to-one  (once more, Theorem  \StRef{Dreduce}{\bf[v]})
	then  implies that  $\bfp \circ \bfq_{Q_1}$  factors through $ \bfq_{P_1}$, that is, 
	there is  a  unique  map $\psi_1:Q_1 \to P_1$ such that the first diagram  in \EqRef{psi1psi2} commutes.  
	A basic  result on the factorization of maps through integral manifolds  (\cite{warner:1983a},  page 47)  states 
	that $\psi_1$ is smooth. 
	Finally,  since $\bfp \circ \bfq_{Q_1}$ and  $\bfq_{P_1}$ both have injective differentials, $\psi_{1*}$ is injective 
	and hence $\psi_1$ is a local diffeomorphism. 

	Let $\hJ$ be an admissible bundle for $(\CalE, \CalI)$ satisfying equation \EqRef{DarbAdmis1}. 
	Since $\bfq_{Q_1}: Q_1 \to N$ is  a maximal integral manifold for $\hZ^\infty$, it follows that $\bfq_{Q_1}^*(\hJ) =  0$.
	Equation \EqRef{IEWs1} (with $J=\hJ$) implies 
\begin{equation*}
	\bfq^*_{Q_1}(\hCalZ) = \bfq^*_{Q_1} ( \bfp^*(\hCalV))
\end{equation*}
	while part {\bf[ii]} in Theorem \StRef{Intro4}  states that
\begin{equation}
	\CalL_1 = \bfq_{Q_1}^*(\hCalZ) \quad \text{and}\quad \CalK_1 = \bfq_{P_1}^*(\hCalV).	
\end{equation}
	In conjunction  with the commutativity of  \EqRef{psi1psi2}  this leads to
\begin{equation}
	\psi^*_1(\CalK_1) =  \psi^*_1 ( \bfq_{P_1}^*(\hCalV)) =   \bfq^*_{Q_1} (\bfp^*(\hCalV)) = \bfq^*_{Q_1}(\hCalZ) = \CalL_1.
\end{equation}
	Similar arguments apply to the second diagram in \EqRef{psi1psi2}  and  the proof of parts {\bf [i]} and  {\bf [ii]} are complete.

	The proof of part {\bf [iii]} is actually  quite straightforward  once all the appropriate notation is fixed. 
	To this end, let $X^Q_i$, $Y^Q_i$  be the infinitesimal generators for the action of $G$ on $Q_1$  and $Q_2$ and
	let  $X^P_ r$, $Y^P_r$ be the infinitesimal generators for the action of $G$ on $P_1$  and $P_2$. 
	The structure equations are
\begin{equation*}
	[X^Q_i, X^Q_j] = C^k_{ij} X^Q_k, \quad [Y^Q_i, Y^Q_j] = C^k_{ij} Y^Q_k,
	\quad [X^P_r, X^P_s] = S^t_{rs} X^P_t,   \quad [Y^P_r, Y^P_s] = S^t_{rs} Y^P_t. 
\end{equation*}
	The infinitesimal generators for the action of $H_\diag$  on $Q_1\times Q_2$ and $G_\diag$ on $P_1\times P_2$ are
\begin{equation*}
	V_i = X^Q_i +  Y^Q_i \quad\text{and}\quad  W_r  = X^P_r+   Y^P_r.
\end{equation*}
	
	Since the actions of  $H$ and $G$ are free and transverse to the differential systems  $\CalL_a$ and $\CalK_a$,  
	we can apply Theorem \StRef{BigLemma} to obtain coframes
\begin{equation}
	\{\, \bftheta_{Q_a},\,  \bfeta_{Q_a},\,  \bfsigma_{Q_a}\,  \}
	\quad \text{and}\quad
	\{\, \bftheta_{P_a},\,  \bfeta_{P_a},\,  \bfsigma_{P_a}\,  \}
\EqTag{Uniq44}
\end{equation}
	on $Q_a$ and $P_a$.  Note that  $\bfpartial_{\bftheta_{Q_1}} = \bfX^Q$,    $\bfpartial_{\bftheta_{Q_2}} = \bfY^Q$,
	$\partial_{\bftheta_{P_1}} = \bfX^P$ and  $\partial_{\bftheta_{P_2}} = \bfY^P$.
	Then from these coframes we construct,  as  in  Section 7.3  (Step 2 and Step 3),  4-adapted  local coframes
\begin{equation}
	\{\bftheta_{\bullet N}, \bfheta_N, \bfceta_N, \bfhsigma_N,  \bfcsigma_N\}
	\quad\text{and}\quad
	\{\bftheta_{\bullet M}, \bfheta_M, \bfceta_M, \bfhsigma_M,  \bfcsigma_M\}
\EqTag{Uniq48}
\end{equation}
	on the quotient manifolds $N$ and $M$.  In the present context  equations  \EqRef{XYpush}
	become
\begin{equation}
	 \bfq_{Q_1*} (X^Q_i) =  -\partial_{\theta^i_{X,N}},
	\quad
	\bfq_{Q_2 *} (Y^Q_i) =   \partial_{\theta^i_{Y,N}} ,
	\quad
	\bfq_{P_1*} (X^P_r) =  -\partial_{\theta^r_{X,M}},
	\quad
	\bfq_{P_2 *} (Y^P_r) =  \partial_{\theta^r_{Y,M}}.
\EqTag{Uniq55}
\end{equation}
	
	We now use \EqRef{Uniq55} to prove that the maps $\psi_a$ induce a  common Lie algebra  homomorphism
	$\psi_{a, *}:\Gamma_H \to  \Gamma_G$, that is,
\begin{equation}
	\psi_{1,*}(X^Q_i) = A^r_{i} X^P_r 
	\quad\text{and}\quad 
	\psi_{2,*}(Y^Q_i) = A^r_{i} Y^P_r,
\EqTag{psihomo}
\end{equation}
	where the $A^r_i$ are constants.
	Since the integrable extension $\bfp$ is maximally compatible $\bfp^*(\CalI) \subset \CalE$, $\bfp^*(\hV) \subset \hZ$ and
	$\bfp^*(\cV) \subset \cZ$ and therefore, by virtue of the remarks made in Section 7.1 (see \EqRef{VessHomo}), there
	are functions $ R^s_i $ and   $S^s_i$ on $N$ such that  
\begin{equation}
	\bfp_*( \partial_{\theta^i_{X, N}})  = R^s_i \, \partial_{\theta^s_{X, M}} \quad\text{and}\quad 
	\bfp_*( \partial_{\theta^i_{Y, N}})  = S^s_i  \,\partial_{\theta^s_{Y, M}}.
\EqTag{Uniq56}
\end{equation}
        If the vector fields  $\partial_{\theta^i_{X, N}}$,  $\partial_{\theta^i_{Y, N}}$   in \EqRef{Uniq56} 
	are evaluated at a point $x \in N$, then the functions $R^s_i $ and $S^s_i$ are  evaluated at $x$.
	The combination of \EqRef{Uniq55} and \EqRef{Uniq56} thus yield
\begin{equation}
	\bfp_* \bigl(\bfq_{Q_1*} (X^Q_i) \bigr)  =  - R^r_{0,i }\, \partial_{\theta^r_{X, M}}
\quad\text{and} \quad 
	\bfp_* \bigl(\bfq_{Q_2*} (Y^Q_i) \bigr)   =   S^r_{0,i } \,\partial_{\theta^r_{Y, M}},
\EqTag{Uniq66}
\end{equation}
	where  $R^r_{0,i} = R^r_i \circ  \bfq_{Q_1}$ and    $S^r_{0,i} = S^r_i \circ  \bfq_{Q_2} $.

	Since $dR^r_i \in \hZ^\infty$ (see equation \EqRef{dRdS}) , the functions $R^r_i$ are constants on every integral manifold of $\hZ^\infty$.
	But the map $\bfq_{Q_1} \: Q_1 \to N$ is an integral manifold of $\hZ^\infty$ (see Theorem \StRef{Dreduce}{\bf[v]})
	and therefore  the $R^r_{0,i}$ are constants.
	Similarly,   $dS^r_i \in \cZ^\infty$, the $S^r_i$ are constants on every integral manifold of $\cZ^\infty$ and therefore  the $S^r_{0,i}$ are constants.
        Consequently we can re-write  \EqRef{Uniq66}  (using the last 2 equations from  \EqRef{Uniq55}) as
\begin{equation}
	\bfp_* \bigl(\bfq_{Q_1*} (X^Q_i) \bigr)  =  \bfq_{P_1*}( R^r_{0,i } X^P_r )
\quad\text{and} \quad 
	\bfp_* \bigl(\bfq_{Q_2*} (Y^Q_i) \bigr)  =   \bfq_{P_2*}( S^r_{0,i } Y^P_r)
\EqTag{Uniq67}
\end{equation}	
	and therefore, by the commutative of the diagrams \EqRef{psi1psi2}, 
\begin{equation}
	\psi_{1,*}(X^Q_i) = R^r_{0,i} X^P_r 
	\quad\text{and}\quad 
	\psi_{2,*}(Y^Q_i) = S^r_{0,i} Y^P_r.
\end{equation}
	
	To complete the proof of \EqRef{psihomo} it suffices to show that $S^r_{0,i}=R^r_{0,i}$. 
	Equation \EqRef{matchTheta} states
\begin{equation}
\bftheta_{N,X}(q) = \bftheta_{N,Y}(q) ,\quad \text{and} \quad \bftheta_{M,X}(p) = \bftheta_{M,Y}(p) ,
\EqTag{SameT}
\end{equation}
	where $q=\bfq_{H_{\diag}}(q_1,q_2)$ and $p=\bfq_{G_\diag}(p_1,p_2)$, and therefore $\bfp^*\bftheta_{M,X}(q)=\bfp^*\bftheta_{M,Y}(q)$.
	Using \EqRef{SameT} in the definitions of $\bfR$ and $\bfS$ from equations \EqRef{Pullback4Adapted} (with $\phi = \bfp$)
	shows $\bfR(q)=\bfS(q)$, and hence  $S^r_{0,i}=R^r_{0,i}$. 

	Equations \EqRef{psihomo} are valid on any local trivialization of $Q_1$ and $Q_2$ and therefore hold globally. 
	Theorem \StRef{HGhomo} from  Appendix D can now be applied  to equations in \EqRef{psihomo}  to conclude that there exists 
	a common homomorphism $\phi \:H \to G$ such that the $\psi_a$ are  $\phi$ equivariant.


	With the assumption that $H' \subset G$ is closed, it follows (see  Remark \StRef{NReduction}) 
	that the action of $H'_\diag$ is regular on $P_1 \times P_2$. 
	The equivariance of the maps $\psi_a$ with respect to the actions of $H$ and $H'$  
	then induces a map $\Psi:N \to (P_1\times P_2)/H'$ such that the diagram
\begin{equation}
\begin{gathered}
\begindc{\commdiag}[3]
\obj(0, 0)[L]{$Q_1 \times Q_2$}
\obj(40, 0)[P]{$P_1 \times P_2$}
\obj(0, -17)[N]{$N$}
\obj(40, -17)[M]{$(P_1\times P_2)/H'$}
\mor{L}{N}{$\bfq_{H}$}[\atright, \solidarrow]
\mor{P}{M}{$\bfq_{H'}$}[\atleft, \solidarrow]
\mor{N}{M}{$\Psi$}[\atleft, \solidarrow]
\mor{L}{P}{$\psi_1\times \psi_2$}[\atleft, \solidarrow]
\enddc
\end{gathered} 
\end{equation}
	commutes. Since $\psi_1\times \psi_2$ is a local diffeomorphism, the function $\Psi$ is as well.  Since $\psi_a^*(\CalK_a)=\CalL_a$, 
	a  final  application of  Theorem 2.2 gives
$$
	\Psi^*\bigl( (\CalK_1+\CalK_2)/H'\bigr) = (\CalL_1+\CalL_2)/H = \CalE .
$$
	and part {\bf[iv]} is established.
\end{proof}	
 
	The next corollary establishes the uniqueness of  the quotient representation of a Darboux integrable system.
\begin{Corollary}  Let $(\CalL_a,Q_a, H)$ and $(\CalK_a,P_a, G)$ 
	satisfy the hypotheses of Theorem \StRef{Dreduce} and suppose, in addition that $P_a$, $Q_a$, $H$ and $G$ are all connected. If  
\begin{equation}
	(Q_1 \times Q_2, \CalL_1+\CalL_2)/H_{\diag} \cong (P_1 \times P_2, \CalK_1+\CalK_2)/G_{\diag},
\end{equation}
	then the manifolds $P_a$ and $Q_a$ are diffeomorphic, the Lie groups $H$ and $G$ are isomorphic and $\CalL_a \cong \CalK_a$. 
\end{Corollary}
\begin{proof} We apply the Theorem \StRef{UIE} using $\bfp = I_M$, 
	the identity map on $M$, to construct smooth maps $\psi_a:Q_a \to P_a$, 
	and $\tilde \psi_a:P_a \to Q_a$ which are inverses to each other. 
	Similarly we have Lie group homomorphisms $\phi:H \to G$ and $\tilde \phi :G \to H$ which are inverses.
\end{proof}
\section{B\"acklund Transformations for Darboux  Integrable Systems}

In this last section we apply our results from the previous sections to produce B\"acklund transformations for Darboux integrable systems by combining Theorem A from the introduction  with Theorem \StRef{Dreduce}. This will provide a proof of Theorem E from the introduction.  We also apply Theorem \StRef{UIE}  which gives sufficient conditions for when a B\"acklund transformation for Darboux integrable systems arises using group quotients.

We begin with the following general  method of
	constructing  B\"acklund transformations for a given Darboux integrable system $\CalI$.

\begin{Theorem} 
\StTag{BigBack}
	Let $(\CalI, M)$ be a Darboux integrable system  with canonical quotient representation 
	$(\CalK_1 + \CalK_2, M_1 \times M_2)/G_{\diag}$. 
	Let $L \subset G\times G$ be a subgroup, let $H_\diag = L \cap G_\diag$ and suppose that the 
	actions of $L$ and $H_\diag$ are regular on $M_1 \times M_2$. Then the 
	commutative diagram of differential systems
\begin{equation}
\begin{gathered}
\begindc{\commdiag}[3]
\obj(0, 34)[I]{$(\CalK_1 +  \CalK_2)$}
\obj(0, 10)[H]{$(\CalK_1 + \CalK_2) /H_\diag$}
\obj(-36, -6)[I1]{$(\CalK_1 + \CalK_2)/L$}
\obj(36,  -6)[I2]{$(\CalK_1 +  \CalK_2) /G_\diag$}
\mor{I}{H}{$\bfq_{H_\diag}$}[\atleft, \solidarrow]
\mor(1, 8)(-35, -4){$\bfp_1$}[\atleft, \solidarrow]
\mor(-1, 8)(35, -4){$\bfp_2$}[\atright, \solidarrow]
\mor(-1, 34)(-38, -5){$\bfq_{L}$}[\atright, \solidarrow]
\mor(1, 34)(38, -5){$\bfq_{G_{\diag}}$}[\atleft, \solidarrow]
\enddc
\end{gathered}
\EqTag{BackDiag9}
\end{equation}
	defines a  B\"acklund transformation $\CalE = (\CalK_1 + \CalK_2)/H_\diag$  
	between $\CalI = (\CalK_1 +  \CalK_2)/G_\diag$ and the Darboux integrable system $\CalJ = (\CalK_1 + \CalK_2)/L$. 
\end{Theorem}

	Theorem \StRef{BigBack} is a restatement of the first portion of Theorem E in the introduction if we identify
	$G_2=G_\diag$, $\CalI_2= (\CalK_1 +  \CalK_2) /G_\diag$, $G_1= L$, and $\CalI_1=(\CalK_1 + \CalK_2)/L$.  
	The B\"acklund transformation is $\CalB =(\CalK_1 + \CalK_2) /H_\diag$. Note that Theorem \StRef{BigBack} 
	constructs a B\"acklund transformation with fibre dimensions $\dim L - \dim H$ for $\bfp_1$ and $\dim G - \dim H$ for $\bfp_2$. 

	Theorem E states that the dimension of the space of Darboux invariants for $\CalI_1$ is larger than that of $\CalI_2$.  We  now
	show why  this is so. 	With  $\hW$ and $\cW$ given by \EqRef{Wsing}, 
Theorem \StRef{Dreduce}{\bf[ii]}, gives the singular Pfaffian systems for 
	$\CalI$ and $\CalJ=(\CalK_1 + \CalK_2)/L$ as
\begin{equation}
	\hV  =\hW/G_\diag, \quad   \cV = \cW/G_\diag, \quad      \hZ =\hW/L, \quad   \cZ = \cW/L.  
\EqTag{B2}
\end{equation}
The next corollary  shows that when  $L \neq G_\diag$, the system $(\CalK_1 + \CalK_2)/L$ has {\it strictly more}  Darboux invariants than
	$\CalI$.
	
\begin{Corollary} \StTag{s9c1}
	The  integrable subsystems  of  the singular Pfaffian systems of   $\CalI$ and $\CalJ$  in \EqRef{B2} satisfy 
\begin{equation}
\rank (\hZ^\infty) \geq \rank (\hV^\infty) \quad\text{and} \quad \rank (\cZ^\infty) \geq \rank (\cV^\infty).
\EqTag{Invcom}
\end{equation}
	Moreover, if $G$ is connected then we have equalities in \EqRef{Invcom} if and only if $\CalJ =\CalI$.
\end{Corollary}
\begin{proof}	
	
	Equations \EqRef{WMODG4},  when  applied  successively  first to 
	the diagonal subgroup $G_{\text{diag}}$  and then to the group $L$ 
        gives
\begin{alignat*}{2}
	\rank (\hV^\infty) &= \dim M_2 - \dim G,  &\quad
	\rank (\cV^\infty) &= \dim M_1- \dim G \quad \text{and}
\\
	\rank (\hZ^\infty) &= \dim M_2 - \dim L_2, &\quad \rank (\cZ^\infty) &= \dim M_1- \dim L_1.
\end{alignat*}
	The combination of these equations then yields
\begin{equation}
	\rank (\hZ^\infty)  = \rank (\hV^\infty)+ \dim G - \dim L_1
	\text{ and }
	\rank (\cZ^\infty)  = \rank (\cV^\infty)+ \dim G - \dim L_2.
\EqTag{Vr3}
\end{equation}
	Since $L_a = \rho_a(L) \subset G$ (see \EqRef{Laction}),  we have $\dim L_a \leq \dim G$.
	This fact, along with \EqRef{Vr3}, proves  \EqRef{Invcom}. 
	Equality holds in  \EqRef{Invcom} only when $\dim L_1= \dim L_2= \dim G$. If $G$ is connected  this implies 
	$ L_1= L_2= G$ and $L$ is the diagonal subgroup. 
\end{proof}

Corollary  \StRef{s9c1} combined with Theorem  \StRef{BigBack} gives a  complete proof of Theorem E in the introduction.

As our last result we use Theorem \StRef{UIE} in Section 8 to give sufficient conditions for when a B\"acklund transformation $\CalB$ between Darboux integrable system can be constructed using quotients as in Theorem \StRef{BigBack}.

\begin{Corollary} \StTag{c9c2} Let 
\begin{equation}
\begin{gathered}
\begindc{\commdiag}[3]
\obj(0, 10)[B]{$\CalB$}
\obj(-13, 0)[I]{$\CalJ$}
\obj(13, 0)[J]{$\CalI$\, .}
\mor{B}{I}{$\bfp_1$}[\atright, \solidarrow]
\mor{B}{J}{$\bfp_2$}[\atleft, \solidarrow]
\enddc
\end{gathered}
\EqTag{BigBackd}
\end{equation}
	define a B\"acklund transformation $\CalB$ between differential systems  $\CalI$ and $\CalJ$.  
	Suppose $\CalI$ is Darboux integrable and that the pair $(\CalB, \CalI)$ is maximally compatible.  
	If  $\CalI  =   (\CalK_1 + \CalK_2)/G_\diag$ then there is subgroup $H \subset G$ such  that
	$\CalB \cong_{\text{\rm loc}}  (\CalK_1 + \CalK_2)/H_\diag$ in which case  the right half of the diagram  \EqRef{BigBackd} 
	coincides with the right  half of \EqRef{BackDiag9}.
\end{Corollary}
\begin{proof}
	This is just a restatement of  Theorem \StRef{UIE}.
\end{proof}	

Corollary \StRef{c9c2} does not imply that $\CalJ$ can be constructed from group quotient  methods, but only $\CalB$. However, in the special cases of a B\"acklund transformations between a Darboux integrable hyperbolic Monge-Amp\`ere PDE in the plane represented by $\CalI$ and the wave equation being represented by $\CalJ$ in equation \EqRef{BigBackd}, the entire B\"acklund transformation  can be constructed using Theorem \StRef{BigBack}, see \cite{anderson-fels:2012a}.

\begin{Remark}
\StTag{ClassicalDI}
Consider  the extreme case where $L = G \times G$ in Theorem \StRef{BigBack}. In this case one has 
$$
	 (\CalK_1+\CalK_2)/L =\CalK_1/G+\CalK_2/G,
$$
$H_\diag = G_\diag$, the projection map $\bfp_2$ is the identity 
	and  the Vessiot group for  $(\CalK_1 + \CalK_2)/L$ is the identity group. 
	As an integrable extension, the integral manifolds of  $\CalI =  (\CalK_1 + \CalK_2)/G_\diag$	
	are obtained directly from integral manifolds  of  $(\CalK_1+\CalK_2)/L$ by solving ODE. 
	{\it This simple remark  effectively  describes the entire classical integration method of  Darboux integrable equations.}
	See also \cite{anderson-fels:2013a}.
\end{Remark}

The starting point for building B\"acklund transformation using Theorem \StRef{BigBack} is the canonical quotient representation of $\CalI$ in Theorem  \StRef{Dreduce}. If however only the local canonical quotient representation of $\CalI$ is given, as assured by Theorem  \StRef{Dreduce}, we can replace the action of $L$ in Theorem \StRef{BigBack} by a local action of $L$ and produce local B\"acklund transformations.  Examples can be found in \cite{anderson-fels:2012a}. Similarly if the B\"acklund transformation\EqTag{BigBack} in Corollary \StRef{c9c2} is only given locally then $\CalB$ can be constructed by taking a subgroup $H$ of the Vessiot group $G$ of $\CalI$ where $H$  acts  only locally.

\appendix

\section{On the Prolongation of Integrable Extensions}
	Let $\bfp:(\CalE, N) \to (\CalI, M)$ be an integrable extension with admissible subbundle  $J\subset T^*N$. 
	Let $\pi_N: G_k(TN) \to N$ and  $\pi_M : G_k(TM) \to M$ be the Grassmann bundles of $k$-planes and  let 
	$G_k(\CalE) \subset  G_k(TN)$ and $G_k(\CalI) \subset G_k(TM)$ be the spaces of $k$-dimensional integral elements for 
	$\CalE$ and $\CalI$ respectively.
	We assume that $\pi_N:G_k(\CalE)\to N$ and $\pi_M: G_k(\CalI)\to M$ are smooth bundles and that the 
	inclusion maps $\iota_N : G_k(\CalE) \to G_k(TN)$ and $\iota_M :G_k(\CalI) \to G_k(TM)$  are smooth immersions.
	The  prolongation spaces for $\CalE$ and $\CalI$ are $N^{[1]} = G_k(\CalE)$ and $M^{[1]} =G_k(\CalI)$ 
	and the prolonged differential systems are  
	$\CalE^{[1]}= \iota_N^* (\CalC_N)$ and $\CalI^{[1]}= \iota_M^*( \CalC_M)$, 	
	where $\CalC_N$ and $\CalC_M$ are the canonical Pfaffian systems on $G_k(TN)$ and $G_k(TM)$ respectively. 
	Our goal is to prove that  $(\CalE^{[1]}, N^{[1]})$ is an integrable extension of $(\CalI^{[1]}, M^{[1]})$.
	This fact is listed in Section 2 as  \IntExtProl. In order to state this result  in a more precise manner  we first need the 
	following preliminary lemma.

\begin{Lemma}\StTag{isoie}  Let $ \bfp :(\CalE, N) \to (\CalI, M)$ be an integrable extension with admissible subbundle  $J\subset T^*N$.
	Then the map $\bfp_*:TN\to TM$ defines a bijection between $k$-dimensional integral elements 
	of $\CalE$ at $x\in N$ and $k$-dimensional integral
	elements of $\CalI$ at $\bfp(x)$.
\end{Lemma}
\begin{proof} 
	By conditions {\bf [i]} and {\bf [ii]} in  \EqRef{IntExtDef},  we have $TN = \ker (\bfp_*)\oplus \ann(J)$. 
	Therefore $\bfp_* : \ann(J) \to TM$ is a bundle map which is an isomorphism on each fibre, that is, $\ann(J)$ is horizontal.
	We shall use this simple fact repeatedly.

	Now let  $x\in N$ and let $P_x\subset T_xN$ be a $k$-dimensional integral element of $\CalE$ at $x$. 
	If $\theta \in \CalI_x$ then $\bfp^* \theta \in \CalE$ and therefore $\theta|_{ \bfp_*(P_x)} = (\bfp^*\theta)|_{P_x}=0$.	
	Therefore $\bfp_*(P_x)$ is an integral element of $\CalI$ at $\bfp(x)$. 
	Since $P_x$ is an integral element of $\CalE$ and $\CalS(J)  \subset \CalE$, we also have $P_x \subset \ann(J_x)$.
	Since the restriction of  $\bfp_* $ to  $\ann(J_x)$ is an isomorphism, the dimensions of $P_x$ and $\bfp_*(P_x)$ are the same.

	Let $P_{1,x}$ and $P_{2, x}$ are two  
	$k$-dimensional integral elements of $\CalE$ at $x$ and suppose $\bfp_*(P_{1,x}) = \bfp_*(P_{2, x})$. Then
	$P_{1,x} \subset \ann(J_x)$ and $P_{2,x} \subset \ann(J_x)$ and again,  because the restriction of  $\bfp_* $ to  $\ann(J_x)$ 
	is an isomorphism,  $P_{1, x} = P_{2,x}$ and therefore $\bfp_*$ is a  one-to-one mapping on integral elements.

	To show that $\bfp_*$ is onto, let $P'_{\bfp(x)}\subset T_{\bfp(x)} M$
	be a $k$-dimensional integral element of $\CalI$ at $\bfp(x)$. Let
\begin{equation}
	P_x = \bfp_*^{-1}(P'_{\bfp(x)}) \cap \ann(J)\subset T_xN,
\EqTag{EE}
\end{equation}
	that is,  let  $P_x$  be  the subspace of $\ann(J_x)$ which maps by $\bfp_*$ to $P'_{\bfp(x)}$. The  foregoing arguments imply that 
	$P_x$ and $P'_x$ have the same dimension so that  we need  only to check that $P_x$ is an integral element of $\CalE$ 
	to finish the proof. Since $\CalE$ is generated algebraically 
	by $\CalS(J)$ and $\bfp^*(\CalI)$ we need to only check $\bfp^*\theta(P_x) =0 $ for all $\theta \in \CalI_{\bfp(x)}$ 
	and $\rho(P_x) = 0$ for all $\rho \in J_x$. These are both trivially true  and
	therefore $\bfp_*$ defines a bijection on $k$-dimensional integral elements.
\end{proof}

	Lemma \StRef{isoie} implies the following.

\begin{Corollary} \StTag{proNM} 
	The map  $ \bfp^{[1]} :N^{[1]}   \to M^{[1]}  $ defined by
\begin{equation}
	\bfp^{[1]} (P_x) = \bfp_* (P_x)
\EqTag{Dprophi}
\end{equation}
	is a  smooth submersion, gives rise to the commutative diagram
\begin{equation}
\begin{gathered}
\begindc{\commdiag}[3]
\obj(0, 0)[I]{$N^{[1]}$}
\obj(30, -15)[H]{$M$}
\obj(30, 0)[H2]{$M^{[1]}$}
\obj(0, -15)[I1]{$N$}
\mor{I}{H2}{$\bfp^{[1]} $}[\atleft, \solidarrow]
\mor{H2}{H}{$\pi_M$\ ,}[\atleft, \solidarrow]
\mor{I1}{H}{$\bfp$}[\atright, \solidarrow]
\mor{I}{I1}{$\pi_N$}[\atright, \solidarrow]
\enddc
\EqTag{proie}
\end{gathered}
\end{equation}
	and  restricts to a diffeomorphism between the fibres of  $\pi_N$  and $\pi_M$.
\end{Corollary}

With these preliminaries in hand, we can now state our theorem on the prolongation of integrable extensions (see \IntExtProl).
 
\begin{Theorem}  \StTag{proIE}  
	Let $ \bfp :(\CalE, N) \to (\CalI, M)$ be an integrable extension.  Then  the submersion $\bfp^{[1]}:N^{[1]} \to M^{[1]}$, given by 
	equation 	\EqRef{Dprophi},  defines  $\CalE^{[1]}$  as an integrable extension of $\CalI^{[1]}$. 
	Moreover, if $J\subset T^*N$ is an admissible bundle for the extension $\CalE$ of $\CalI$, then 
	$\pi_N^*(J)$ is an admissible bundle for the extension $\CalE^{[1]}$ of $\CalI^{[1]}$.
\end{Theorem}

	To prove Theorem \StRef{proIE},   we will check the three
	conditions of equation \EqRef{IntExtDef} for  $(\CalE^{[1]}, N^{[1]})$ to be an integrable extension of $(\CalI^{[1]}, M^{[1]})$
	by showing that  $\pi^*_N(J)$ is an admissible  subbundle. 
	To this end, we shall need the  following two lemmas.

\begin{Lemma} \StTag{isos} 
	The maps $\bfp^{[1]}_* : \ker (\pi_{N,*}) \to \ker (\pi_{M,*})$ and $\pi_{N,*} : \ker (\bfp^{[1]}_*) \to \ker (\bfp_*)$ 
	are isomorphisms, and
\begin{equation}
	\ker (\pi_{N,*}) \cap \ker (\bfp^{[1]}_*) = 0 .
\EqTag{ICC1}
\end{equation}
\end{Lemma}

\begin{proof}  We use Corollary \StRef{proNM} and the commutative diagram   \EqRef{proie}.
	The first claim follows because $\bfp^{[1]}$ is a diffeomorphism on each fibre. 
	To prove   \EqRef{ICC1}, let $ Y \in \ker \pi_{N,*} \cap \ker \bfp^{[1]}_* $. 
	Then  $Y $ is $\pi_N$ vertical and $\bfp^{[1]}_* Y =0$. But $\bfp^{[1]}_*$ is
	an isomorphism on vertical vectors, therefore $Y = 0$.

	To  prove that $\pi_{N,*} : \ker (\bfp^{[1]}_*) \to \ker (\bfp_*)$ is an isomorphism, we begin with the fact 
	that $\ker \bfp_* = \pi_{N,*} (\ker ( (\pi_M \circ\bfp^{[1]})_*)) $ and use this fact to show that $\pi_{N,*}$ is onto. 
	If $X \in \ker \bfp_*$,  then  there exists $ Y \in \ker ((\pi_M \circ\bfp^{[1]})_*) $ with $\pi_{N,*}(Y) = X$.  By definition 
	$\bfp_*^{[1]} Y $ is $\pi_M$ vertical. Since, by Corollary \StRef{proNM}, $\bfp^{[1]}$ is a diffeomorphism on the fibres,
	there exists $Z \in T_p N^{[1]}$ which is $\pi_N$ vertical and such that
	$\bfp_*^{[1]}(Z) = \bfp_*^{[1]}(Y)$. Then, with $Y' = Y-Z$, we have  that $\pi_{N,*}(Y') = \pi_{N,*}(Y)=X$
	(because $Z$ is  $\pi_{N, *}$ vertical), and $\bfp^{[1]}_*(Y') =\bfp^{[1]}_*(Y) - \bfp^{[1]}_*(Z) = 0$. 
	This proves that $\pi_{N,*} : \ker (\bfp^{[1]}_*) \to \ker (\bfp_*)$  is onto. The fact that $\pi_{N,*}$ is one-to-one follows easily from equation  \EqRef{ICC1}.
\end{proof}

\begin{Lemma} 
\StTag{DS} 
	The bundles $E^{[1]} $ and $I^{[1]}$ satisfy
\begin{equation}
	E^{[1]} = \pi_N^*(J) \oplus \bfp^{[1],*}(I^{[1]}).
	\EqTag{ESS}
\end{equation}
\end{Lemma}

\begin{proof} 
	Let $p =P_x\in N^{[1]}$ be an integral $k$-plane for $\CalE$ at $x\in N$.  To prove this lemma we first  note that 
	$\bfp^* (\ann (\bfp_* P_x)) \cap J_x  = 0 $ in which case the annihilator of \EqRef{EE} gives
\begin{equation}
	\ann(P_x) = \ann (\bfp_*^{-1}(\bfp_*(P_x)) + J_x = \bfp^* (\ann (\bfp_* P_x)) \oplus J_x .
\EqTag{annE}
\end{equation}
	By the definition of the canonical Pfaffian system,
\begin{equation}
	E^{[1]}_p = \pi_N^*( \ann(P_x))
\EqTag{E11}
\end{equation}
	while by equation \EqRef{Dprophi},
\begin{equation}
	I^{[1]}_{\bfp^{[1]}(p)} = \pi_M^*\left( \ann( \bfp_*(P_x)) \right) .
\EqTag{E12}
\end{equation}
	Applying   \EqRef{annE} and the commutative diagram  \EqRef{proie} to 
	equation \EqRef{E11}, and using equation \EqRef{E12}, we arrive at
$$
	E^{[1]}_p = \pi_N^*(J_x \oplus  \bfp^* (\ann( \bfp_* P_x)) ) 
	=  \pi_N^*(J_x )  \oplus  \bfp^{[1],*} \pi_M^* (\ann( \bfp_* P_x)) 
	= \pi_N^* (J_x) \oplus \bfp^{[1],*} I^{[1]}_{\bfp^{[1]}(p)},
$$
	which proves the lemma.
\end{proof}

\begin{proof}[Proof of  Theorem  \StRef{proIE}]

	To prove that $\bfp^{[1]}: (\CalE^{[1]}, N^{[1]}) \to  (\CalI, M^{[1]})$ is an integrable extension, we must check the three
	conditions of  equation \EqRef{IntExtDef} with  $\pi^*_N(J)$  in place of  $J$.
	By the remark at the beginning of the proof of  Lemma \StRef{isoie} 
	we have that  $\rank J = \dim N - \dim M$ and therefore, in view  of  Corollary
	\StRef{proNM},
\begin{equation}
	 \rank \pi_N^*(J)  = \rank J =  \dim N - \dim M =  \dim N^{[1]} -\dim M^{[1]}.
\end{equation}
	This proves condition {\bf [i]} in  \EqRef{IntExtDef}.  

	To prove the transversality condition {\bf [ii]} in  \EqRef{IntExtDef},  let $ Y \in \ann(\pi_N^* J_x) \cap \ker \bfp^{[1]}_*$.
	This implies $ \pi_{N,*}Y\in \ann(J_x)$, while Lemma  \StRef{isos} implies $\pi_{N,*} Y \in \ker \bfp_*$.
	By the transversality of the admissible  subbundle  $J$  for $\CalE$ (see {\bf [ii]}  in  equation  \EqRef{IntExtDef}) 
	these two conditions on $\pi_{N,*} Y$ show that  $\pi_{N,*} Y = 0$. The fact that $\pi_{N,*}: \ker \bfp^{[1]}_* \to \bfp_*$ is 
	an isomorphism (Lemma \StRef{isos}) implies $Y=0$. Thus
	$\ann(\pi_N^* J_x) \cap \ker \bfp^{[1]}_* = 0$ and condition {\bf [ii]}  is verified for the prolongation.

	It remains to prove condition {\bf [iii]} in   \EqRef{IntExtDef}, which in this case reads,
	\begin{equation}
	\CalE^{[1]} = \langle \bfp^{[1],*} \CalI^{[1]} + \CalS( \pi_N^* J) \rangle_{\text{\rm alg}}.
\EqTag{IE1p}
	\end{equation}
Lemma \StRef{DS} shows that \EqRef{IE1p} holds at the level of 1-forms. 
Since $\CalE^{[1]}$ is a Pfaffian system it remains to be shown that the right-hand side of \EqRef{IE1p} is a Pfaffian system. Since $\CalI^{[1]}$ is a Pfaffian system, it is then sufficient to verify condition \EqRef{IntExt4} for a local basis of $ \CalS( \pi_N^* J)$. 
		
	Let  $ \{\, \xi^u \,\}$  be a set of  1-forms on $N$
	which define a local basis for $\CalS(J)$. Then  the 1-forms $\sigma^u=  \pi_{N}^* \xi^u$  define  a local basis for $\CalS(\pi_{N}^*(J))$. 
	Since  $d \xi^u \equiv  0 \mod  J + \bfp^* (I^2)$,  we immediately have
\begin{align*}	
	d \sigma^u &  \equiv0  \mod \pi^*_N(J) + \pi^*_N( \bfp^*( I^2)) 
\\
	&\equiv0 \mod    \pi^*_N(J) +  \bfp^{[1],*} ( \pi^*_M(I^2)). 
\end{align*}
	We remark that it is not generally true that $\pi^*_M(I^\ell ) \subset I ^{[1], \ell}$. However, for $\ell = 2$
	one can readily check that  $\pi^*_M(I^2 ) \subset I ^{[1], 2}$ and therefore
\begin{equation}
	d \sigma^u = 0  \mod  \pi^*_N (J) + \bfp^{[1],*} (I^{[1],  2}).
\end{equation}
This shows that the right hand side of \EqRef{IE1p} is differentially closed. This establishes \EqRef{IE1p} and   hence {\bf [iii]} in   \EqRef{IntExtDef} is proved.
\end{proof}

\section{A Remark on Involutivity}

	In this Appendix we establish a simple property (see Corollary \StRef{I1D} in particular) of the prolongation of an involutive linear Pfaffian systems which we use in the 	proof of Theorem \StRef{ProDI}.

	Let $\CalI$ be an involutive linear Pfaffian system on a manifold $M$ with independence condition $\omega \in \Omega^n(M)$. 
	Then, about each point $x\in M$, there  is an open set $U$, a local basis of $1$-forms $\{ \theta^\alpha \} _{1\leq \alpha \leq m}$ for $\CalI$,
	and a local coframe  $\{ \theta^\alpha, \, \omega^i, \, \pi^\epsilon   \}_{1\leq i \leq n, 1 \leq \epsilon \leq p }$ for $M$ 
	such that $\omega = \omega^1 \wedge\cdots \wedge \omega^n$, and 
	(see page 128 in \cite{bryant-chern-gardner-griffiths-goldschmidt:1991a})
\begin{equation}
	d \theta^ \alpha = A^\alpha _{\epsilon i} \,\pi^\epsilon \wedge \omega^i \quad \mod I .
\EqTag{stre1}
\end{equation}
Let  $M^{[1]}$ be the space of integral $n$-planes of $\CalI$ and let $\CalI^{[1]}$ be the linear Pfaffian system on $M^{[1]}$ which is  the prolongation of $\CalI$. The prolongation  of $\CalI$ is determined locally as follows. Let  $S^\epsilon_{i, v}$be smooth functions on the open set $U$ from above which form a basis for the $t$-dimensional solution space of the linear system
\begin{equation}
	A^\alpha _{\epsilon i} \Sigma^\epsilon_{j}\, \omega^j\wedge  \omega^ i = 0
\EqTag{appb1}
\end{equation}
in the  unknowns $\Sigma^\epsilon_j$. Then on the open set $U^{[1]}=U\times \real^t\subset M^{[1]}$  define the $1$-forms
\begin{equation}
	\tilde \theta^ \epsilon = \pi^\epsilon - s^v S^\epsilon _{i,v} \omega^i,
\EqTag{pieps}
\end{equation}
where $ 1\leq v \leq t$ and $s^v$ are coordinates on $\real^t$. On $U^{[1]}$ the prolongation  of $\CalI$ is then given by
\begin{equation}
\CalI^{[1]}= \langle \theta^\alpha, \tilde \theta^\epsilon \rangle_{{\rm diff}}.
\EqTag{bpI11}
\end{equation}

	Let $\pi:M^{[1]}\to M$ be the projection map. Then it is clear from equation \EqRef{bpI11} that $\pi^*(\CalI^1) \subset \CalI^{[1]}$ (see also page 150 in \cite{bryant-chern-gardner-griffiths-goldschmidt:1991a}). A form $\eta \in \Lambda^*(M^{[1]})$ is $\pi$ semi-basic if $X\hook \eta=0$ for all $X \in \ker(\pi_*)$. Equations \EqRef{pieps} and \EqRef{bpI11} show that the  one-forms  in $I^{[1],1}$ are $\pi$ semi-basic.

The technical result  that we need in Theorem \StRef{ProDI} is now given  by the following. 

\begin{Theorem} 
\StTag{kerS} 
	Let $\CalI $ be an involutive linear Pfaffian system, and let $ S^\epsilon_{i, v}$ be the $t$-dimensional basis for the solution space of \EqRef{appb1}. If $k_\epsilon S^\epsilon_{i, v} = 0$,  then $k_\epsilon =0$. 		
\end{Theorem}
	
The proof of this theorem is based upon the following lemma.
	
\begin{Lemma}
 \StTag{NOpiB}  
	Let $\CalI$ be an involutive linear Pfaffian system. If $\eta \in \CalI^{[1],1}$ is a one-form 
	and $d \eta $ is $\pi$ semi-basic, then $\eta \in \CalS (\pi^*( I^1))$.
\end{Lemma}

\begin{proof}  First we note that this condition only needs to be checked point-wise. Therefore let $x\in M$ and let $U^{[1]}\subset M^{[1]}$ be an open set as chosen above. We now use the condition of  involutivity to 
	 assume that a coframe on $U$ is chosen so that the tableaux  $A^\alpha_{\epsilon i}$ has the form given by equation (90) of Chapter IV 
	in \cite{bryant-chern-gardner-griffiths-goldschmidt:1991a}. This means we replace the $1$-forms $\pi^\epsilon$ 
	by the so-called principal components  $\bar \pi^a_{i}$. Let $s'_k$ be the last non-zero Cartan character.
	Then by involutivity $s'_1+\ldots +s'_k = t  $, where $t$ is dimension of the kernel in equation \EqRef{appb1}, and we set
\begin{equation}
	\theta^a_i = \bar \pi^a_i - p^a_{ij} \,\omega ^j,\quad    \quad i \leq k, \ a \leq s_i' ,
\EqTag{FB}
\end{equation}
where $p_{ij}=p_{ji}$. 
	The number of independent functions $p^a_{ij}$ is  $t=s'_1+ 2 s'_2+ \ldots + k s'_k$
	and these define the fibre coordinates for the projection $\pi:U^{[1]}\to U$. A local basis for the prolongation of 
	$\CalI^{[1]}$ is then given by
\begin{equation}
	\CalI^{[1]}|_{U^{[1]}} = \langle \theta^\alpha , \, \theta^a _ i \rangle_{\text{\rm diff}}, \quad \quad i \leq k, \ a \leq s_i' ,
\EqTag{appBprI1}
\end{equation}
	and a local coframe on $U^{[1]}$ is
\begin{equation}
	\{\, \theta^\alpha, \, \theta^a_i, \omega^i, \,  dp^a_{ij} \, \} , \quad 1\leq i \leq j \leq k, \ a \leq s_i' .
\EqTag{CFU1}
\end{equation}

Suppose that $\eta$ is a one-form in $\CalI^{[1]}$ satisfying the conditions in the statement of the lemma. 
	Then using the local coframe as in equation \EqRef{CFU1}, there exists smooth functions $T^i_a$ and $R_\alpha$ on $U^{[1]}$ such that
\begin{equation}
	\eta= T^i_a \, \theta ^a_i + R_\alpha \, \theta^ \alpha.
\EqTag{alD1}
\end{equation}
	In order for $d \eta $ in equation \EqRef{alD1} to be $\pi$ semi-basic at $y\in \pi^{-1}(x)$, it is necessary and sufficient that
\begin{equation}
\partial_{p^b_{jk}} \hook(  dT^i_a\,  \wedge \theta^a_i + T^i_a \, d\theta^a_i + dR_\alpha \wedge  \theta^\alpha + R_\alpha\, d\theta^\alpha)|_y = 0 .
\EqTag{SBC1}
\end{equation}
	Since $\theta^a_i$, $\omega^i$, $\theta^\alpha $, $d\omega^i$, $d\theta^\alpha$ and $d \bar \pi^a_i$ are $\pi$ semi-basic we find using this fact along with  equation \EqRef{FB}, that  \EqRef{SBC1} simplifies to 
\begin{equation}
	\left( (\partial_{p^b_{jk}} \hook dT^i_a) \theta^a_i + T_b^j \omega^k +( \partial_{p^b_{jk}} \hook  dR_\alpha ) \theta^\alpha \right)|_y  = 0 .
\end{equation}
	However $\theta^a_i, \omega^i$ and $\theta^\alpha$ are point-wise linearly independent 
	and so this equation implies $T_b^j(y)=0$ and therefore from the expression for $\eta$ in equation \EqRef{alD1} we have $\eta|_y \in\pi^*( I_x )$. 
\end{proof}

\noindent
\begin{proof}[Proof of  Theorem \StRef{kerS}]
	Suppose that $k_\epsilon S^\epsilon_{i, v}(x) = 0 $ at some point $x\in U$ where $k_\epsilon \in \real$. Let
	$\eta \in \CalI^{[1],1}$ be the $1$-form given using equation \EqRef{bpI11} by
\begin{equation}
	\eta = k_\epsilon \tilde \theta^\epsilon =  k_\epsilon( \pi^\epsilon - s^v S^\epsilon _{i, v} \omega^i).
\EqTag{appBalpha}
\end{equation}
	The exterior derivative of $\eta$	evaluated at $y\in \pi^{-1}(x)$ is then
\begin{equation}
\begin{aligned}
	d\eta|_y & = 
	k_\epsilon( d\pi^\epsilon - d s^v  \wedge S^\epsilon _{i, v}\, \omega^i  - s^vdS^\epsilon _{ i, v} \wedge \omega^i - s^vS^\epsilon_{i, v} \,d \omega^i  )|_y  \\
		&=k_\epsilon( d \pi^\epsilon-  s ^vdS^\epsilon _{i, v} \wedge \omega^i)|_y.
\end{aligned}
\end{equation}
	In this computation we have, as is customary, identified $S^\epsilon _{i, v}$ with $\pi^*S^\epsilon _{i, v}$. 
	Therefore $d\eta|_y$ is $\pi$ semi-basic which
	by Lemma \StRef{NOpiB} implies $\eta|_y \in \pi^*(I_x)$. From equations \EqRef{appBalpha} and \EqRef{bpI11}, this is clearly possible only if  $k_\epsilon  = 0$. This proves Theorem \StRef{kerS}.
\end{proof}

The geometric content of Theorem \StRef{kerS} is contained in the following corollary.

\begin{Corollary} 
\StTag{I1D} 
	If $\CalI$ is an involutive linear Pfaffian system, then ${I^{ [1] }}' =\pi^*(I ^1 )$.
\end{Corollary}

\begin{proof} Again note that this a point-wise condition. Since $\pi^*(I)\subset {I^{[1]}}'$, we need only show that if $\eta \in {I^{[1]}}'$ and has the form
$\eta =k _\epsilon \tilde \theta^\epsilon $,  (using the forms in \EqRef{bpI11}) then $\eta=0$. We therefore compute $d \eta \ \mod I^{[1]}$ which is
\begin{equation*}
	d\eta \equiv  k_\epsilon( d \pi^\epsilon -  S^\epsilon_{i,v}  d s^v \wedge \omega^i ) \quad  \mod I^{ [1] }. 
\end{equation*}
Since the forms $d\pi^\epsilon$ are $\pi$ semi-basic,  in order for the right hand side  to be  zero at any point $y\in U^{[1]}$ it is necessary that   $k_\epsilon S^\epsilon_{i,v}    = 0$ at $y$. Therefore by Theorem \StRef{kerS}, $k_\epsilon=0$ and hence $\eta=0 $.
\end{proof}

\section{On the definition of Darboux integrability}
	In this appendix we prove Theorem \StRef{DPD}.
	The key step to proving Theorem \StRef{DPD} is to show that the decomposability of $\CalI$, 
	together with conditions {\bf[i] }of Definition \StRef{Intro2}, 
	implies that   $\hV^\infty \cap \cV^\infty $ is an integrable Pfaffian system.  Then, since we are assuming that  
	$(\hV \cap \cV)^\infty  = \{\,0 \}$,
	we deduce that  $\hat V^\infty \cap \check V^\infty  \subset  (\hV \cap \cV)^\infty = \{\, 0\,\}$.

	To prove that $\hV^\infty \cap \cV^\infty$ is an integrable Pfaffian system  we shall need a 
	generalization of the  $0$-adapted coframe  defined in \cite{anderson-fels-vassiliou:2009a} (page 1917) or by \EqRef{ZeroC}.   	
	This new coframe is defined locally in a neighborhood of any given point and is constructed as follows. 
	First choose  independent one-forms $\bftau = \{\,\tau^1, \tau^2, \dots ,\tau^\ell\,\}$ such that 
\begin{equation*}
	\hV^\infty \cap \cV^\infty= \spn \{ \bftau\}.
\end{equation*} 
	Extend these by vector-valued (independent)  1-forms $\bfheta$ and $\bfceta$ in such manner that
\begin{equation}
	\hV^\infty \cap \cV=\spn \{\bftau, \, \bfheta  \} \quad \text{and} \quad  \hV \cap \cV^\infty= \spn \{\ \bftau, \, \bfceta \ \}.
\EqTag{DI23}
\end{equation}
	Then, just as in  \cite{anderson-fels-vassiliou:2009a}, these forms may in turn be extended (by conditions \EqRef{Intro7}) to a 
	local coframe  of vector-valued 1-forms $\{\ \bftheta,\, \bfhsigma, \, \bfcsigma,  \bfheta, \, \bfceta, \, \bftau \, \}$ on $M$ such that
\begin{equation}
\begin{aligned}
	\hV =  \spn \{\ \bftheta,\,  \bfheta,\,  \bfceta,\,  \bftau,\,  \bfhsigma \ \}, & \quad \hV^\infty = \spn \{\ \bfheta,\,  \bftau,\, \bfhsigma \ \}, 
\\
	\cV=  \spn \{\ \bftheta,\,  \bfheta,\, \bfceta,\, \bftau,\,  \bfcsigma \ \}, & \quad \cV^\infty = \spn \{ \ \bfceta,\, \bftau,\, \bfcsigma  \ \} .
\end{aligned}
\EqTag{ZeroC2}
\end{equation}
	We will call such a coframe $0$-adapted.  The first step in the proof of Theorem \StRef{DPD} is given by the next lemma. 
	In what follows we will use the convention that bold face Roman letters such as $\bfa$, $\bfA$, $\bfalpha, \dots$ 
	denote array-valued functions and differential forms of the appropriate  rank and dimensions.
\begin{Lemma} 
	If $\{\ \bftheta,\, \bfhsigma, \, \bfcsigma, \bfheta, \, \bfceta, \, \bftau \, \}$ is a 0-adapted coframe, 
	then the forms $\bftau$  (which span  $\hV^\infty \cap \cV^\infty$) satisfy the structure equations
\begin{equation}
	d\bftau  = \bfalpha \wedge \bftau + \bfa_1\, \bfceta \wedge \bfheta +\bfa_2\, \bfcsigma \wedge \bfheta +  \bfa_3\, \bfceta  \wedge \bfhsigma.
\EqTag{dtau1}
\end{equation}
\end{Lemma}

\begin{proof} The conditions $\bftau \in \hV^\infty $ and $\bftau \in \cV^\infty$ imply, 
	by the Frobenius condition for integrability and equations \EqRef{ZeroC2}, 
	that there exists 1-forms $\bfalpha$, $\bfbeta$, $\bfgamma$, $\bfmu$, $\bfnu$, $\bfxi$ such that
\begin{equation}
	d\bftau  = \bfalpha \wedge \bftau + \bfbeta \wedge \bfheta + \bfgamma \wedge \bfhsigma 
	\quad\text{and}\quad
	d\bftau  = \bfmu \wedge \bftau + \bfnu \wedge \bfceta + \bfxi \wedge \bfcsigma .
\EqTag{dbt1}
\end{equation}
	Now the 1-forms $\bftau \in \CalI^1 = \spn \{\, \tilde \theta^e\}$ and so, by the decomposability condition \EqRef{DecStEq},
	there are  no $\bfhsigma\wedge \bfcsigma$ terms in either of these structure equations.   
	Hence, after absorbing the $\bftau$, $\bfheta$  terms in $\bfgamma$ and the  $\bftau$, $\bfceta$  terms in $\bfxi$ into the other coefficients,
	we may re-write equations \EqRef{dbt1} in expanded form as
\begin{equation}
\begin{aligned}
	d\bftau 
&	= \bfalpha \wedge \bftau + \bfbeta \wedge \bfheta + (\bfc_1 \bftheta  + \bfc_2 \bfceta +\bfc_3 \bfhsigma ) \wedge \bfhsigma \quad \text{and} 
\\	
	d\bftau
&	= \bfmu \wedge \bftau + \bfnu \wedge \bfceta + (\bfC_1 \bftheta +\bfC_2 \bfheta  +\bfC_3 \bfcsigma ) \wedge \bfcsigma,
\end{aligned}
\EqTag{dbt2}
\end{equation}
	where the $\bfc_i$ and $\bfC_i$ are $3$-dimensional arrays of  locally defined smooth functions.
	For the right-hand sides of these equations to be equal we find immediately  that $\bfc_1=\bfc_3=\bfC_1=\bfC_3=0$.  
	Upon absorbing the   $\bftau$  terms in  $\bfbeta$ and $\bfnu$ into  $\bfalpha$ and  $\bfmu$,  
	we can re-write equations  \EqRef{dbt2}   in expanded form as 
\begin{equation}
\begin{aligned}
	d\bftau 
& 	= \bfalpha \wedge \bftau + (\bfb_1 \bftheta+\bfb_2 \bfheta +\bfb_3\bfceta+\bfb_4 \bfcsigma) \wedge \bfheta +  \bfc_2 \bfceta  \wedge \bfhsigma, \quad \text{and} 
\\
	d\bftau 
&	 = \bfmu \wedge \bftau +  (\bfB_1 \bftheta+\bfB_2 \bfheta +\bfB_3\bfceta+\bfB_4 \bfhsigma) \wedge \bfceta + \bfC_2\bfheta  \wedge \bfcsigma .
\end{aligned}
\EqTag{dbt3}
\end{equation}
	For  the right-hand sides of these equations  to be equal we further find that  $ \bfb_1= \bfb_2=\bfB_1=\bfB_3=0$ in which case equation \EqRef{dbt3}  reduces to \EqRef{dtau1},  as required.
\end{proof}

	To complete the proof  of Theorem \StRef{DPD}  we must show that the 
	coefficients  $\bfa_1$,  $\bfa_2$ and $\bfa_3$ in \EqRef{dtau1} vanish and, to this end, we shall need a refined coframe,  
	generalizing the 1-adapted coframe of  \cite{anderson-fels-vassiliou:2009a}.

\begin{Lemma}\StTag{I3}  
	Let $\{\ \bftheta,\, \bfhsigma, \, \bfcsigma,\,  \bfheta, \, \bfceta, \, \bftau \, \}$ be a 0-adapted coframe.
	There exists  vector-valued functions $\bfI_1,\bfI_2,\bfI_3$, and 
	$\bfJ_1,\bfJ_2,\bfJ_3$  and matrix-valued functions  $\bfE$ and $\bfF$
	such that the forms
\begin{equation}
	\bfhsigma_0 = d \bfI_3, \quad \bfheta_0 = d\bfI_2 + \bfE \bfhsigma_0 ,  \quad  \bfcsigma_0 = d \bfJ_3 \quad  \bfceta_0  = d \bfJ_2+\bfF \bfcsigma_0
\EqTag{cf0}
\end{equation}
	may be  used to define a  0-adapted coframe 
	$\{\, \bftheta, \, \bfheta_0, \, \bfceta_0, \, \bfhsigma_0,\, \bfcsigma_0, \, \bftau \, \}$.  Moreover, the 1-forms $\bftau$ can be expressed as 
\begin{equation}
	\bftau = \bfR\,( d\bfI_1 + \bfS \bfheta_0+ \bfT \bfhsigma_0) = \bfA\, ( d\bfJ_1 + \bfB \bfceta_0+ \bfC \bfcsigma_0),
\EqTag{taudef}
\end{equation}
	where $\bfR, \bfA, \bfS, \dots$ are smooth matrix-valued functions. The matrices   $\bfR, \bfA$
	are  invertible matrices of dimension $\rank (\hV^\infty \cap \cV^\infty)$.
\end{Lemma}

\begin{proof}
	We shall use the following simple observation to choose the functions $\bfI$ and $\bfJ$ 
	-- if  $K = \{\ \omega^1,\, \omega^2,\,\dots, \,\omega^k\ \}$ is any completely integrable Pfaffian (where the $\omega^i$ are independent 1-forms),
	then there are  (locally defined)  functions whose differentials will complete any given subset of the generators $\{\omega^{i_1},\dots, \omega^{i_p}\}$ to a basis of $K$.
 
	Accordingly, choose independent functions $\bfI_3$ such that $\hV^\infty = \spn \{ \bftau,\bfheta, d\bfI_3 \}$ 
	and let $\bfhsigma_0 = d\bfI_3$.  The first equation in \EqRef{cf0} therefore holds.
	Then choose independent functions $\bfI_2$  such that $\hV^\infty = \spn \{ \bftau, d \bfI_2, \bfhsigma_0 \} $. 
	The 1-forms $\bfheta$  can therefore be written as 
\begin{equation*}
	\bfheta = \bfP d\bfI_2 + \bfG \bftau+ \bfH \bfhsigma_0,
\end{equation*}
	where $\bfP$ is an invertible matrix of functions. Let
\begin{equation*}
	\bfheta_0 = \bfP^{-1}(\bfheta - \bfG \bftau) = d\bfI_2 +  \bfE \bfhsigma	
\end{equation*}
	so that the second equation in \EqRef{cf0} holds. Note, by \EqRef{DI23},  that  
	$\bfheta_0 \in \hV^\infty \cap \cV$ and $\hV^\infty = \spn\{ \, \bftau, \bfheta_0, \bfhsigma_0 \, \}$. 
	Finally,  if we choose $\bfI_1$ such that $\hV^\infty=\spn \{ d \bfI_1, \bfheta_0, \bfhsigma_0 \}$ 
	then the first equation in
	\EqRef{taudef} holds.  Similar arguments allow us to choose  $\bfJ_1,\bfJ_2,\bfJ_3$ so that the remaining equations \EqRef{cf0}  and \EqRef{taudef} hold.
\end{proof}

\begin{Remark} 
	The number of invariants $\bfI_1$ and $\bfJ_1$ are the same and equals the rank of  the bundle $\hV^\infty\cap \cV^\infty$.
\end{Remark}

\begin{Lemma} 
	For the coframe $\{\,\bftheta,  \bfhsigma_0,\bfcsigma_0, \bfheta_0,\bfceta_0, \bftau \, \}$ we have structure equations
\begin{equation}
	d \bfheta_0 = \bfa\,  \bftau \wedge \bfhsigma_0 + \bfb\, \bfheta_0 \wedge \bfhsigma_0 + \bfc\, \bfhsigma_0 \wedge \bfhsigma_0 .
\EqTag{dheta0}
\end{equation}
\end{Lemma}
\begin{proof} 
	Equation \EqRef{cf0} in Lemma \StRef{I3} gives
\begin{equation*}
	d \bfheta_0 = d\bfE \wedge \bfhsigma_0 .
\EqTag{deta}
\end{equation*}
	But, by definition, the 1-forms $\bfheta_0$ belong to $\CalI^1$ and therefore, by the decomposability condition \EqRef{DecStEq},  
	the 1-forms  $d \bfE$  contain no $ \bfcsigma_0$ terms and hence  $d\bfE \in \CalS(\hV)$.  Therefore $d\bfE \in \CalS(\hV^\infty)$  
	(see Lemma \StRef{dfL})
	and since   $\CalS(\hV^\infty) =  \spn  \{ \, \bfheta_0,\,  \bfhsigma_0,\, \bftau \, \}$
	 we have $d \bfE = \bfa\bftau + \bfb \bfheta_0 + \bfc \bfhsigma_0$ and  \EqRef{dheta0} holds.
\end{proof}

\begin{proof}[Proof of Theorem \StRef{DPD}] 
	Equation \EqRef{dtau1} holds for any 0-adapted coframe, in particular it holds for the 
	0-adapted coframe $\{\, \bftheta, \bfhsigma_0, \bfcsigma_0, \bfheta_0,\bfceta_0, \bftau\, \}$ constructed in Lemma \StRef{I3}. We now compare equation \EqRef{dtau1} with
	the result of  taking the exterior derivative of equation \EqRef{taudef} and utilizing \EqRef{cf0} 
	and \EqRef{dheta0},
\begin{equation}
\begin{aligned}
	d \bftau 
&	 = (d\bfR)\bfR^{-1} \wedge \bftau + \bfR(\  d\bfS \wedge \bfheta_0+\bfS d\bfheta_0 + d \bfT  \wedge\bfhsigma_0\ )
\\
& = 	(d\bfR)\bfR^{-1} \wedge \bftau + \bfR(\  d\bfS \wedge \bfheta_0+\bfS (\bfa \bftau \wedge \bfhsigma_0 + \bfb \bfheta_0 \wedge \bfhsigma_0 + \bfc \bfhsigma_0 \wedge \bfhsigma_0) + d \bfT \wedge \bfhsigma_0\ ).
\end{aligned}
\EqTag{dtau2}
\end{equation}
	Now decomposability implies $d \bfT$ is independent of the forms $\bfcsigma_0$. 
	Therefore   $d \bfT \in \CalS(\hV^\infty)$  $=  \spn  \{ \, \bfheta_0,\,  \bfhsigma_0,\, \bftau \, \}$. 
	Hence, from \EqRef{dtau2}, we deduce that  there are no terms of the form $\bfceta_0 \wedge 	\bfhsigma_0$ in  $d \bftau$.
	Therefore $\bfa_3=0$ in equation  \EqRef{dtau1} (in the current 0-adapted coframe).  
	A similar argument using the expression for $\bftau $ in terms of the invariants $\bfJ$ 
	gives $\bfa_2=0 $ in equation \EqRef{dtau1}, leaving
\begin{equation}
	d\bftau  = \bfalpha \wedge \bftau  +\bfa_1 \bfceta_0 \wedge \bfheta_0 .
\EqTag{df4}
\end{equation}
	Returning  to equation \EqRef{dtau2},  
	we  now conclude, from the absence of  the terms $\bfheta_0 \wedge \bfcsigma_0$ 
	in \EqRef{df4}, that $d\bfS$ is free of $\bfcsigma_0$ terms and so, 
	again by Lemma \StRef{dfL},  $d\bfS \in \CalS(\hV^\infty)$.
	Consequently equation \EqRef{dtau2}  implies that $d \bftau$ contains no terms $\bfheta_0\wedge \bfceta_0$ and  therefore $\bfa_1=0$ in equation \EqRef{df4}. 
	This proves that $\hV^\infty \cap \cV^\infty $ is  an integrable Pfaffian system  which, by the remarks made at the beginning of this
	appendix,  shows  that $\hV^\infty \cap \cV^\infty = \{\,0\,\}$.
\end{proof}

\section{Infinitesimal Equivariant Mappings}
	In this appendix  we prove the following theorem  which is required for the proof of Theorem \StRef{UIE}.

\begin{Theorem}  
\StTag{HGhomo}
	Let $H$ be a connected Lie group acting freely and regularly on $N$.  Let 
	$G$ be a Lie group acting freely and regularly on $M$ and  let   $\Gamma_H$ and  $\Gamma_G$ 
	be the corresponding Lie algebras of infinitesimal generators. Suppose $\Phi:N \to M$ is a smooth map
	whose differential   satisfies
\begin{equation}
\Phi_* : \Gamma_H \to \Gamma_G  
\EqTag{Phimono}
\end{equation}
	and is a monomorphism of the Lie algebras. Then there exists a  unique 
	homomorphism $\phi:H \to G$ such that $\Phi$ is $H$ equivariant,  that is, $\Phi(h \cdot p) = \phi(h) \cdot  \Phi(p)$
	for all $p\in N$ and $h \in G$. Moreover $\phi_*:\lieh \to \lieg$ is a Lie algebra monomorphism.
\end{Theorem}

\begin{proof} 
	Write $\mu: H \times N \to N$ and $\tau: G  \times M \to M$  for the actions of $H$ on $N$ and $G$ on $M$
	and fix a point $p\in N$. Then,  by the freeness and  regularity hypothesis for the actions of $H$ and $G$, 
	the maps $\mu_p: H \to N$ and $\tau_{\Phi(p)}:G \to  M$ given by (see \EqRef{actionmaps})
\begin{equation}
	\mu_p(h) = h\cdot p  \quad \text{and} \quad \tau_{\Phi(p)}(g) = g \cdot\Phi(p) 
\end{equation}
	are imbeddings. Equation \EqRef{Phimono}  implies $ \Phi\circ \mu_p :H \to M$ satisfies
\begin{equation}
(\Phi\circ \mu_p)_* TH \subset \bfGamma_G.
\EqTag{pfTH}
\end{equation}
	Since $H$ is connected, equation \EqRef{pfTH} shows that $\Phi(\mu_p(H))$ is contained
	in the maximal integral manifold of the  integrable distribution  $\bfGamma_G$ through $\Phi(p)$. 
	But the maximal integral manifold 	 of $\Gamma_G$ through $\Phi(p)$ is just the connected component  of the orbit of 
	$G$ through $\Phi(p)$ and therefore
\begin{equation}
	\Phi(\mu(H)) \subset G\cdot \Phi(p).
\EqTag{equiPhi}
\end{equation}
	Equation \EqRef{equiPhi}, together  with the  fact that the action of $G$ is free on $M$,
	implies there exists a unique function $\phi: H \to G$ such that the following diagram commutes
\begin{equation}
\begin{gathered}
\begindc{\commdiag}[3]
\obj(0, 0)[H]{$H$}
\obj(30, 0)[G]{$G$}
\obj(0, -15)[N]{$N$}
\obj(30, -15)[M]{$M \ .$}
\mor{H}{G}{$\phi$}[\atleft, \solidarrow]
\mor{H}{N}{$\mu_p$}[\atright, \solidarrow]
\mor{G}{M}{$\tau_{\Phi(p)}$}[\atleft, \solidarrow]
\mor{N}{M}{$\Phi$}[\atright, \solidarrow]
\enddc 
\end{gathered}
\EqTag{CDappD}
\end{equation}
	Since the maps $\mu_p:H\to N$ and $\tau_{\Phi(p)}:G \to M$ are imbeddings, 
	the function $\Phi$ restricts to a smooth map from the orbit of $H$ through $p$  to the orbit of $G$ though $\Phi(p)$. 
	This implies that $\phi:H \to G$ is  also smooth.

	To prove the theorem we will prove that  
	the map $\phi:H \to G$ defined in diagram \EqRef{CDappD} is a Lie group homomorphism with injective differential. We also 
	show that the definition of  $\phi:H \to G$ is independent of the choice of the point $p$.

	From the commutative diagram \EqRef{CDappD} we have $\Phi(h \cdot p) = \phi(h) \cdot \Phi(p)$
	which, on account of the freeness of the action of $G$, implies  that $\phi(e) = \epsilon$, where $e$ is the identity element of $H$
	and $\epsilon$ is the identity element of $G$.
	By taking differentials of the maps in  \EqRef{CDappD} and applying the chain rule we  
	therefore obtain the following commutative diagram
\begin{equation}
\begin{gathered}
\begindc{\commdiag}[3]
	\obj(0, 0)[H]{$T_{e}H$}
	\obj(30, 0)[G]{$T_{\epsilon}G$}
	\obj(0, -15)[N]{$T_pN$}
	\obj(30, -15)[M]{$T_{\Phi(p)}M \ .$}
	\mor{H}{G}{$\phi_{*}$}[\atleft, \solidarrow]
	\mor{H}{N}{$\mu_{p,*}$}[\atright, \solidarrow]
	\mor{G}{M}{$\tau_{\Phi(p), *}$}[\atleft, \solidarrow]
	\mor{N}{M}{$\Phi_{*}$}[\atright, \solidarrow]
\enddc 
\end{gathered}
\EqTag{CDappD2}
\end{equation}
	
	Recall that $ \lieh=T_{e}H$ and $\lieg=T_{\epsilon}G$ are the Lie algebras of $H$ and $G$ 
	defined by brackets of right invariant vector-fields and that the maps $\mu$ and $\tau$ 
	induce Lie algebra isomorphisms $\rho:\lieh \to \Gamma_H$  and $\lambda : \lieg \to \Gamma_G$ by  
	$\rho(Z_{e}) = \mu_{x*}( Z_{e}) $ and  $\lambda(W_{\epsilon}) = \mu_{y*} (W_{\epsilon})$ (see \EqRef{infgen}).
	Hence there is a uniquely defined Lie algebra monomorphism $\psi : \lieh \to \lieg$ such that  the following diagram
\begin{equation}
\begin{gathered}
\begindc{\commdiag}[3]
	\obj(0, 0)[H]{$\lieh$}
	\obj(30, 0)[G]{$\lieg$}
	\obj(0, -15)[N]{$\Gamma_H$}
	\obj(30, -15)[M]{$\Gamma_G$}
	\mor{H}{G}{$\psi$}[\atleft, \solidarrow]
	\mor{H}{N}{$\rho$}[\atright, \solidarrow]
	\mor{G}{M}{$\lambda$}[\atleft, \solidarrow]
	\mor{N}{M}{$\Phi_*$}[\atright, \solidarrow]
\enddc 
\end{gathered}
\EqTag{CDappD3}
\end{equation}
	commutes. The combination of the diagrams \EqRef{CDappD2} and \EqRef{CDappD3}  yields 
\begin{align*}
	\tau_{\Phi(p), *} (\phi_{*}(Z_e))  
&	= \Phi_{*}(\mu_{p, *}(Z_e))  =  \Phi_{*}(\rho((Z_e)_p) =  (\Phi_*(\rho(Z_e))(p)
\\
&	= \lambda(\psi(Z_e))(\Phi(p))  = \tau_{\Phi(p), *}(\psi(Z_e))
\end{align*}	
	and therefore, by the injectivity of  $\tau_{\Phi(p), *} $,   we have $\phi_{*}(Z_e)  = \psi(Z_e)$ and  therefore $\phi_{*} = \psi$.
	This proves that $\phi_{*} : T_eH \to T_\epsilon G$ is a Lie algebra monomorphism. 
	 It  then follows from  
	Exercise 9, page 134 of  \cite{warner:1983a} that $\phi:H \to G$ is a homomorphism (with injective differential). 
	Moreover, since the map $\psi$ does not depend upon the point $p$,  then   
	$\phi_{*}(Z_e)  = \psi(Z_e)$ does not  depend upon $p$.  By the uniqueness theorem for
	homomorphisms (Theorem 3.16, page 92  \cite{warner:1983a}), the homomorphism  $\phi$ does not depend on $p$. 
\end{proof}


\end{document}